\documentclass[a4wide,10pt]{article}
\usepackage[utf8]{inputenc}
\usepackage{flushend}
\usepackage{graphicx}
\usepackage[misc]{ifsym}
\usepackage{footmisc}
\usepackage{enumitem}
\usepackage{subfigure}
\usepackage{url}
\usepackage{multirow}
\usepackage[noadjust]{cite}
\usepackage[numbers,sort&compress]{natbib}
\usepackage{amsmath,amsthm}
\usepackage{amssymb,amsfonts}
\usepackage{booktabs}
\usepackage{color}
\usepackage{booktabs}
\usepackage{float}
\usepackage{fancyhdr}
\usepackage{xcolor,stfloats}
\usepackage{comment}
\usepackage{cuted}  
\usepackage{epstopdf}
\usepackage{multicol}
\usepackage{algorithmic}
\usepackage{algorithm}
\usepackage{bm}
\usepackage{footnote}
\newcommand\norm[1]{\left\lVert#1\right\rVert}
\title{Fast and Accurate\\ Surface Normal Integration \\on Non-Rectangular Domains}
\author{M. B\"ahr, M. Breu\ss{}, Y. Qu\'eau,\\ A.~S. Boroujerdi and J.-D. Durou}

\begin{document}

\maketitle

\begin{abstract}
The integration of surface normals for the purpose of computing the 
shape of a surface in 3D space is a classic problem in computer vision.
However, even nowadays it is still a challenging task to devise a method 
that combines the flexibility to work on non-trivial computational 
domains with high accuracy, robustness and computational efficiency.

By uniting a classic approach for surface normal integration 
with modern computational techniques we construct a solver that fulfils these requirements.
Building upon the Poisson integration model we propose to use an iterative Krylov subspace 
solver as a core step in tackling the task. While such a method can be very efficient, it 
may only show its full potential when combined with a suitable numerical preconditioning 
and a problem-specific initialisation. 

We perform a thorough numerical study in order to identify an appropriate preconditioner for 
our purpose. To address the issue of a suitable initialisation we propose to compute this initial 
state via a recently developed fast marching integrator. Detailed numerical experiments illuminate 
the benefits of this novel combination. In addition, we show on real-world photometric stereo datasets
that the developed numerical framework is flexible enough to tackle 
 modern computer vision 
applications.
\end{abstract}

\section{Introduction}

The integration of surface normals is a fundamental task in computer vision. 
Classic examples for processes where this technique is often applied are image edition~\cite{Perez2003},
shape from shading as analysed by Horn~\cite{Horn} and photometric stereo (PS), 
see the pioneering work of Woodham~\cite{Woodham1980}. 
Regarding PS, some modern applications include facial 
recognition~\cite{zafeiriou2013}, industrial product quality control~\cite{Smith2000}, 
object preservation in digital heritage~\cite{Esteban2008},
or new utilities for potential use in the area of the video (game) industry 
or robotics~\cite{Haque2014}, among many others. 

In this paper we illustrate surface 
normal integration at hand of the PS problem which serves as a role model for potential 
applications. The task of PS is to compute the 3D surface of an object 
from multiple images of the same scene under different illumination conditions.
The standard method of PS is divided into two stages. In a first step
the depth of an unknown surface in 3D space is computed in terms of a field of surface normals, 
or equivalently, a corresponding gradient field. In a subsequent step this needs to be integrated 
in order to obtain the depth of the surface. 

To handle the integration step, many different 
approaches and methods have been developed during the last decades. However, despite all these 
developments there is still the need for approaches that combine a high {\em accuracy} of the 
reconstruction with {\em robustness} against noise and outliers and a reasonable
computational {\em efficiency} for working with high-resolution cameras and 
corresponding imagery.

\paragraph{Computational Issues}
Let us briefly elaborate on the demands on an ideal integrator.
As discussed e.g.\ in~\cite{Harker} a practical issue is the robustness with respect
to noise and outliers. Since computer vision processes such as PS rely on 
simplified assumptions that often do not capture realistic illumination and surface reflectance, 
such artefacts may often arise when estimating surface normals at hand of real-world input images. 
Therefore, the determined depth gradient field is not noise-free and it may also contain outliers.
This may have a strong influence on the integration process.  

Secondly, objects of interest for 3D reconstruction are typically in the centre of a 
photographed scene. Therefore, they are only displayed on some part of an input image.
The sharp gradient usually representing the transition from foreground to background is a difficult 
feature for most surface normal integrators
and generally influences the shape estimation of the object of interest. Because of this it is desirable to consider only image segments that represent the object of interest 
and not the background. Although similar difficulties (sharp gradients) also arise for discontinuous surfaces that may appear at self-occlusions of an object~\cite{Durou2009}, we do not tackle this issue in the present article. Instead, we neglect self-occlusions and focus on reconstructing a \emph{smooth} surface on the possibly non-rectangular subset of the image domain representing the object of interest.
Related to this point, another important aspect is the computational 
time and cost saving that can be achieved by the simultaneous decreasing of the number of elements of the 
computational domain. With respect to reconstruction quality and efficiency, an ideal solver for 
surface normal integration should thus work on non-rectangular domains. 

Finally, to capture as much detail as possible of a 3D reconstruction, camera technology evolves 
and the resolution of images tends to increase continually. This means, the 
 integrator 
has to work accurately and fast for various sizes of input images, including large images
at least of size $1000\times 1000$ pixels. Consequently, the computational efficiency of a solver 
is a key requirement for many possible modern and future applications.

To summarise, one can identify the desirable properties of \emph{robustness} with respect to noise and outliers, 
ability to work on \emph{non-rectangular domains} and the \emph{efficiency} of the method meaning that
one aims for an accurate solution using reasonable computational resources (such as time and memory).

\paragraph{Related Work}
In order to solve the problem of surface normal integration many methods have been developed during 
the last decades that take into account the abovementioned issues in individual aspects. 

According to Klette and Schl\"uns~\cite{Klette1996}, the integration methods can be 
classified in two categories, local and global integration methods. 

The most basic local method, also 
referred to as direct line-integration scheme~\cite{Coleman1982,Wu,Robles-Kelly}, is based on the 
line integral technique and the fact that a closed path on a continuous surface should be zero. 
These methods are in general quite fast, however as by their local nature, the solution of the 
reconstruction depends on the integration path. 
Another more recent local approach is based on an eikonal-type equation for normal integration, 
which can be solved by applying the computationally quite efficient fast marching 
(FM) method~\cite{Ho,Galliani2012,BB15}.
However, a common disadvantage of all the local approaches is the sensitivity with 
regard to noise and discontinuities which eventually leads to error accumulation in the reconstruction. 

In order to minimise error accumulation it is preferable to adopt a global approach based on the 
calculus of variations. Horn and Brooks~\cite{Horn} proposed the classic and most natural variational 
method for the intended task by casting the corresponding functional in a format that 
results in a least-squares approximation. The necessary optimality condition represented by the 
Euler-Lagrange equation of the classic functional is given by the Poisson equation
which is an elliptic partial differential equation (PDE). 
This approach to the surface normal integration problem is often called \emph{Poisson integration}. 
The arising practical task amounts to solve the linear system of equations that corresponds to the discretised
Poisson equation. 

Direct methods for solving the latter system, as for instance Cholesky factorisation, can be fast, 
however this type 
of solvers may use substantial memory and appears to be rather impractical for images larger 
than $1000\times 1000$. Moreover, if based on matrix factorisation the factorisation itself is relatively 
expensive to compute. Generally, direct methods offer an extremely high accuracy of the result
but one must pay computationally a price for that.
In contrast, \emph{iterative methods} are naturally not tweaked for extremely high accuracy
but are very fast in computing approximate solutions. They
require less memory space and are thus 
inherently more attractive candidates for this application, however, they involve 
some non-trivial aspects (upon which we elaborate) which makes a tool that is not straightforward to use. 

An alternative approach to solve the least-squares functional has been introduced by 
Frankot and Chellappa~\cite{Frankot1988}. 
The main idea is to transform the problem to the frequency domain where a solution can be computed 
very fast in linear time through the Fast Fourier Transform (FFT), if 
periodic boundary conditions are assumed. 
The latter unfavourable condition can be resolved with the use of the Discrete Cosine Transform (DCT) as shown 
by Simchony {\em et al.}~\cite{Simchony1990}. However, these methods remain limited to rectangular domains. 
To apply these methods on non-rectangular 
 domains require introducing a zero-padding in the gradient field which may lead to an 
unnecessary bias in the solution. Some conceptually related basis function approaches include the use of 
wavelets~\cite{WeiKlette2001} and shapelets~\cite{Kovesi2005}. The method of Frankot and Chelappa was enhanced by 
Wei and Klette~\cite{WeiKlette2003} to improve the accuracy and robustness against noise. 
Another approach is proposed by Karacali and Snyder~\cite{Karacali} who make use of an additional adaptive 
smoothing for noise reduction. 

Among all the mentioned global techniques \emph{variational methods} offer a high robustness with 
respect to noise and outliers. Therefore, many extensions have been developed in modern 
works~\cite{Agrawal,Badri,Du,Reddy,Durou2009,Queau2015,Harker2015}.
Agrawal {\em et al.}~\cite{Agrawal} use anisotropic instead of isotropic weights for the gradients during 
the integration. The paper~\cite{Durou2009} shows a numerical study of several functionals, in 
particular non-quadratic and non-convex regularisations. To reduce the influence of outliers 
the $L_1$ norm has also become an important regularisation instrument, see~\cite{Du,Reddy}. 
In~\cite{Badri} the extension to $L_p$ minimisation with $0<p<1$ is presented.
Two other recent works are~\cite{Queau2015} where the use of alternative optimisation 
schemes is explored and~\cite{Harker2015} where the proposed formulation leads
to the task of solving a Sylvester equation. Nevertheless, the mentioned methods have 
some drawbacks. By the application of additional regularisations 
as e.g.\ in~\cite{Agrawal,Badri,Du,Reddy,Durou2009,Queau2015} the reconstruction of the depth is quite
time-consuming 
and
 the correct setting of parameters is more difficult,
whereas the approach of Harker and O'Leary in~\cite{Harker2015} is only efficient for 
rectangular domains $\Omega$. 

As a conclusion of the achievements of previous works, the problem of surface normal integration on non-rectangular domains
is not perfectly solved by now. The main challenge is still to find a balance between the quality 
and the computational time to generate the result. Thereby, one should take into account that
for achieving a high quality in 3D object reconstruction the resolution of images tends to increase 
continually and thus the computational efficiency is surely a key requirement for many 
potential current and future applications.

\paragraph{Our Contributions}
To bring the aspects of quality, robustness and computational efficiency into balance, 
we go back to the powerful classic approach of Horn and Brooks as the variational framework
has the benefit of high modeling flexibility. 
In detail, our contributions when extending this classic path of research are:

{\em (i)} Building upon a recent conference paper where we compared several Krylov subspace methods for surface normal integration~\cite{BQBD16}, we investigate the use of 
the preconditioned conjugate gradient (PCG) method 
for performing Poisson integration over non-trivial computational domains. While such methods constitute advanced yet standard methods in numerical computing~\cite{Meister2,Saad}, they still do not belong to the standard tools in image processing, computer vision and graphics. To be more precise, we propose to employ the conjugate gradient (CG) scheme as the iterative solver and we explore modern variations of the incomplete Cholesky (IC) decomposition for preconditioning. The thorough numerical investigation performed here represents a significant extension of our conference paper. 

{\em (ii)} For computing a good initialisation for the PCG solver we propose to employ 
a recent FM integrator~\cite{Galliani2012} already mentioned above. 
The main advantages of the FM integrator are its flexibility for 
use with non-trivial domains coupled with a low computational complexity.
While we proposed this means of initialisation already in the work~\cite{BQBD16} which we extend here, 
let us note that due
to our other numerical extensions the conclusion we can draw in this paper is a much sharper one.

{\em (iii)} We prove experimentally that our resulting, combined novel method 
unites the advantages of {\em flexibility} and {\em robustness} of variational 
methods with {\em low computational times} and {\em low memory requirements}. 
{\em (iv)} We propose a simple yet effective modification of gradient fields containing severe outliers for use with Poisson integration methods.

Let us note that the abovementioned building blocks of our method represent a pragmatic choice among the current tools  
in numerical computing. Moreover, as demonstrated by our new integration model that is specifically designed for 
tackling data with outliers, our numerical proceeding is readily adapted to possible new Poisson-based integration models.
This together with the well-engineered algorithm for our application, 
i.e.\ the FM initialisation and fine-tuned algorithmic parameters, makes the method we 
propose a unique, efficient and flexible procedure.

\section{Surface Normal Integration}\label{FM-Krylov}

The mathematical set-up of surface normal integration (SNI) can be described as follows. We assume that for a domain $\Omega$ a normal field 
$\mathbf{n}:=\mathbf{n}(x,y)=\left[n_1(x,y),n_2(x,y),n_3(x,y)\right]^\top$ is given for each grid point $(x,y)\in\Omega$. 
The task is to recover a surface $S$, which can be represented as a depth map $v(x,y)$ over $(x,y)\in\Omega$, 
such that $\mathbf{n}$ is the normal field of~$v$. 
Assuming orthographic projection\footnote{The perspective integration problem can be formulated in a similar way, provided that the change of variable $v = \log v$ is introduced~\cite{Durou2}.}, a
 normal field $\mathbf{n}$ of a surface at $(x,y,v(x,y)) \in \mathbb{R}^3$ can be written as
\begin{equation}
 \mathbf{n}(x,y):=\frac{\left[-v_x,-v_y,1\right]^\top}{\sqrt{\Vert\nabla v\Vert^2+1}}\label{n}
\end{equation}
with $v_x:=\frac{\partial v}{\partial x}$, $v_y:=\frac{\partial v}{\partial y}$ and $\nabla v := \left[v_x,v_y\right]^\top$. 
Moreover, the components of $\mathbf{n}$ are given by
partial derivatives of $v$ 
with
\begin{equation}
\left(
v_x, v_y
\right)
=
\left(
-\frac{n_1}{n_3},
-\frac{n_2}{n_3}
\right)
=
\left(
p,q
\right)
\label{nComp}
\end{equation}
where we think of $p$ and $q$ as given data.

In this section, we present the building blocks of our new algorithm in two steps,
first the fast marching integrator and afterwards the iterative Poisson solver 
relying on the conjugate gradient method supplemented by (modified) incomplete 
Cholesky preconditioning. In the presentation of the Poisson integration
we also demonstrate by a novel adaptation for handling data with outliers
 the flexibility of the arising discrete computational model.

Let us note that the detailed description of the fast marching integrator can be found in 
the conference papers~\cite{Galliani2012,BB15}, and the presentation of the 
components of the CG scheme can be found in literature on Krylov subspace solvers,
see e.g.~\cite{Saad}. We still recall here the algorithms in some detail
because there are important parameters that need to be set
and some choices to make: since the efficiency of integrators depends largely on such practical implementation details, 
we think that our explanations generate an additional value beyond a plain description of the methods. 

While our discretisation of the Poisson equation is a standard one, 
we deal with non-trivial boundary conditions
in our application which necessitates a thorough description of the latter. The construction of our
non-standard numerical boundary conditions, which is often overlooked in the literature,  is another technical contribution to the field.

\subsection{The Fast Marching Integrator}

We recall for the convenience of the reader some of the
relevant developments from~\cite{Ho,Galliani2012,BB15}. In these works it has been shown that it is possible to tackle the problem
of surface normal integration via the following PDE-based model in $w = v + \lambda f$:
\begin{equation}
\Vert \nabla w\Vert=\sqrt{(p+\lambda f_x)^2+(q+\lambda f_y)^2} 
\label{FMpde}
\end{equation} 
where $\lambda >0$ and $f: \mathbb{R} \to \mathbb{R}$ are user-defined. Let us note that the FM integrator requires parameter $\lambda$ to be tuned, yet is not a crucial choice as any large number $\lambda \gg 0$ will work~\cite{Galliani2012}.\footnote{In our experiments, we used the value $\lambda = 10^5$.}

Using PDE~\eqref{FMpde} we do not compute the depth function $v$ directly, but instead we solve in a first step for $w$. 
Let us note that this intermediate step is necessary for the successful application of the FM method, in order to avoid local minima and ensure that any initial point can be considered~\cite{Ho}.

It turns out that a natural candidate for $f$ is the squared Euclidean distance function with the minimum 
in the centre of the domain $(x_0,y_0)=(0,0)$, i.e. 
\begin{equation}
f:=f(x,y)=x^2+y^2
\end{equation}
Note that other choices for $f$ are possible~\cite{BB15}.
Then in \eqref{FMpde} on the right hand side, $(p,q)$ are the given data and $\nabla f$ is known.
As boundary condition we may employ $w(0,0)=0$. After the computation of $w$ we easily compute the sought 
depth map $v$ via $v=w-\lambda f$.

\paragraph{Numerical Upwinding}

A crucial point of the FM integrator is the correct discretisation
of the derivatives of $f$ in (\ref{FMpde}). In order to obtain a stable method, 
an upwind discretisation for the partial derivatives of $f$ is required, reading as
\begin{equation}
f_x :=
\Big[\max\Big(\frac{f_{i,j}-f_{i-1,j}}{\Delta x},\frac{f_{i,j}-f_{i+1,j}}{\Delta x},0\Big)\Big]^2
\label{rt-discrete}
\end{equation}
and analogously for $f_y$, where we make use of grid widths $\Delta x$ and $\Delta y$. 

Making use of the same discretisation for the components of $\nabla w$,
one obtains a quadratic equation that needs 
to be solved for every pixel except at the initial pixel $(0,0)$ where some depth value is prescribed.

Let us note that the initial point can be chosen in practice anywhere, i.e.\ there is no restriction 
to $(0,0)$.

\paragraph{Non-convex Domains}

If the above method is used without modification over 
non-convex domains, the FM 
integrator eventually fails to reconstruct the solution.
The reason is that the original, unmodified squared Euclidean distance
does not yield a meaningful distance from 
the starting point to pixels without a connection by a direct line in the integration domain. 
In other words, the unmodified scheme just works over convex domains. 

To overcome the problem, a suitable distance which calculates the shortest path from 
the starting point to every point on the computational domain is necessary. 

To this end, the use of a geodesic distance function $d$ is advocated~\cite{Galliani2012}. 
The proceeding is as follows, relying on similar ideas as used e.g.\ in~\cite{kimmelsethian01}
for path planning. In a first step we solve an eikonal equation 
$\left\Vert \nabla d \right\Vert=1$ over all the points of the domain 
with $d:=0$ at the chosen start point. This can of course be done again with the FM method.
Then, in a second step we are able to compute the depth map $v$. 
Therefore, we use equation (\ref{FMpde}) for $w$, with the squared geodesic distance function $d$ 
instead of $f$ and using (\ref{rt-discrete}). Afterwards we 
recover $v$ via $v=w-\lambda d$.

\paragraph{The Fast Marching Algorithm}

The idea of FM goes back to the works~\cite{Tsitsiklis,Sethian2,Helmsen}. 
For a com\-pre\-hensive introduction see e.g.~\cite{Sethian}. The benefit of FM 
is its relatively low complexity of $\mathcal{O}(n\log{} n)$ where $n$ is the number of points in the 
computational domain.\footnote{When using the \emph{untidy priority queue} structure~\cite{Yatziv2006} the 
complexity may even be lowered to $O(n)$.}

Let us briefly describe the FM strategy. The principle behind FM is that information advances from smaller 
values $w$ to larger values $w$, meeting thereby each point of the computational domain just once. 
To this end, one may employ three disjoint sets of nodes as discussed in detail in~\cite{Sethian,CCF14}: 
$(s1)$ accepted nodes, $(s2)$ trial nodes and $(s3)$ far nodes. The values $w_{i,j}$ of set $(s1)$
are considered known and will not be changed. 
A member $w_{i,j}$ in set $(s2)$ is always at a neighbour of an accepted node. This is the set 
where the computation actually takes place and the values of $w_{i,j}$ can still change. 
In set $(s3)$ are the nodes $w_{i,j}$ where 
an approximate solution has not yet been computed as these are not in a neighbourhood of a member of $(s1)$. 

The FM algorithm can then be described by the following procedure until all nodes are accepted:
\begin{enumerate}
 \item[(a)] Find the grid point $A$ in $(s2)$ with the smallest value and change it to $(s1)$.
 \item[(b)] Place all neighbours of $A$ into $(s2)$ if they are not there already and compute the arrival time 
for all of them, if they are not already in $(s1)$.
 \item[(c)] If the set $(s2)$ is not empty, return to (a).
\end{enumerate}
Let us finally note that for initialisation, one may take the 
node at $(0,0)$ which bears the boundary condition of the PDE \eqref{FMpde} and put it into set $(s1)$.

An efficient implementation amounts to store the nodes in $(s2)$ in a heap data structure, 
so the smallest element in step $(a)$ can be chosen as fast as possible.

\subsection{Poisson Integration}
\label{krylov}

The first part of this section is dedicated to \emph{modeling}.
We first briefly review the classic variational approach to the Poisson integration problem, 
cf.~\cite{Durou2,Harker2015,Horn,Queau2015,Simchony1990}. 
The handling of extremely noisy data eventually 
motivates modifications with regard to the underlying energy functional 
(\ref{energy}), compare e.g.~\cite{Queau2015}. We demonstrate 
by proposing a new model dealing with outliers that the Poisson integration 
framework is flexible enough to deal with such modern models.

The second part is devoted to the \emph{numerics}. We propose a dedicated, to some
part non-standard discretisation for our application.

\paragraph{Classic Poisson Integration (PI) Model}
In order to recover the surface it is common to minimise the least squares 
error between the input and the gradient field of $v$ via minimising
\begin{align}
 J(v)&=\iint_\Omega \Vert \nabla v-\mathbf{g}\Vert^2 \,dx \,dy\nonumber\\
 &=\iint_\Omega (v_x-p)^2+(v_y-q)^2 \,dx\,dy\label{energy}
\end{align}
where we denote $\mathbf{g} = \left[p,q\right]^\top$.

A minimiser $v$ of $(\ref{energy})$ must satisfy the associated Euler-Lagrange equation 
which is equivalent to the following {\em Poisson equation}
\begin{equation}
 \Delta v= \mathrm{div}(p,q)=p_x+q_y\label{poisson}
\end{equation}
that is usually complemented by (natural) Neumann boundary conditions \newline $(\nabla v-\mathbf{g})\cdot{\bm \mu}=0$,
where the vector ${\bm \mu}$  
is normal to $\partial\Omega$. In this case, the uniqueness of the solution
is guaranteed, apart from an additional constant. 
Thus, one recovers the shape but not an account of absolute depth (as by FM integration).

\paragraph{A Modified PDE for Normal Fields with Outliers}

We will now demonstrate exemplarily that the PI framework is flexible enough to allow
to deal also with gradient fields featuring strong outliers. To this end, we propose a simple, yet effective way to modify the PI model in order to limit the influence of outliers.
Other variations for different applications, e.g.\ self-occlusions~\cite{Durou2009,Queau2015}, are of course also possible.

Let us briefly recall that the classic model \eqref{energy} which leads to the Poisson equation
\eqref{poisson} is based on a simple least-squares approach.
On locations $(x,y)$ corresponding to outliers, the values $p(x,y)$ and $q(x,y)$ are not reliable, and 
one would prefer to limit the influence of these corrupted data.

Therefore, we modify the Poisson equation \eqref{poisson}
by introducing a space-dependent fidelity term $\nu := \nu (x,y)$ via
\begin{equation}
  \Delta v = \nabla \cdot \left( \frac{1}{1+\nu} \left[p,q\right]^\top \right)
\label{poisson-modified-equation}
\end{equation}
Let us note that a similar strategy, namely to introduce modeling improvements in a PDE 
that is originally the Euler-Lagrange equation of an energy functional instead of modifying the latter, 
is occasionally employed in computer vision, see e.g.~\cite{ZBVBWRS08}. However, we do not tinker here 
with the core of PDE i.e.\ the Laplace operator $\Delta$, we merely install a preprocessing
by modifying the right hand side of the Poisson equation.

The key to an effective preprocessing is of course to consider the role 
of $\nu$ so that it smoothes the surface only at locations where the input gradient is not reliable. 
Thus, we have to seek a function $\nu(x,y)$ which 
would be close to zero if the input gradient is reliable, and takes high values when it is not.

Let us remark that the \emph{integrability} term
\begin{equation}
  \mathcal{I}(x,y) := p_y - q_x = \nabla \cdot \left[ -q,p \right]^\top 
\end{equation}
should vanish if the surface is $\mathcal{C}^2$-smooth. This argument was used in~\cite{Queau2015} to 
suggest an integrability-based weighted least-squares functional able to recover discontinuity jumps, 
which generally correspond to a high integrability absolute value. 

Since integrability not only indicates the location of discontinuities, but also that of the outliers, 
we suggest to use this integrability term in the context of finding a smooth surface explaining a 
corrupted gradient. For this purpose, we introduce the following choice for our regularisation parameter:
\begin{equation}
  \nu(x,y) = \exp\left(\mathcal{I}(x,y)^2\right)-1
\label{lambda-varies}
\end{equation}
which holds the desired properties {\em (i)} to vanish when integrability is low (reliable gradients), and 
{\em (ii)} to take a high value when integrability is high (outliers). 

Putting it altogether, our new model amounts to the resolution of the following equation:
\begin{equation}
\Delta v = \nabla \cdot \left[ \overbrace{\frac{p}{\exp\left\{\left(p_y-q_x\right)^2\right\}}}^{\overline{p}},\, 
\overbrace{\frac{q}{\exp\left\{\left(p_y-q_x\right)^2\right\}}}^{\overline{q}} \right]^\top
\label{poisson-outlier}
\end{equation}
which is another Poisson equation, where the right hand side can be computed \emph{a priori} from the input gradient. 

Let us emphasise that \emph{all} methods dedicated to SNI by resolution of such a Poisson equation
are straightforward to adapt: it is enough to 
replace $(p,q)$ by $(\overline{p},\overline{q})$. The algorithmic complexity for all
of such approaches remains exactly the same. The practical validity of this simple new model
and the benefit of better numerics for it are demonstrated in Section \ref{PCG}.

In the main part of our paper, we will stick for simplicity of presentation to the classic model
\eqref{poisson} and come back to the proposed modification in Section~\ref{PCG}.

\paragraph{Discretisation of the Poisson Equation}

A useful standard numerical approach to solve the Poisson PDE as in \eqref{poisson}
or \eqref{poisson-outlier} makes use of 
finite differences. Often $\mathrm{div}(p,q)$ and $\Delta v=v_{xx}+v_{yy}$ are approximated
by central differences. For simplicity we suppose that the grid 
size is $\Delta x=\Delta y=1$ as is common practice in image processing. 
Then, a suitable discrete version of the Laplacian is given in stencil notation by
\begin{equation}
 \Delta v(x_i,y_j)\approx 
 \begin{tabular}{|c|c|c|}
 \hline
 & 1 & \\
\hline
1 & -4 & 1\\
\hline
 & 1 & \\
\hline
\end{tabular}
\cdot v_{i,j}
\label{discLapl}
\end{equation}
and for the divergence by
\begin{equation}
 \mathrm{div}(p_{i,j},q_{i,j})\approx \frac{1}{2}
 \begin{tabular}{|c|c|c|}
\hline
-1 & 0 & 1\\
\hline
\end{tabular}
\cdot p_{i,j}+\frac{1}{2}
 \begin{tabular}{|c|}
 \hline
 1\\
\hline
0\\
\hline
 -1 \\
 \hline
\end{tabular}
\cdot q_{i,j}
\label{discDiv}
\end{equation}
with the measured gradient $\mathbf{g}=\left[p,q\right]^\top$. Making use of (\ref{discLapl}) and (\ref{discDiv}) 
for the discretisation of (\ref{poisson}) leads to
\begin{align}
 -4v_{i,j}+&(v_{i+1,j}+v_{i-1,j}+v_{i,j+1}+v_{i,j-1})=\nonumber\\
 &\frac{p_{i+1,j}-p_{i-1,j}+q_{i,j+1}-q_{i,j-1}}{2}
\label{system-1}
\end{align}
which corresponds to a linear system $A\mathbf{x}=\mathbf{b}$,
where the vectors $\mathbf{x}$ and $\mathbf{b}$ are obtained by stacking the unknown values $v_{i,j}$ and 
the given data, respectively. The matrix $A$ contains the coefficients arising by the discretisation
of the Laplace operator $\Delta$. 

We employ in all experiments here the above discretisation, as it is very simple and still gives results of 
high quality. While other discretisations as e.g.\ of higher order
are of course possible~\cite{Harker2015}, let us note that this implies that one should change the parameter settings 
of the method as we propose it. One would also have to adapt the dedicated numerical boundary conditions.

\paragraph{Non-standard Numerical Boundary Conditions}
At this point it should be noted that the stencils in (\ref{discLapl}), (\ref{discDiv}) and 
the subsequent equation (\ref{system-1}) are only valid for inner points of the computational domain. 
Indeed, when pixel $(i,j)$ is located near the border of $\Omega$, some of the four neighbour 
values $\{v_{i+1,j},v_{i-1,j},\linebreak v_{i,j+1},v_{i,j-1}\}$ 
in~\eqref{system-1} refer to the depth outside $\Omega$. The same holds for the data values 
$\{ p_{i+1,j},p_{i-1,j},q_{i,j+1},q_{i,j-1} \}$: some of these values are unknown when $(i,j)$ is near the border. 
To handle this, a numerical boundary condition must be invoked. 

Using empirical Dirichlet (e.g., using discrete sine transform~\cite{Simchony1990}) or homogeneous Neumann boundary conditions~\cite{Agrawal} may result in biased 
3D reconstructions near the border. The so-called ``natural" condition  
$(\nabla v-\mathbf{g})\cdot{\bm \mu}=0$~\cite{Horn} is preferred, because it is the only one which is justified. 

Let us emphasise that it is not a trivial task to define suitable 
boundary conditions for $\{ p_{i+1,j},p_{i-1,j},q_{i,j+1},q_{i,j-1} \}$.
As we opt for a common strategy for discretising values of $p,q,v$, we employ the following
non-standard procedure which has turned out to be preferable in experimental evaluations.
Whenever $p,q,v$ values outside $\Omega$ are involved in~\eqref{system-1}, we discretise 
this boundary condition using the mean of forward and backward first-order finite differences. 
This allows us to express the values outside $\Omega$ in terms of values inside $\Omega$. To clarify 
this aspect, let us distinguish the boundaries according to the number of 
missing neighbours. 

\medskip

{\em When only one neighbour is missing} \quad There are four types of boundary pixels having exactly 
one of the four neighbours outside $\Omega$ (respectively lower, upper, right and left borders). 
Let us first consider the case of a ``lower boundary", i.e.\ a pixel $(i,j) \in \Omega$ such 
that $(i-1,j),(i+1,j),(i,j+1) \in \Omega^3$ but $(i,j-1)   \notin \Omega$. Then, Eq.~\eqref{system-1} 
involves the undefined quantities $v_{i,j-1}$ and $q_{i,j-1}$. Yet, discretisation of the natural 
boundary condition, at pixel $(i,j-1)$, by forward differences, provides the following equation:
\begin{equation}
  v_{i,j}-v_{i,j-1} = q_{i,j-1}
  \label{eq:BC1}
\end{equation}
On the other hand the natural boundary condition can be also discretised at pixel $(i,j)$ by 
backward differences, leading to:
\begin{equation}
  v_{i,j} - v_{i,j-1} = q_{i,j}
  \label{eq:BC2}
\end{equation}
Taking the mean of the forward~\eqref{eq:BC1} and backward~\eqref{eq:BC2} discretisations, we obtain:
\begin{equation}
  v_{i,j}-v_{i,j-1} =  \frac{q_{i,j-1}+q_{i,j}}{2}
  \label{eq:BC3}
\end{equation}
Now, plugging~\eqref{eq:BC3} into~\eqref{system-1}, the undefined quantities actually vanish, and one obtains:
\begin{align}
 -3v_{i,j}+&(v_{i+1,j}+v_{i-1,j}+v_{i,j+1})=\nonumber\\
 &\frac{p_{i+1,j}-p_{i-1,j}+q_{i,j+1}+q_{i,j}}{2}
\label{eq:system-2}
\end{align}
In other words, the stencil for the Laplacian is replaced by
\begin{displaymath}
 \Delta v(x_i,y_{j})\approx 
 \begin{tabular}{|c|c|c|}
 \hline
  & 1 & \\
\hline
 1 & -3 & 1\\
\hline
   &    & \\
\hline
\end{tabular}
\cdot v_{i,{j}}
\label{eq:stencil_lap_1}
\end{displaymath}
and that for the divergence by
\begin{equation}
 \mathrm{div}(p_{i,j},q_{i,j})\approx \frac{1}{2}
 \begin{tabular}{|c|c|c|}
\hline
-1 & 0 & 1\\
\hline
\end{tabular}
\cdot p_{i,j}+\frac{1}{2}
 \begin{tabular}{|c|}
 \hline
 1\\
\hline
1\\
\hline
0 \\
 \hline
\end{tabular}
\cdot q_{i,j}
\label{eq:stencil_div_1}
\end{equation}
The corresponding stencils for upper, left and right borders are obtained by straightforward 
adaptations of this procedure. 

\medskip

{\em When two neighbours are missing} \quad
Boundary pixels having exactly two neighbours outside $\Omega$ are either ``corners" 
(e.g., $(i,j-1)$ and $(i+1,j)$ inside $\Omega$ but $(i-1,j)$ and $(i,j+1)$ outside $\Omega$) 
or ``lines" (e.g. $(i-1,j)$ and $(i+1,j)$ inside $ \Omega$, but $(i,j-1)$ and $(i,j+1)$ outside $\Omega$). 
For ``lines", the natural boundary condition must be discretised four times 
(both forward and backward, on the two locations of missing data). 
Applying a similar rationale as in the previous case, we obtain the following stencils for ``vertical" lines:
  \begin{displaymath}
 \Delta v(x_i,y_{j})\approx 
 \begin{tabular}{|c|c|c|}
 \hline
  & 1 & \\
\hline
  & -2 & \\
\hline
  & 1 & \\
\hline
\end{tabular}
\cdot v_{i,{j}}
\end{displaymath}
and
\begin{equation}
 \mathrm{div}(p_{i,j},q_{i,j})\approx \frac{1}{2}
 \begin{tabular}{|c|c|c|}
\hline
0 & 0 & 0\\
\hline
\end{tabular}
\cdot p_{i,j}+\frac{1}{2}
 \begin{tabular}{|c|}
 \hline
 1\\
\hline
0\\
\hline
-1 \\
 \hline
\end{tabular}
\cdot q_{i,j}
\end{equation}
A straightforward adaptation provides the stencils for the ``horizontal" lines. 
Applying the same procedure for corners, we obtain, for instance for the ``top-left" corner:
  \begin{displaymath}
 \Delta v(x_i,y_{j})\approx 
 \begin{tabular}{|c|c|c|}
 \hline
  &  & \\
\hline
  & -2 & 1 \\
\hline
  & 1 & \\
\hline
\end{tabular}
\cdot v_{i,{j}}
\end{displaymath}
and
\begin{equation}
 \mathrm{div}(p_{i,j},q_{i,j})\approx \frac{1}{2}
 \begin{tabular}{|c|c|c|}
\hline
0 & 1 & 1\\
\hline
\end{tabular}
\cdot p_{i,j}+\frac{1}{2}
 \begin{tabular}{|c|}
 \hline
 0\\
\hline
-1\\
\hline
-1 \\
 \hline
\end{tabular}
\cdot q_{i,j}
\end{equation}
Again, it is straightforward to find the other three discretisations for the other corner types. 

\medskip

{\em When three neighbours are missing} \quad
In this last case, we discretise the boundary condition six times 
(forward and backward, for each missing neighbour). Most quantities actually vanish. For instance, 
for the case where only the right neighbour $(i+1,j)$ is inside $\Omega$, we obtain the following stencils:
\begin{displaymath}
 \Delta v(x_i,y_{j})\approx 
 \begin{tabular}{|c|c|c|}
 \hline
  &  & \\
\hline
  & -1 & 1 \\
\hline
  &  & \\
\hline
\end{tabular}
\cdot v_{i,{j}}
\end{displaymath}
and
\begin{equation}
 \mathrm{div}(p_{i,j},q_{i,j})\approx \frac{1}{2}
 \begin{tabular}{|c|c|c|}
\hline
0 & 1 & 1\\
\hline
\end{tabular}
\cdot p_{i,j}+\frac{1}{2}
 \begin{tabular}{|c|}
 \hline
 0\\
\hline
0\\
\hline
0 \\
 \hline
\end{tabular}
\cdot q_{i,j}
\end{equation}
In the end, we obtain explicit stencils for all the fourteen types of boundary pixels. 
Let us emphasise that, apart from $4$-connectivity, we introduced no assumption on the shape of $\Omega$.

\paragraph{Summarising the Discretisation}
The discretisation procedure defines a sparse linear system of equations $A\mathbf{x}=\mathbf{b}$. 
Incorporating Neumann boundary conditions the matrix $A$ is 
symmetric, positive semidefinite, diagonal dominant and 
its null space contains the vector $e:=[1,\dots,1]^T$.
In other words, A is a rank-$1$ deficient 
and singular matrix.

\subsection{Iterative Krylov Subspace Methods}
\label{sec:krylovsubspaces}

As indicated, in consequence of the enormous computer memory costs the application 
of a direct solver to deal with the above linear system appears to be impractical
for large images. Therefore, we propose an iterative solver to handle this problem. 

A modern class of iterative solvers 
designed for use with large sparse linear systems is the 
class of {\em Krylov subspace solvers}, for a detailed exposition see e.g.~\cite{Meister14,Saad}.
The main idea behind the Krylov approach is to search for an approximate solution of $A\mathbf{x}=\mathbf{b}$, 
with $A\in\mathbb{R}^{n\times n}$ a large regular sparse matrix and $\mathbf{b}\in\mathbb{R}^{n}$, in a suitable 
low-dimensional (affine) subspace of $\mathbb{R}^{n}$ that is constructed iteratively. 

This construction is in general not directly visible in the formulation
of a Krylov subspace method, as these are often described in terms of a reformulation of
solving $A\mathbf{x}=\mathbf{b}$ as an optimisation task.
An important example for this is given by the classic conjugate gradient (CG)
method of Hestenes and Stiefel~\cite{Hestenes1952} which is still an adequate iterative solver for
problems involving sparse symmetric matrices as in (\ref{system-1})\footnote{While in general also positive definiteness is required, this point is more delicate. We comment later on 
the applicability in our case.}.

\paragraph{About the Conjugate Gradient Method}
As it is of special importance for this work, let us briefly recall some properties of
the CG method; a more technical, complete exposition can be found in many textbooks on 
numerical computing, see e.g.~\cite{Meister14,Meurant1999,Meurant06,Saad}.

Note that a useful implementation of CG is given in the Matlab 
package. However, some knowledge about the technique is useful in order to understand some 
algorithmic properties. Moreover, it is a crucial point in the effective application of the
CG method to be aware of its critical parameters. In our exposition we aim to make 
clear the corresponding points.

The CG method requires a symmetric and positive definite 
matrix $A\in\mathbb{R}^{n\times n}$. In its construction it combines the gradient descent method with the method of 
conjugate directions. It can be derived from making use of the fact, that 
for such a matrix the solution of $A\mathbf{x}=\mathbf{b}$ is exactly the minimum of the function
\begin{equation} 
F(\mathbf{x})=
\frac{1}{2}\langle \mathbf{x},A\mathbf{x}\rangle_2-\langle \mathbf{b},\mathbf{x}\rangle_2\label{CG-F}
\end{equation}
since
\begin{equation} 
\nabla F(\mathbf{x})=
\mathbf{0}
\quad
\Leftrightarrow
\quad
A\mathbf{x}=\mathbf{b}
\label{einschub-1}
\end{equation}
Thereby, $\langle \cdot,\cdot \rangle_2$ means the Euclidean 
scalar product.

Let us now denote the $k$-th Krylov subspace by $\mathcal{K}_k$. Then,  
$\mathcal{K}_k:= \mathcal{K}_k (A,\mathbf{r}_0)$ is a subspace of $\mathbb{R}^n$ defined as
\begin{equation} 
\mathcal{K}_k
:=
\mathsf{span}
\left(
\mathbf{r}_0,
A\mathbf{r}_0,
A^2\mathbf{r}_0,
\ldots
A^{k-1}\mathbf{r}_0
\right)
\label{einschub-2}
\end{equation}
This means $\mathcal{K}_k$ is generated from an initial residual vector $\mathbf{r}_0= b-A\mathbf{x}_0$ 
by successive multiplications with the system matrix $A$. 

Let us briefly highlight some important theoretical considerations.
The nature of an iterative Krylov subspace method is, that the computed approximate solution 
$\mathbf{x}_k$ is in $\mathbf{x}_0 + \mathcal{K}_k (A,\mathbf{r}_0)$,
i.e.\ it is determined by the $k$-th Krylov subspace. Thereby, the index $k$ 
is also the $k$-th iteration of the iterative scheme. 

For the CG method, 
one can show that the approximate solutions $\mathbf{x}_k$
are optimal in the sense that they minimise the so-called energy norm of the error vector. This means,
if $\mathbf{x}^\star$ is a solution of the system $A \mathbf{x} = \mathbf{b}$, that $\mathbf{x}_k$ minimises $\norm{\mathbf{x}^{\star} - \mathbf{x}_k}_A$
for the $A$-norm $\norm{\mathbf{y}}_A := \sqrt{\mathbf{y}^\top A \mathbf{y}}$. Note again that $\mathbf{x}_k$ is restricted to
the search in the $k$-th Krylov subspace. In other words, the CG method gives in the $k$-th iteration the
best solution available in the generated subspace. Since the dimension of the Krylov subspace is
increased in each step of the iteration, theoretical convergence is achieved at latest after the 
$n$-th step of the method if the sought solution is in $\mathbb{R}^n$. 

\paragraph{Practical Issues}
An useful observation on the Krylov subspace methods is, that they can obviously benefit a lot from
a good educated guess of the solution which could be used as the initial iterate $\mathbf{x}_0$.
Therefore, we consider $\mathbf{x}_0$ as an important open parameter of the method that can be addressed
in a proper way.

Moreover, an iterative method also requires the user to set a tolerance defining the usual \emph{stopping criterion}: if the norm of the relative 
residual is below the tolerance, the algorithm stops. 

However, there is \emph{a priori} no means to say in which regime the tolerance has to be chosen.
This is one of the issues that make a reliable and efficient application of the method not so straightforward.
It will be one of the aims of our experiments to determine for our application a reasonable tolerance.

While our presentation of the CG method relates to ideal theoretical properties, in practice numerical rounding errors
appear and one may suffer from severe convergence problems for very large systems. Thus, a \emph{preconditioning}
is recommended to enforce all the beneficial properties of the algorithm, along with fast convergence. 
However, as it turns out it will require a thorough study to identify the most useful parameters in the 
preconditioning method.

Let us note that the CG method is applicable even though our matrix $A$ is just positive \emph{semidefinite}. 
The positive definiteness is useful for avoiding division by zero within the CG algorithm. 
If $A$ is positive semidefinite, theoretically it may happen that one needs to restart the 
scheme using a different initialisation. In practice this situation rarely occurs.

\paragraph{Preconditioning}

The basic idea of preconditioning is to multiply
the originally arising system $A\mathbf{x}=\mathbf{b}$ from the left with a matrix $P$ such 
that $P$ approximates $A^{-1}$. The modified system $PA\mathbf{x}=P\mathbf{b}$ is in general better conditioned and 
much more efficient to solve. For sparse matrices $A$, typically preconditioners are 
defined over the same sparse structure of entries of $A$.

Dealing with symmetric matrices such as in our case, the incomplete Cholesky 
(IC) decomposition~\cite{Golub96} is often used for constructing 
a common and very efficient preconditioner for the CG method~\cite{Meijerink1977,Kershaw1978,Benzi2002}. 
As a consequence of~\cite{BQBD16} we study here the application of the IC preconditioner 
and its modified version MIC.

Let us briefly describe the underlying ideas. The complete decomposition of $A$ will be given by $A=LL^T+F$.
If the lower triangular matrix $L$ is allowed to have non-zero entries anywhere in the lower matrix, then
$F$ is the zero matrix and the decomposition is the standard Cholesky decomposition.
However, in the context of sparse systems only the structure of entries in $A$ is used in defining $L$, so 
that the factorisation will be incomplete. Thus, in our case the lower triangular matrix $L$ keeps the same 
non-zero pattern as that of the lower triangular part of $A$. The general form of the preconditioning 
then amounts to the transformation from $A\mathbf{x}=\mathbf{b}$ to $A^p\mathbf{x}^p=\mathbf{b}^p$ with
\begin{equation}
A^p=L^{-1}AL^{-T}, \, 
x^p=L^{-T}x
\quad
\text{and}
\quad
b^p=L^{-1}b
\label{precon-1}
\end{equation}

\paragraph{Practical Issues}
Let us identify another important computational parameter.
The mentioned approach to take as the sparsity pattern of $L$ the existing pattern in $A$
is often called IC($0$). If one extends the sparsity pattern of $L$ by additional non-zero elements 
(usually in the vicinity of existing entries) then the closeness between the product $LL^T$ and $A$ 
is potentially improved. 
This proceeding is often denoted as numerical fill-in strategy IC($\tau$), called ''drop tolerance``, where the parameter $\tau>0$ 
describes a dropping criterion, cf.~\cite{Saad}. The approach can be described as follows: new fill-ins are accepted only if the elements are greater than a local drop tolerance $\tau$.
It turns out that the corresponding PCG method is applicable for positive semidefinite matrices, 
see for example~\cite{Kaasschieter1988,Tang2007}.

When dealing with a discretised elliptic PDE as in~\eqref{poisson} or~\eqref{poisson-modified-equation}, 
the \emph{modified IC (MIC)} 
factorisation can lead to an even better preconditioner, for an 
overview on MIC see~\cite{Benzi2002, Meurant1999}. 
The idea behind the modification is to force the preconditioner to have the same 
row sums as the original matrix $A$. 
This can be accomplished by adding dropped fill-ins to the diagonal.
The latter is known as MIC($0$) and can be combined with the abovementioned drop tolerance strategy to MIC($\tau$).
Let us mention, that MIC can lead to possible pivot breakdowns. This problem can be circumvented by a global diagonal shift applied to $A$ 
prior to determining  the incomplete factorization, see \cite{Manteuffel}. Therefore, the factorization\footnote{We denote the combined methods of MIC($\tau$) and the shifted incomplete Cholesky version as 
MIC($\tau,\alpha$).} of $\tilde{A}=A+\alpha \operatorname{diag}(A)$ is performed, where $\alpha>0$
and $\operatorname{diag}(A)$ is the diagonal part of $A$. Note that the diagonal part of $A$ never contains a zero value.

Adding fill-ins may obviously lead to a better preconditioner and a potentially
better convergence rate. On the other hand, it becomes computationally more expensive to compute
the preconditioner itself. Thus, there is a trade-off between the latter effect and the 
improved convergence rate, an important issue upon which we will elaborate for our application.

\subsection{On the FM-PCG Normal Integrator}

Caused by its local nature the reconstructions computed by FM often have a lower quality compared to results of 
global approaches. On the other hand, the empowering effect of preconditioning the Poisson integration
may still not suffice for achieving a high efficiency. The basic idea we follow now is that if one starts 
the PCG with a proper initialisation $\mathbf{x}_0$ obtained by FM integration, instead of the standard case $\mathbf{x}_0=0$, 
the PCG normal integrator could benefit from a significant speed-up. 
This idea together with a dedicated numerical evaluation leading to the definition of a well-engineered choice
of computational parameters of the numerical PCG solver is the core of our proposed method.

In the following, we first determine the important building blocks and parameters of our algorithm.
This is done in Section \ref{evaluation}. It will become evident how the \emph{individual} methods 
perform and how they compare.

After that we will show in Section \ref{PCG} that our proposed FM-PCG normal integrator in which 
\emph{suitable building blocks are put together} is highly competitive and in many instances superior to the 
state-of-the-art methods for surface 
normal integration.


\section{Numerical Evaluation}\label{evaluation}

We now demonstrate relevant properties of several 
state-of-the-art
 methods for  
surface normal integration. For this purpose, we give a careful evaluation regarding the 
accuracy of the reconstruction, the influence of boundary conditions, the flexibility to handling 
non-rectangular domains, the robustness against noisy data and the computational efficiency 
-- the main challenges for an advanced surface normal integrator. 
On the technical side, let us note that the experiments were 
conducted on a I7 processor at $2.9GHz$.

\paragraph{Test Datasets}
To evaluate the proposed surface normal integrators, we represent examples of possible
applications in gradient-domain image reconstruction (PET imaging, Poisson image editing) and 
surface-from-gradient (photometric stereo).
The gradients of the ``Phantom'' and ``Lena'' images were constructed using finite differences, 
while both the surface and the gradient of the ``Peaks'', ``Sombrero'' and ``Vase'' datasets are 
analytically known, preventing any bias due to finite difference approximations. 

Let us note that our test datasets also address fundamental
aspects one may typically find in gradient fields obtained from real-world problems: 
sharp gradients (``Phantom''), highly fluctuating gradients oriented in all grid directions
representing for instance textured areas (``Lena'') and smoothly varying gradient fields (``Peaks''). 
The gradient field of the ``Vase'' dataset shows difficulties like a non-trivial computational domain.

\subsection{Existing Relevant Integration Methods}

The fast \emph{and} accurate surface normal integrators are not abundant.
For a meaningful assessment we have to compare our novel FM-PCG approach
with the Fast Fourier Transform (FFT) method of Frankot and Chellappa \cite{Frankot1988} and the 
Discrete Cosine Transform (DCT) extended by Simchony {\em et al.} \cite{Simchony1990} which are 
two of the most popular methods in science and industry. Furthermore, we include 
the recent method of Harker and O'Leary \cite{Harker2015} which relies on the formulation of
the integration problem as a Sylvester equation. It will be helpful to consider in a first step
the building blocks of our approach, i.e.\ FM and CG-based Poisson integration, in a separate way. Hence, we also include in our comparison the FM method from~\cite{Galliani2012}. As for Poisson integration, only Jacobi~\cite{Durou2,Durou2009} and Gauss-Seidel~\cite{Queau2015} iterations were employed so far, hence we consider \cite{Hestenes1952} as reference for CG-Poisson integration.

To highlight the differences between the methods, we establish for the beginning a
comparison with respect to their algorithmic complexity, the type of admissible boundary conditions 
and the requirement for the computational domain $\Omega$, see Table \ref{tab:1}. While the algorithmic 
complexity is an indicator for the speed of a solver, the admissible boundary conditions and the handling of 
non-rectangular domains influence 
its accuracy. The ability to handle non-rectangular domains 
improves also its computational efficiency.
\begin{savenotes}
\begin{table}[!h]
\caption{
Comparison of five existing fast and accurate surface normal integration methods
based on three 
criteria: their algorithmic complexity w.r.t.\ the number $n$ of pixels inside the computational 
domain $\Omega$ (the lower the better), the type of boundary condition (BC) they use 
(free boundaries are expected to reduce bias), and the requirement for $\Omega$ to be rectangular or not
(handling non-rectangular domains can be important for accuracy and algorithmic speed). }
  \vspace {2.5mm}
    \label{tab:1}
    \setlength{\tabcolsep}{4.5pt}
    \centering\begin{tabular}{ccccc}
        \toprule
        Method & Ref. & Complexity & BC & Non-rect. \\ \hline
        FFT & \cite{Frankot1988} & $n \log n$ & periodic & no \\
        DCT & \cite{Simchony1990} & $n \log n$ & \textbf{free} & no \\
        FM & \cite{Galliani2012} & $n \log n$
 & \textbf{free} & \textbf{yes} \\
        Sylvester & \cite{Harker2015} & $n^{\frac{3}{2}}$~\footnotemark\footnotetext{Assuming $\Omega$ is square. For rectangular domains of size $n_r \times n_c$, the complexity is $O(n_c ^3)$.} & \textbf{free} & no \\
        CG-Poisson & \cite{Hestenes1952} & $n^3$~\footnotemark\footnotetext{Without using preconditioning techniques. } & \textbf{free} & \textbf{yes} \\
        \bottomrule
    \end{tabular}
\end{table}
\end{savenotes}

The findings in Table \ref{tab:1} already indicate the potential usefulness of a mixture of FM and CG-Poisson as both are free of 
constraints in the last two criteria and so their combination may lead to a reasonable computational efficiency of the 
Poisson solver. Although the other methods have their strengths in the algorithmic complexity and in the application 
of boundary conditions (apart from FFT), we will see that the flexible handling of domains is a fundamental task and 
a key requirement of an ideal solver for surface normal integration.

\subsection{Stopping Criterion for CG-Poisson}\label{stop}

Among the considered methods, the Poisson solver (conjugate gradient method), where solving the discrete Poisson 
equation (\ref{poisson}) corresponds to a linear system $A\mathbf{x}=\mathbf{b}$, is the only iterative scheme. 
As indicated in Section~\ref{sec:krylovsubspaces}, a practical solution can be reached quickly after a small number $k$ of iterations, whereby $k$ cannot be 
verified exactly. The general stopping criterion of an iterative method can be based on the \emph{relative residual} 
$\frac{\|\mathbf{b}-A\mathbf{x}\|}{\|\mathbf{b}\|}$ which we analyse in this paragraph. 

To 
guarantee the efficiency of the CG-Poisson solver it is necessary to define the number $k$ of iterations depending on 
the quality of the reconstruction in the iterative process. To tackle this issue we examined the 
comparison between the 
 MSE\footnote{The mean squared error (MSE) is used to quantify the error of the reconstruction. 
We employed it to estimate the amount of the error contained in the reconstruction compared to the original.} 
and the relative residual during each CG iteration. 
Let us note that due to the fact that the solution of the linear system (\ref{system-1}) is not unique an 
additive ambiguity $v \mapsto v+c,\, c \in \mathbb{R}$ in the integration problem ($c$ is the ``integration constant") occurs. 
Therefore, in each numerical experiment we chose the additive constant $c$ which minimises the MSE, for fair comparison.  

\begin{figure}[!h]
  \caption{
MSE vs.\ relative residual during CG iterations, for the ``Sombrero" dataset. Although arbitrary relative 
accuracy can be reached, it is not useful to go beyond a $10^{-3}$ residual, since such refinements have very 
few impact on the quality of the reconstruction, as shown by the MSE graph. Similar results were obtained for all 
the datasets used in this paper. Hence, we set as stopping criterion a $10^{-4}$ relative residual, which can be 
considered as ``safe". }
  \label{fig:MSE_residual}  
  \vspace {-4.5mm}    
  \begin{center}
    \begin{tabular}{cc}
      \includegraphics[width = 0.49\linewidth]{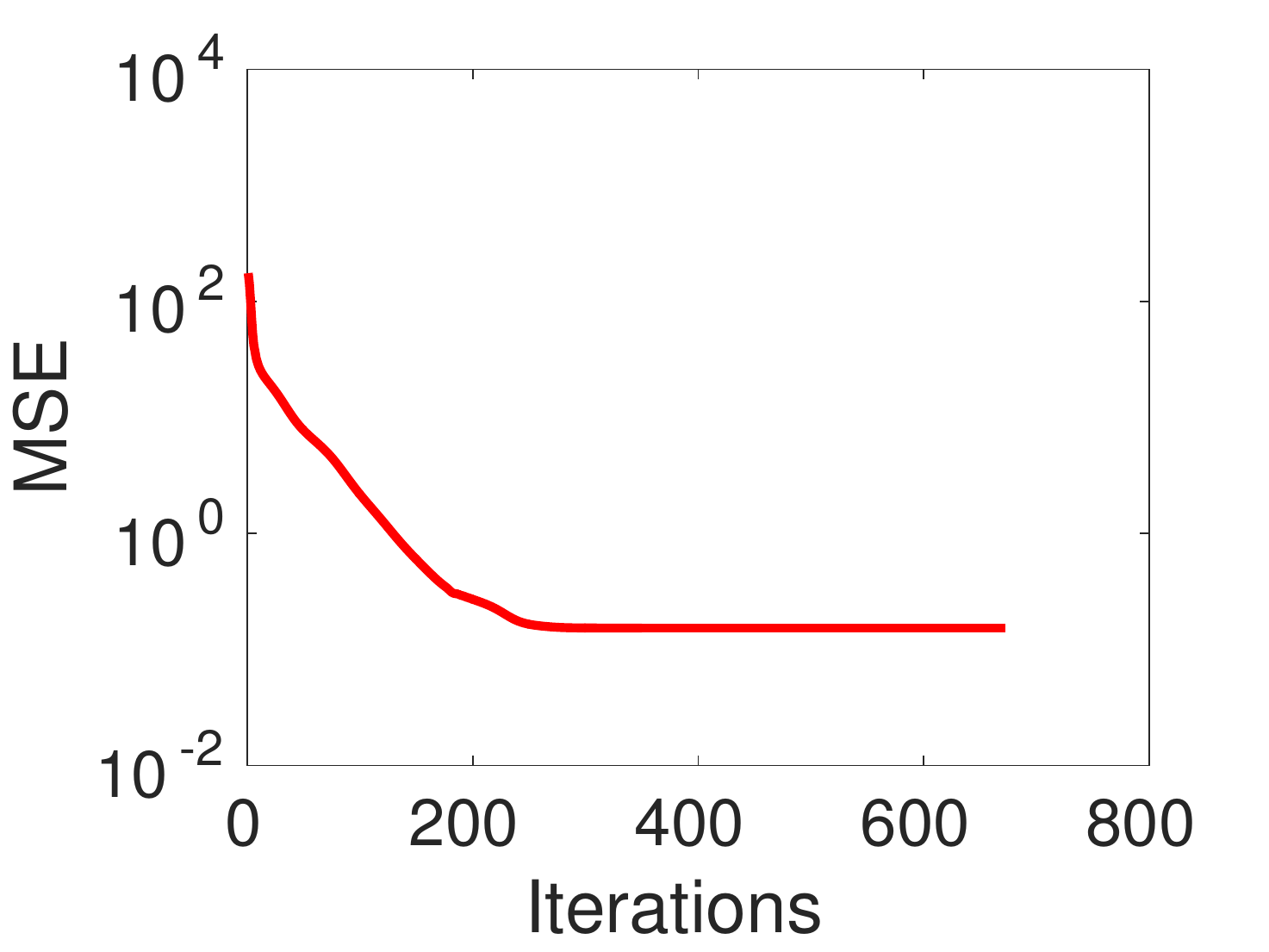} &\hspace{-0.3cm}
      \includegraphics[width = 0.49\linewidth]{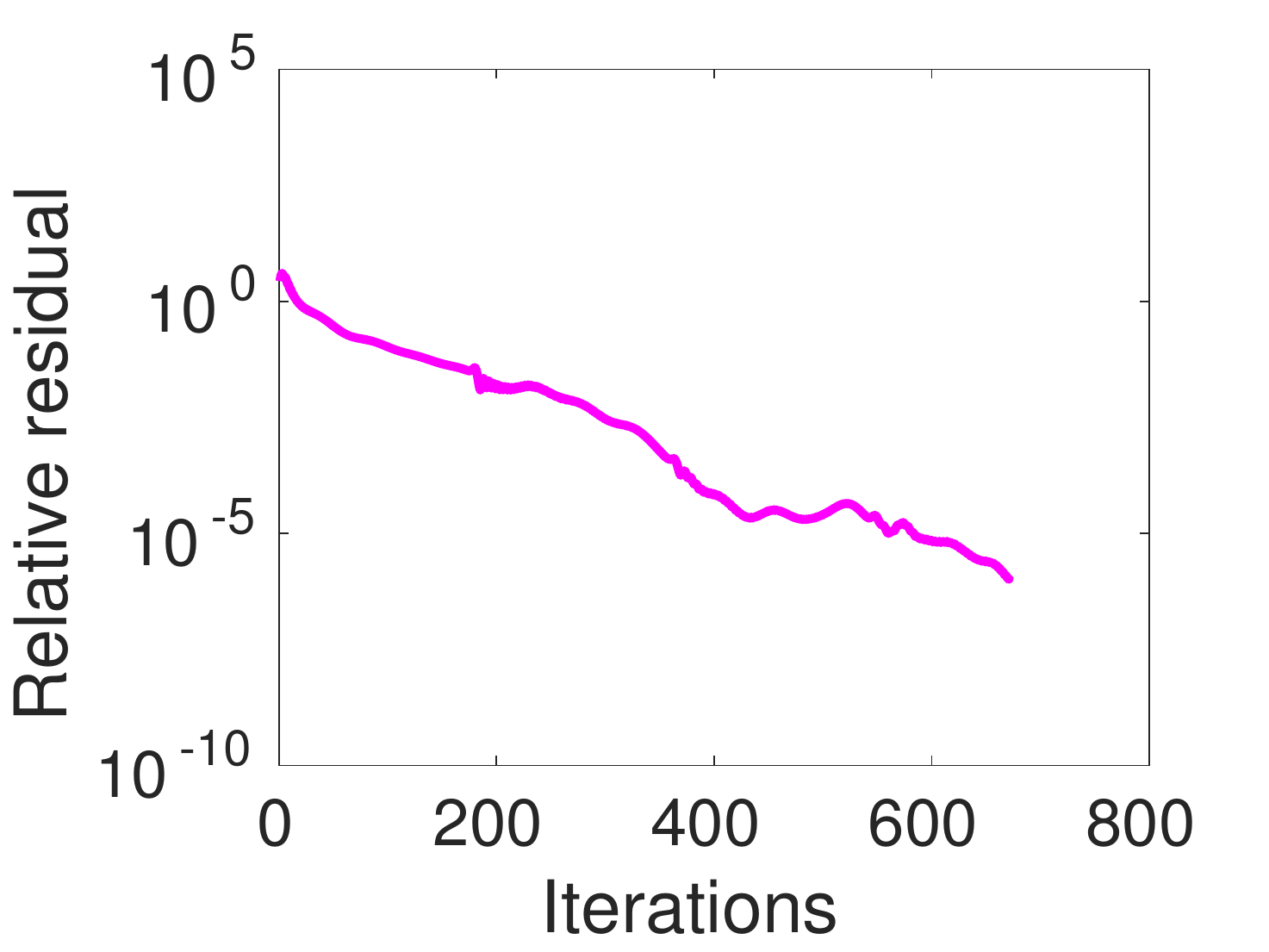} 
    \end{tabular}
  \end{center}
\end{figure}

\begin{figure}[!h]
  \caption{
The ``Sombrero" surface ($256 \times 256$) used in this experiment, whose gradient can be calculated analytically. 
Note that the depth values are periodic on the boundaries.}
  \label{fig:MSE_residual_Sombrero}
  \vspace {2.5mm}  
 \centering \includegraphics[width = 0.8\linewidth]{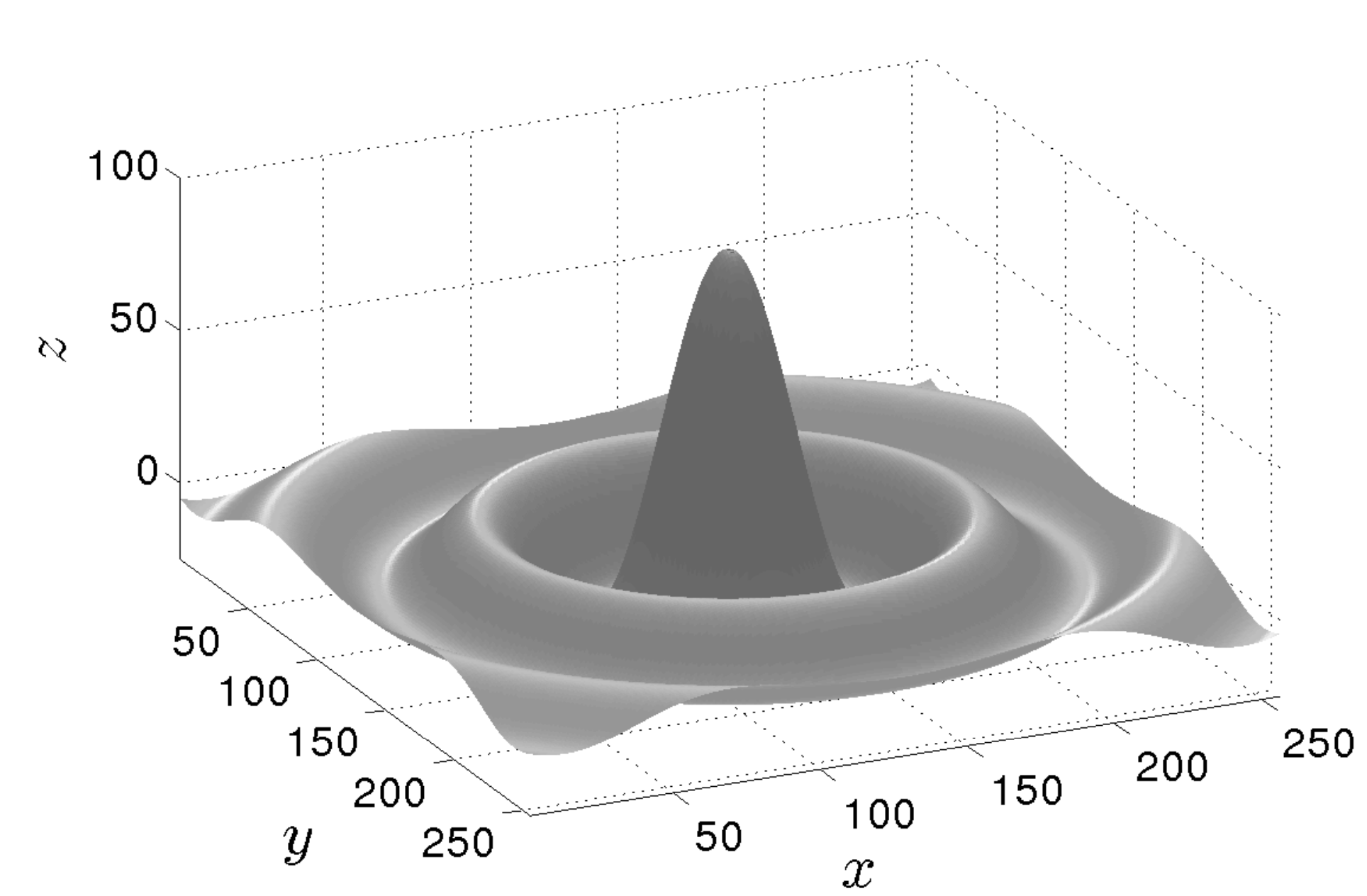}
\end{figure}

To determine a proper relative residual, we verified the datasets ``Lena'', ``Peaks'', ``Phantom'', ``Sombrero'' and ``Vase'' 
on rectangular and non-rectangular domains. All test cases showed similar results like the graphs 
in Figure \ref{fig:MSE_residual} for the reconstruction of the ``Sombrero'' 
surface (see Figure \ref{fig:MSE_residual_Sombrero}).

In 
this experiment the iterative solver CG-Poisson 
is stopped when the relative residual 
is lower than $10^{-6}$. However, it can be 
seen clearly, that after around iteration 250 the quality measured by the MSE cannot be improved and therefore 
more than 400 iterations are without success. This numerical steady state of the MSE and therefore of the residual will be 
found in regions 
where the relative residual is between $10^{-3}$ and $10^{-4}$, thus we consider the stopping criterion of $10^{-4}$ for 
the subsequent experiments as suitable and ''safe``.


\subsection{Accuracy of the Solvers}

First of all we analyse the general quality of the methods listed in Table \ref{tab:1} for the 
 ``Sombrero" 
dataset over a quadratic domain of size $256\times 256$. In this example the gradient can be calculated in 
an analytical way and furthermore 
the boundary condition is periodic.
Consulting Table \ref{tab:1} it 
is obvious that all methods have no restrictions and consequently no discrimination. 

Basically all methods, see Figure {\ref{fig:results_sombrero}, provide a satisfactory reconstruction, only 
FM produces a less accurate solution. 
This can be viewed more easily in Table \ref{tab:results_sombrero}, where the 
values of the measurements\footnote{We tested the two common measurements MSE and SSIM. 
A superior reconstruction has a small value (tends to zero) for MSE and a value close to one 
(tends to one) for the SSIM.} 
MSE and SSIM\footnote{The structural similarity (SSIM) index is a method for predicting the perceived 
quality of an image, see~\cite{Wang2004}.} and the CPU time (in seconds) are illustrated.  The 
accuracy of all methods is similar,
merely the solution of FM 
with $1.16$ for MSE and $0.98$ for SSIM 
 is slightly worse. In contrast the CPU times vary strongly, in 
this case FFT and DCT, which need around $0.01~s$,
are unbeatable. The times for FM and Sylvester are in
a reasonable range, solely the standard CG-Poisson with around $1.06~s$
 shows
unacceptable running costs and is consequently quite inefficient. 
For problems of surface reconstruction with the indicated conditions the choice of a solver is fairly 
easy, namely the frequency domain methods FFT and DCT. 
\begin{table}[!h]
    \caption{Results on the ``Sombrero" dataset ($256 \times 256$). As expected, all methods provide 
reasonably accurate solutions. Yet, one can remark that the FM result is slightly less accurate: 
this is due to error accumulation by the \emph{local} nature of FM, while the other methods are 
\emph{global}. The reconstructed surfaces are shown in Figure \ref{fig:results_sombrero}.}
  \label{tab:results_sombrero}    
    \vspace {2.5mm}
    \setlength{\tabcolsep}{4.5pt}
    \centering\begin{tabular}{cccc}
        \toprule
        Method & MSE ($px$) & SSIM & CPU ($s$) \\ \hline
        FFT \cite{Frankot1988} & 0.17 & \textbf{1.00} & $\mathbf{<}$ \textbf{0.01} \\
        DCT \cite{Simchony1990} & 0.19 & \textbf{1.00} & 0.01 \\
        FM \cite{Galliani2012} & 1.16 & 0.98 & 0.07  \\
        Sylvester \cite{Harker2015} & \textbf{0.15} & \textbf{1.00} & 0.18 \\
        CG-Poisson \cite{Hestenes1952} & \textbf{0.15} & \textbf{1.00} & 1.06 \\
        \bottomrule
    \end{tabular}
\end{table}

\begin{figure}[!h]
  \caption{Results on the ``Sombrero" dataset (cf. Table~\ref{tab:results_sombrero}).}
  \label{fig:results_sombrero}  
  \vspace {1mm}    
  \begin{center}
    \begin{tabular}{cc}
      \includegraphics[width = 0.45\linewidth]{2_sombrero-eps-converted-to.pdf} & 
      \includegraphics[width = 0.45\linewidth]{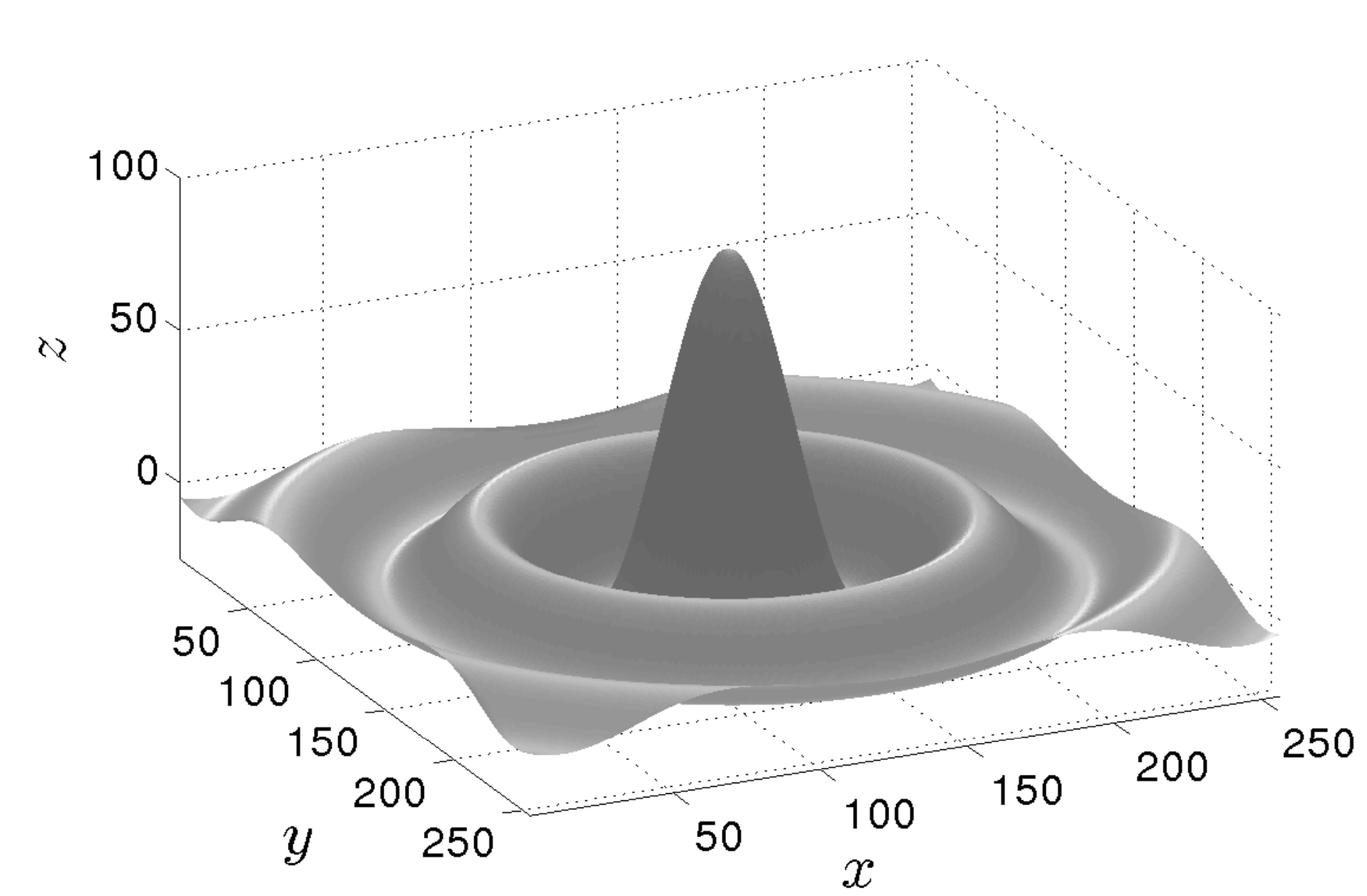} \\
      Ground truth & FFT~\cite{Frankot1988} \\
      \includegraphics[width = 0.45\linewidth]{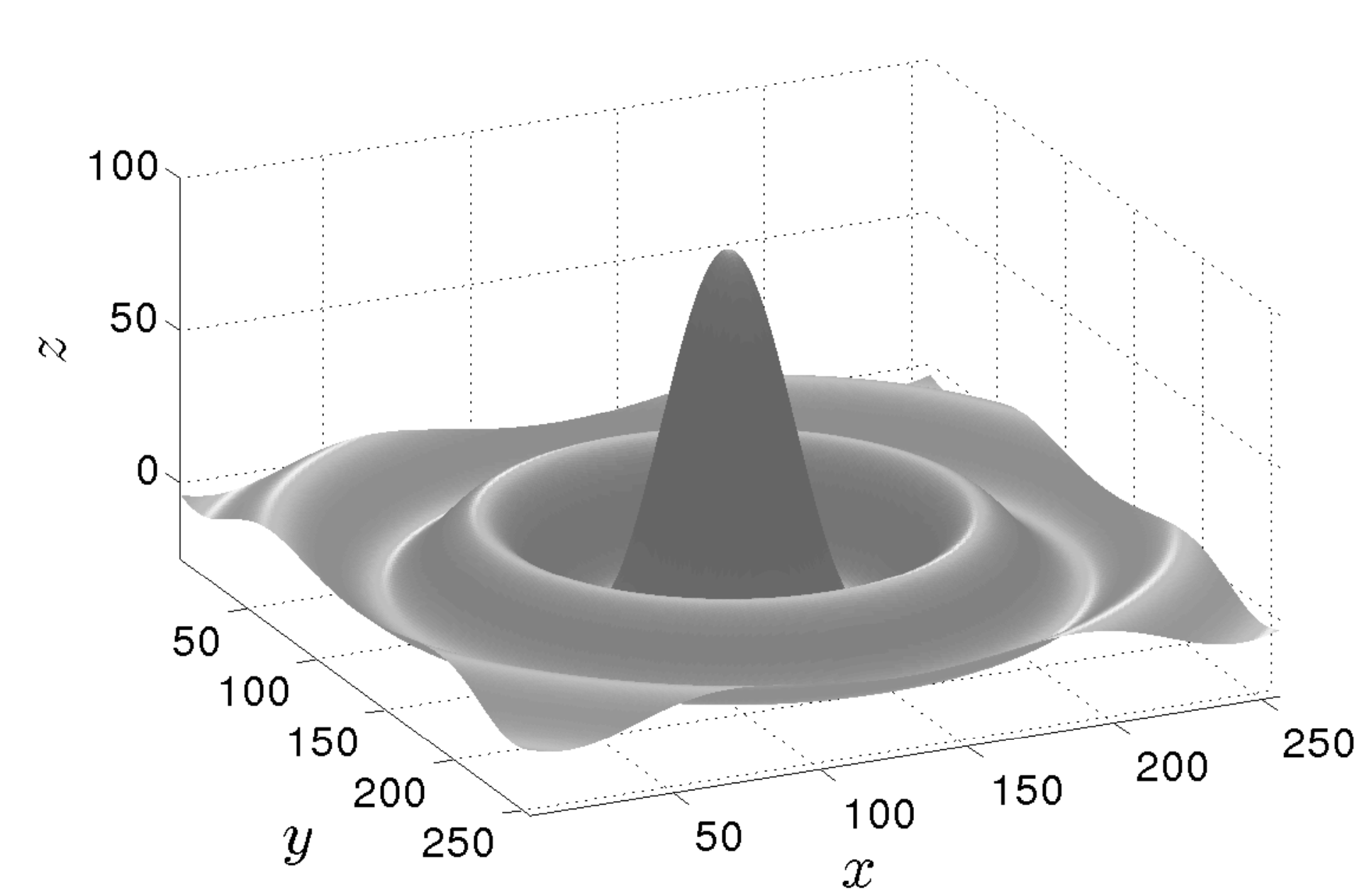} &
      \includegraphics[width = 0.45\linewidth]{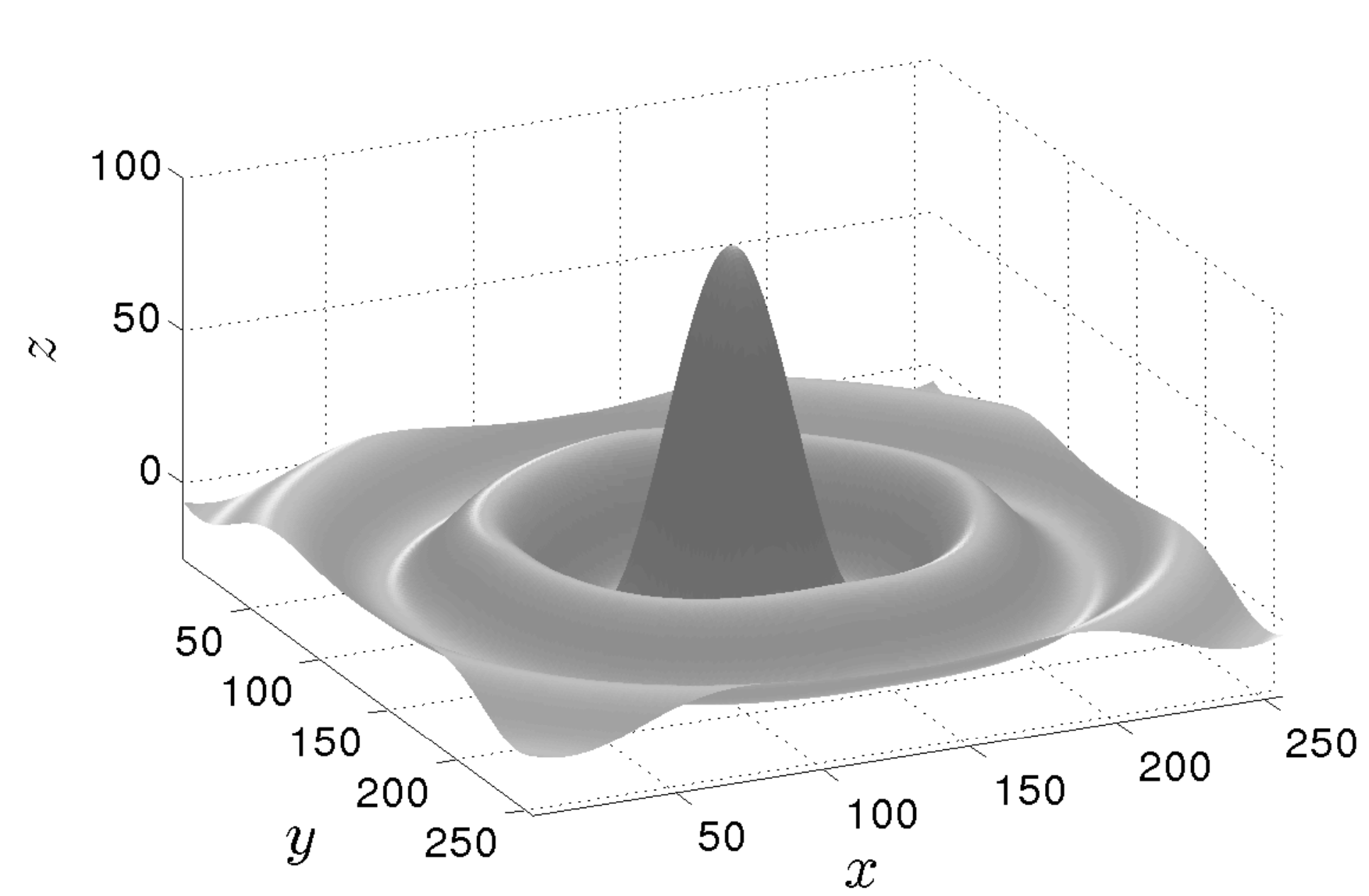} \\      
      DCT~\cite{Simchony1990} & FM~\cite{Galliani2012} \\
      \includegraphics[width = 0.45\linewidth]{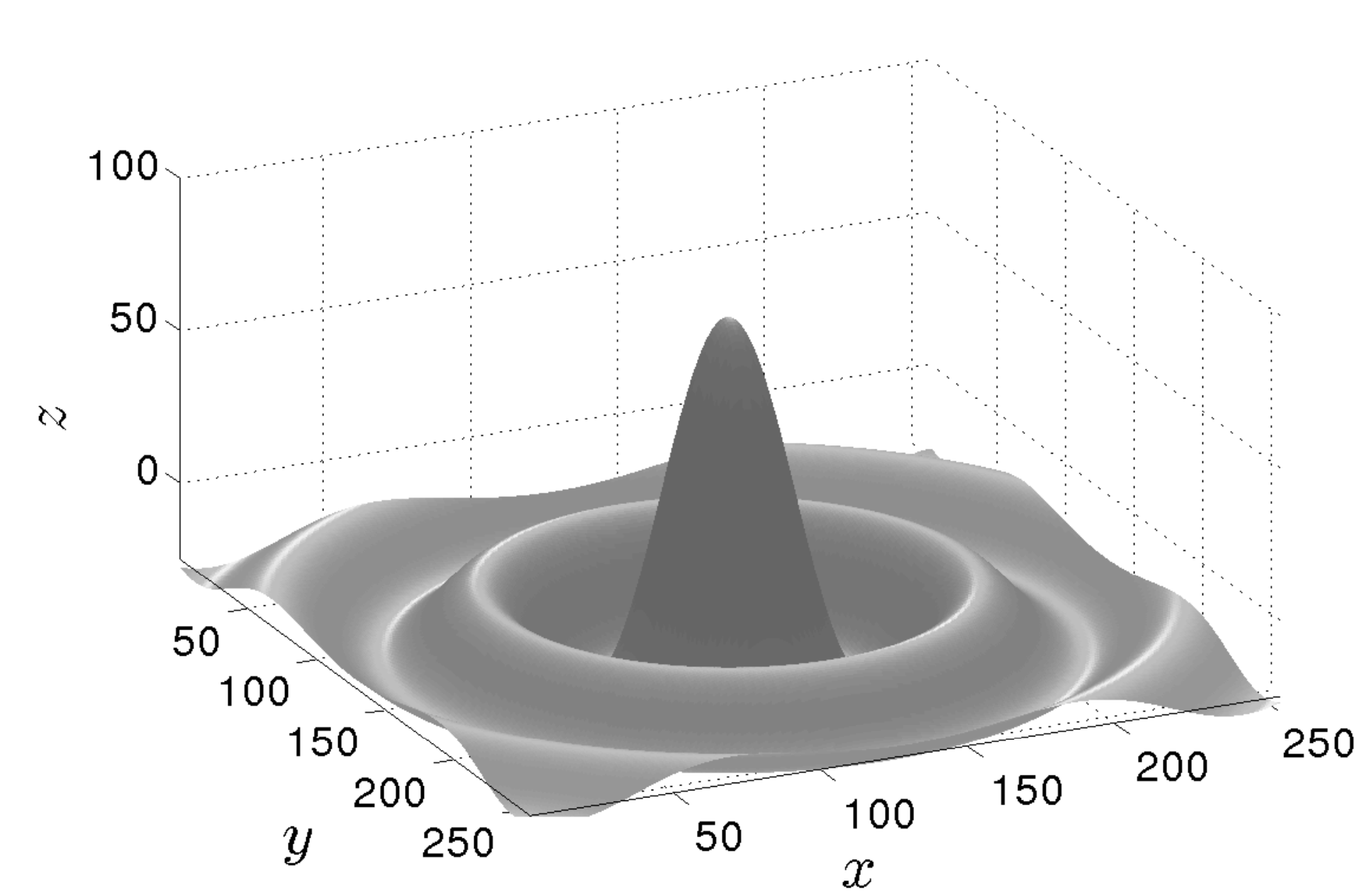} &
      \includegraphics[width = 0.45\linewidth]{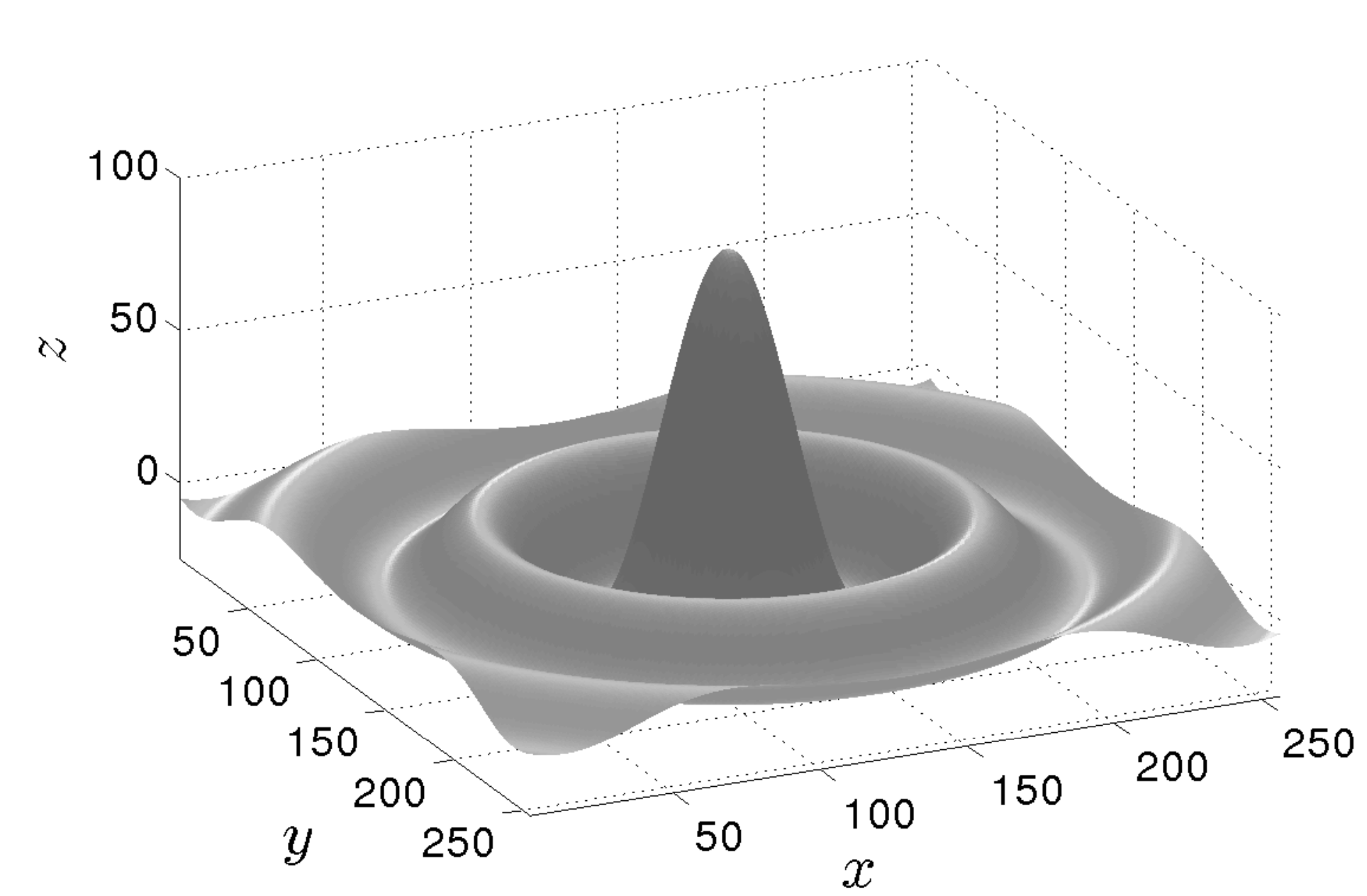} \\
      Sylvester~\cite{Harker2015} & CG-Poisson~\cite{Hestenes1952}           
    \end{tabular}
  \end{center}
\end{figure}



\subsection{Influence of Boundary Conditions}

The handling of boundary conditions is a necessary issue which cannot be ignored. As we will show, different 
boundary conditions lead to surface reconstructions of different accuracy. The assumption of 
Dirichlet, periodic or homogeneous Neumann boundary conditions is often not justified and may even be unrealistic in some 
applications. A better choice is to use the ``natural'' boundary condition
 (see \cite{Horn}) which 
corresponds to the Neumann 
type. 

The behaviour of the discussed solvers for unjustified boundary conditions, particularly for FFT, 
is illustrated by the ``Peaks'' dataset in Figure \ref{fig:results_peaks} and the associated 
Table \ref{tab:results_peaks}. Almost all methods 
provide good reconstructions, also the FM result is acceptable. Only FFT, with $7.19$ for MSE and $0.96$ for SSIM,
is strongly inferior and unusable for real surfaces, since the accuracy will be lost. 

Basically our results show that FFT-based methods enforcing periodic boundary conditions can be discarded
from the list of candidates for an ideal solver. Once again CG-Poisson is a very accurate integrator, but 
DCT and Sylvester are much faster and provide useful results. However, we may point out again in advance
that 
enforcing the domain $\Omega$ to be rectangular
may lead to difficulties w.r.t.\ the transition from foreground
to background of an object.
\begin{table}[!h]
    \caption{Results on the ``Peaks" dataset ($128 \times 128$). Methods enforcing periodic BC fail at providing a good reconstruction. The reconstructed surfaces are shown in Figure \ref{fig:results_peaks}.}
    \label{tab:results_peaks}
    \vspace {2.5mm}
    \setlength{\tabcolsep}{4.5pt}
    \centering\begin{tabular}{cccc}
        \toprule
        Method & MSE ($px$) & SSIM & CPU ($s$) \\ \hline
        FFT \cite{Frankot1988} & 7.19 & 0.96 & $\mathbf{<}$ \textbf{0.01} \\
        DCT \cite{Simchony1990} & 0.09 & \textbf{1.00}  & $\mathbf{<}$ \textbf{0.01} \\
        FM \cite{Galliani2012} & 0.80 & 0.99 & 0.03 \\
        Sylvester \cite{Harker2015} & \textbf{0.02} & \textbf{1.00} & 0.05 \\
        CG-Poisson \cite{Hestenes1952} & \textbf{0.02} & \textbf{1.00} & 0.29 \\
        \bottomrule
    \end{tabular}
\end{table}

\begin{figure}[!h]
  \caption{Results on the ``Peaks" dataset (cf. Table~\ref{tab:results_peaks}).}
  \label{fig:results_peaks}  
  \vspace {1mm}    
  \begin{center}
    \begin{tabular}{cc}
      \includegraphics[width = 0.45\linewidth]{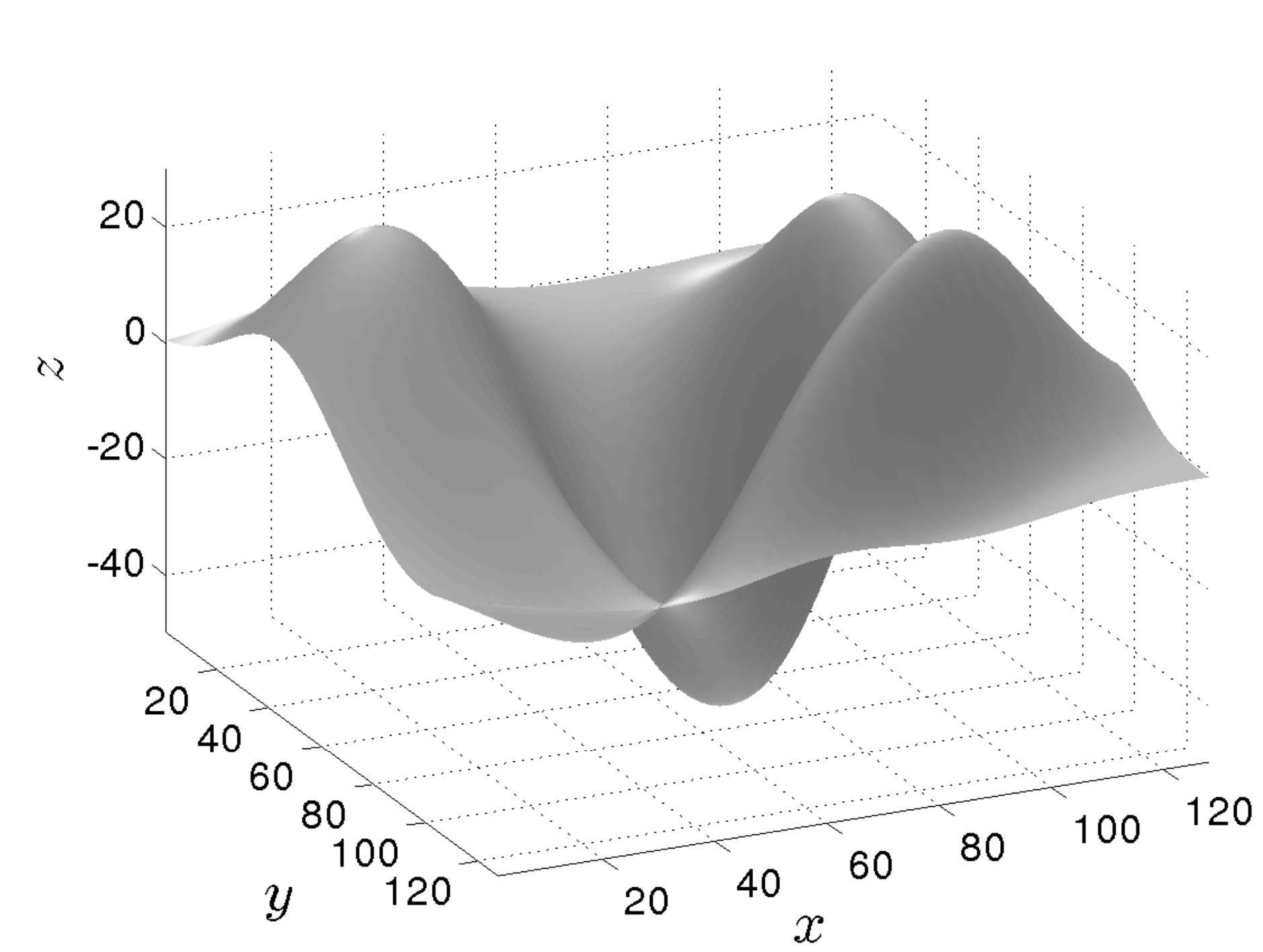} &
      \includegraphics[width = 0.45\linewidth]{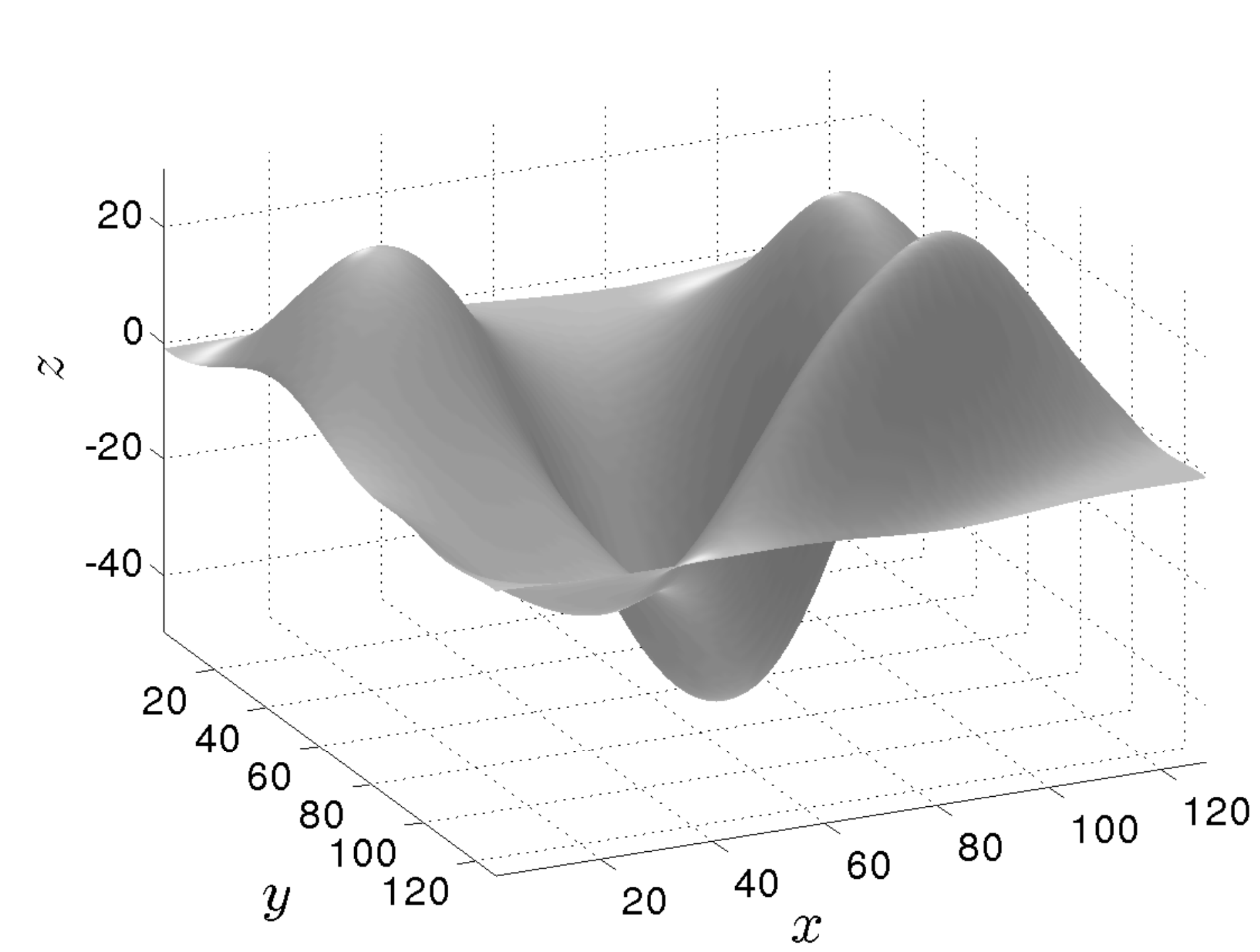} \\
      Ground truth & FFT~\cite{Frankot1988} \\
      \includegraphics[width = 0.45\linewidth]{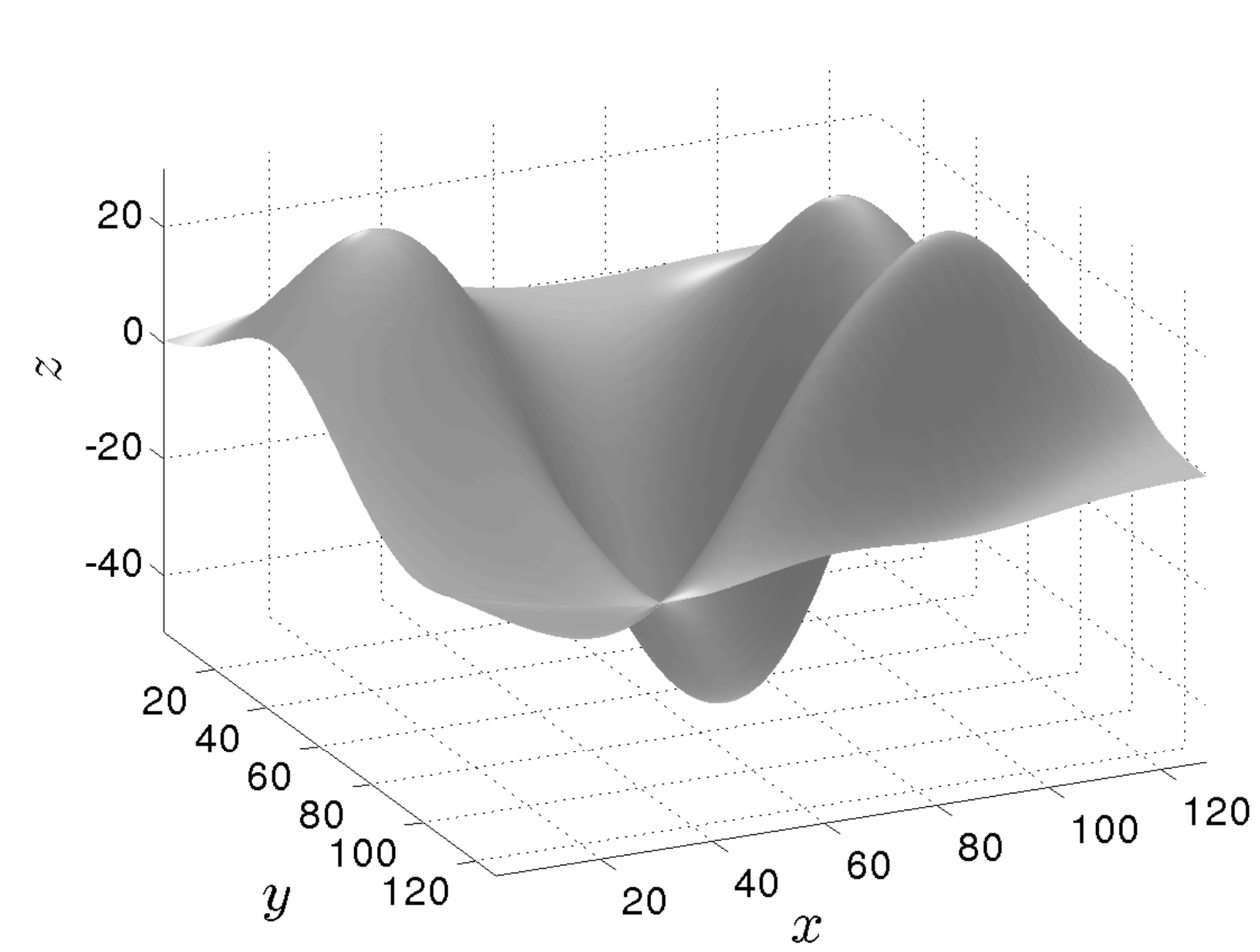} 
    & \includegraphics[width = 0.45\linewidth]{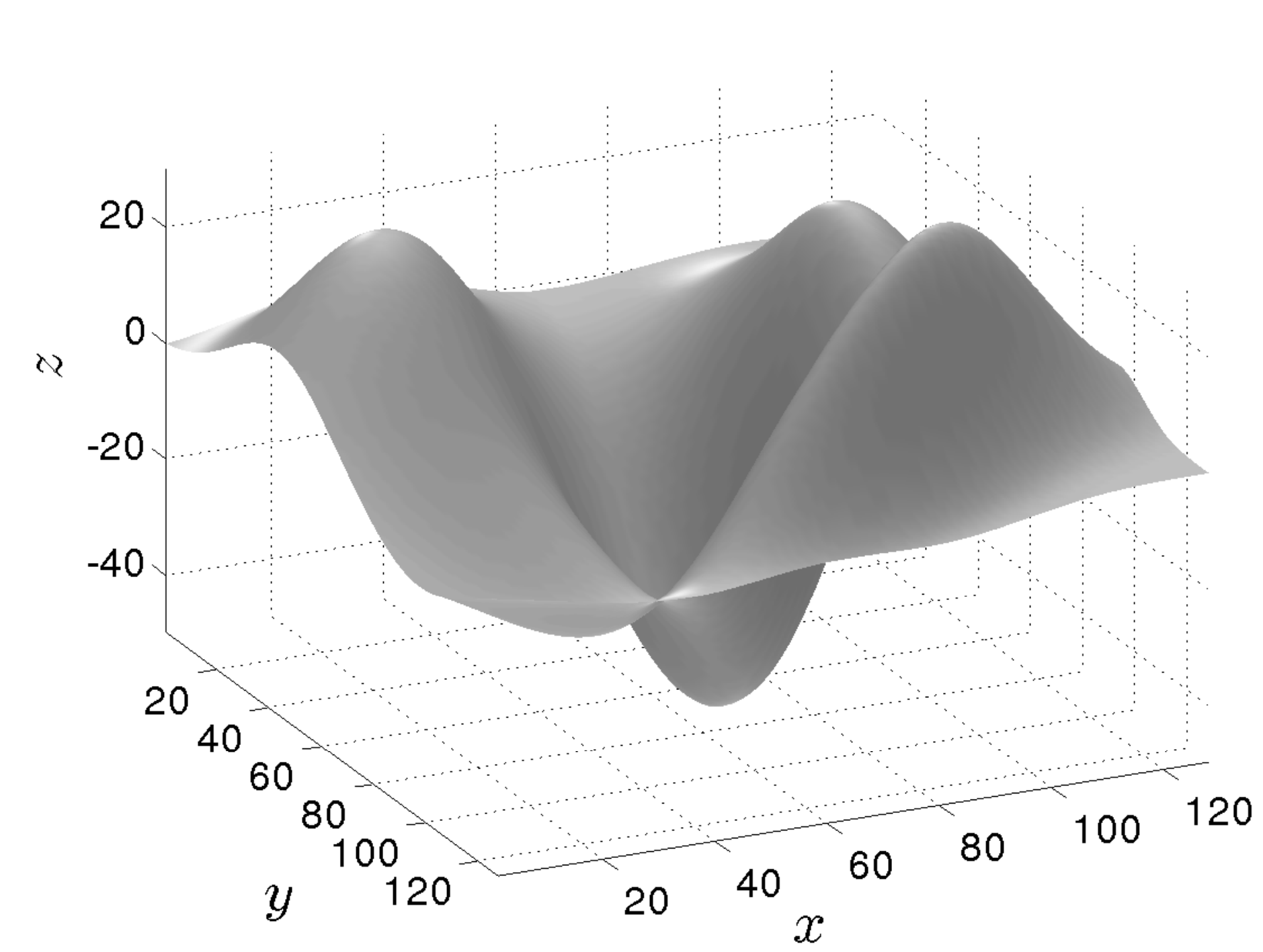} \\
       DCT~\cite{Simchony1990} & FM~\cite{Galliani2012} \\
      \includegraphics[width = 0.45\linewidth]{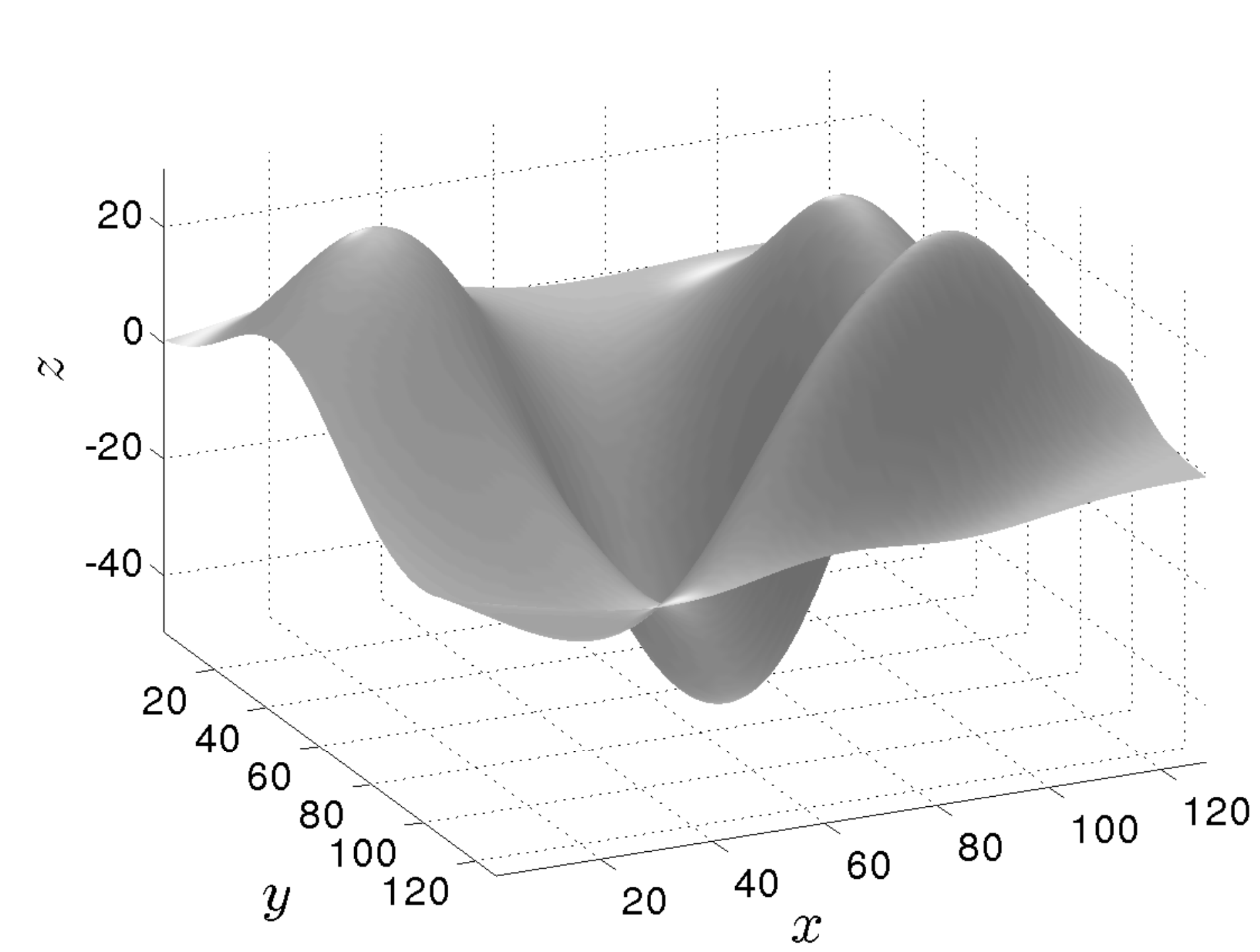} &
      \includegraphics[width = 0.45\linewidth]{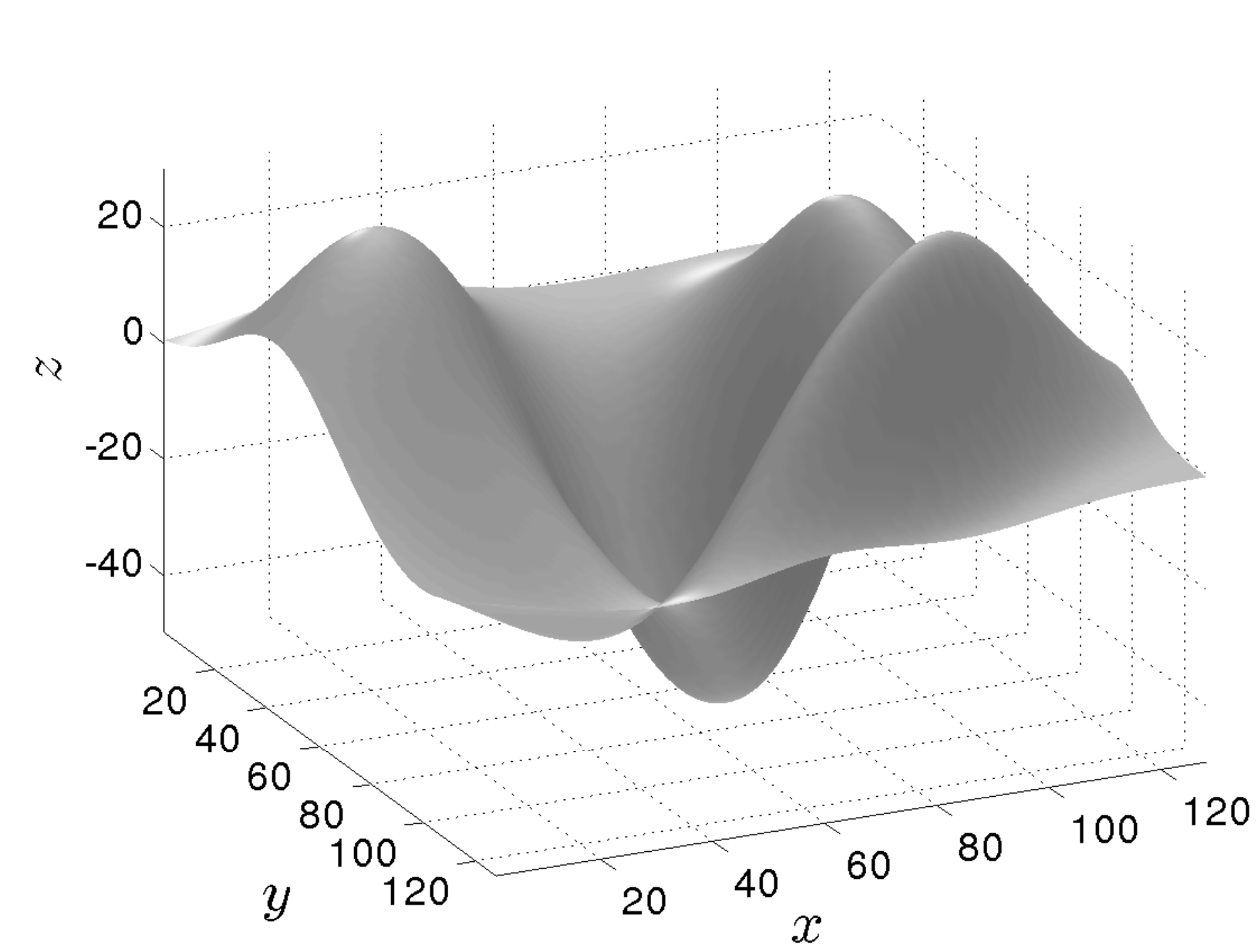} \\
      Sylvester~\cite{Harker2015} & CG-Poisson~\cite{Hestenes1952}           
    \end{tabular}
  \end{center}
\end{figure}



\subsection{Influence of Noisy Data}

An interesting point in many applications is the question of the influence of noise 
on the quality of the reconstructions referring to the different methods. Usually, the correctness 
of the given data, without noise, cannot be guaranteed. 
Therefore, it is essential to have a robust surface normal integrator with respect to noisy data. 

To study the influence of noise, we should consider a dataset which, apart from noise, is perfect. 
Based on this aspect, a very reasonable test example is the ``Sombrero'' dataset, see Figure \ref{fig:results_sombrero}. 
The advantage of ``Sombrero'' is that the gradient of this object is known analytically and not only in 
an approximated way. Furthermore, the computational domain $\Omega$ is rectangular and the boundary conditions 
are periodic. For this test we added to the known gradients\footnote{In the context of photometric stereo, it would be more realistic to add noise to the input images rather than to the gradient~\cite{Noakes2003a}. Nevertheless, evaluating the robustness of integrators to noisy gradients remains useful in order to compare their intrinsic properties.} a Gaussian noise with a standard deviation $\sigma$ varying from $0\%$ to $20\%$ of $\|[p,q]^\top\|_\infty$. 

\begin{figure}[!h]
  \caption{The MSE
  as a function of the standard deviation of a Gaussian noise, expressed in 
percentage of the maximal amplitude of the gradient, added to the gradient. The methods FFT, DCT and CG-Poisson 
  provide the best results for different levels of Gaussian noise. The Sylvester method leads to reasonable 
results for a noise level lower than $10 \%$. 
  Since the FM approach propagates information in a single pass, it obviously also propagates errors, 
inducing a reduced robustness compared to all other approaches.}
  \label{fig:noise}
 \centering \includegraphics[width = 0.7\linewidth]{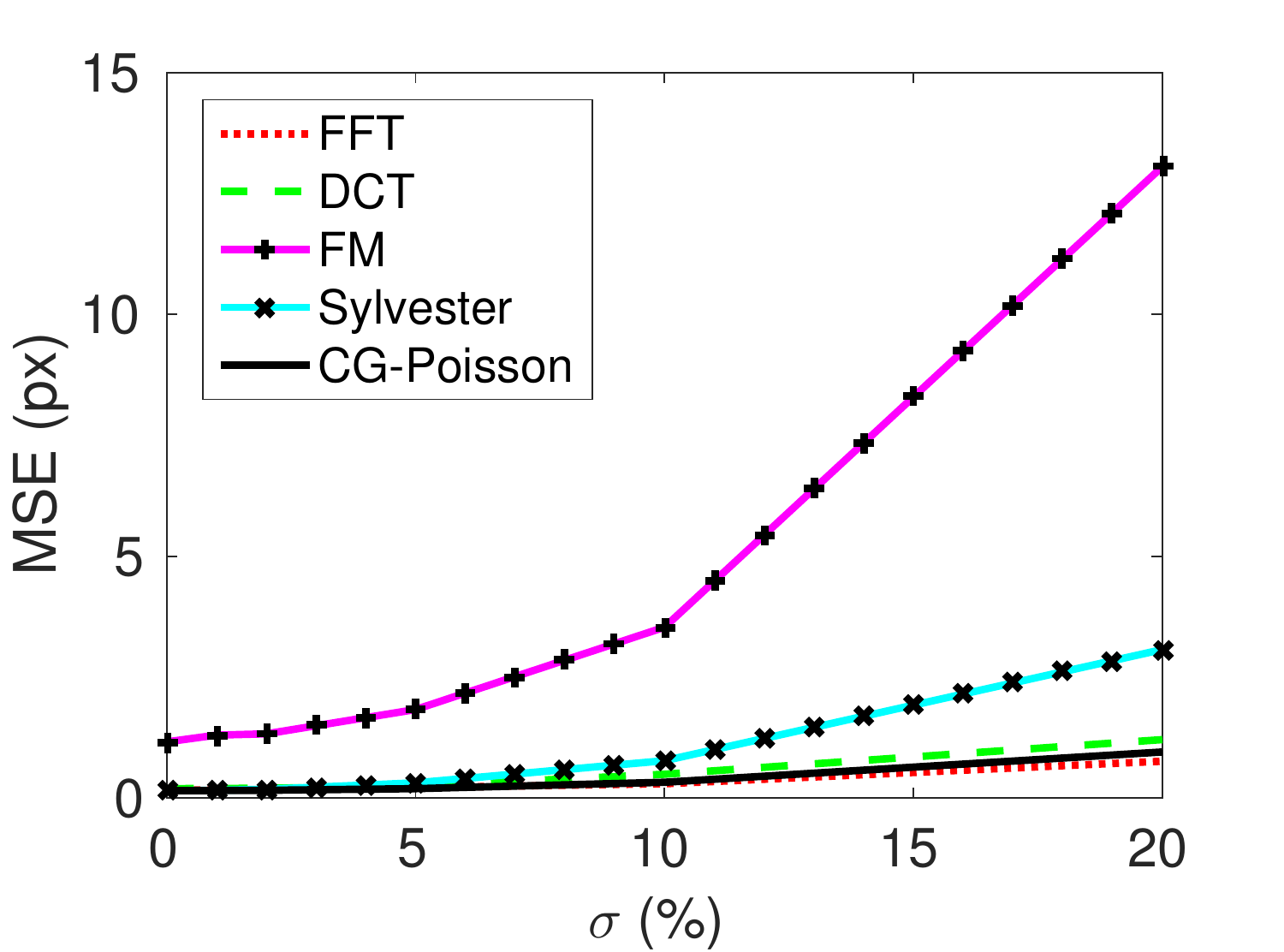}
 \vspace {-2.5mm}  
\end{figure}
The graph in Figure \ref{fig:noise}, which compares the MSE versus the 
standard deviation of Gaussian noise, illustrates the robustness of the tested methods.
The best performance is achieved by FFT, DCT and CG-Poisson, even for strong noisy data with a standard deviation of 20\%. 

The results of Sylvester are similar, however the method contains weaknesses  
in examples with noise of higher standard deviation i.e.\ larger than 10\%. 

As the FM integrator accumulates errors during the front propagation, we observe as expected when increasing the noise 
that this integrator is not a useful choice for highly noisy data. 

To conclude, if the accuracy of the 
given data is not available then FFT, DCT and CG-Poisson are the most secure integrators.

\subsection{Handling of Non-rectangular Domains}
In this experiment we consider the situation when the gradient values are only known on a non-rectangular part 
of the grid. Applying methods dedicated to rectangular grids~\cite{Frankot1988,Simchony1990,Harker2015} requires 
to empirically fixing the values $[p,q]:=[0,0]$ outside $\Omega$, see Figure~\ref{fig:mask_Vase}, inducing a bias.

\begin{figure}[!h]
  \caption{Mask for the evaluation of the ``Vase'' dataset. The gradient values are only known on a non-rectangular 
part $\Omega$ of the grid, which is represented by the white region. The FM and the CG-Poisson integrators can handle easily any form of domain $\Omega$.
In contrast FFT, DCT and Sylvester rely on a rectangular domain and therefore the 
values
$[p,q]:=[0,0]$ need to be fixed outside 
$\Omega$, which is represented by the black region.}
  \label{fig:mask_Vase}  
  \begin{center}
      \includegraphics[width = 0.2\linewidth]{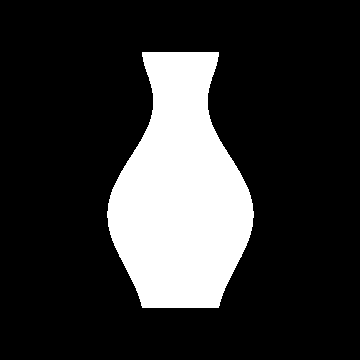}       
  \end{center}
\end{figure}

This can be explained as follows: filling the gradient with null values outside $\Omega$ creates discontinuities 
between the fore- and the background, preventing one from obtaining reasonable results since all the solvers considered here are intended to reconstruct smooth surfaces.
This problem is illustrated in Figure \ref{fig:results_vase} for the ``Vase'' 
dataset with the corresponding Table \ref{tab:results_vase} for the selected measurements MSE and SSIM.

\begin{table}[!h]
    \caption{Results on the ``Vase" dataset ($320 \times 320$). Methods dedicated to rectangular 
    domains are clearly 
biased if $\Omega$ is not rectangular. The corresponding reconstructed surfaces are shown in Figure \ref{fig:results_vase}.}
  \label{tab:results_vase}      
    \vspace {2.5mm}
    \setlength{\tabcolsep}{4.5pt}
    \centering\begin{tabular}{cccc}
        \toprule
        Method & MSE ($px$) & SSIM & CPU ($s$) \\ \hline
        FFT \cite{Frankot1988} & 5.71 & 0.99 & \textbf{0.01} \\
        DCT \cite{Simchony1990} & 5.69 & 0.99 & 0.02 \\
        FM \cite{Galliani2012} & 0.71 & \textbf{1.00} & 0.15 \\
        Sylvester \cite{Harker2015} & 5.99 & 0.98 & 0.36 \\
        CG-Poisson \cite{Hestenes1952} & \textbf{0.01} & \textbf{1.00} & 0.52 \\
        \bottomrule
    \end{tabular}
\end{table}

\begin{figure*}[!ht]
  \caption{Results on the ``Vase" dataset (cf. Table~\ref{tab:results_vase}).}
  \label{fig:results_vase}  
  \vspace {1mm}    
  \begin{center}
    \begin{tabular}{ccc}
      \includegraphics[width = 0.3\linewidth]{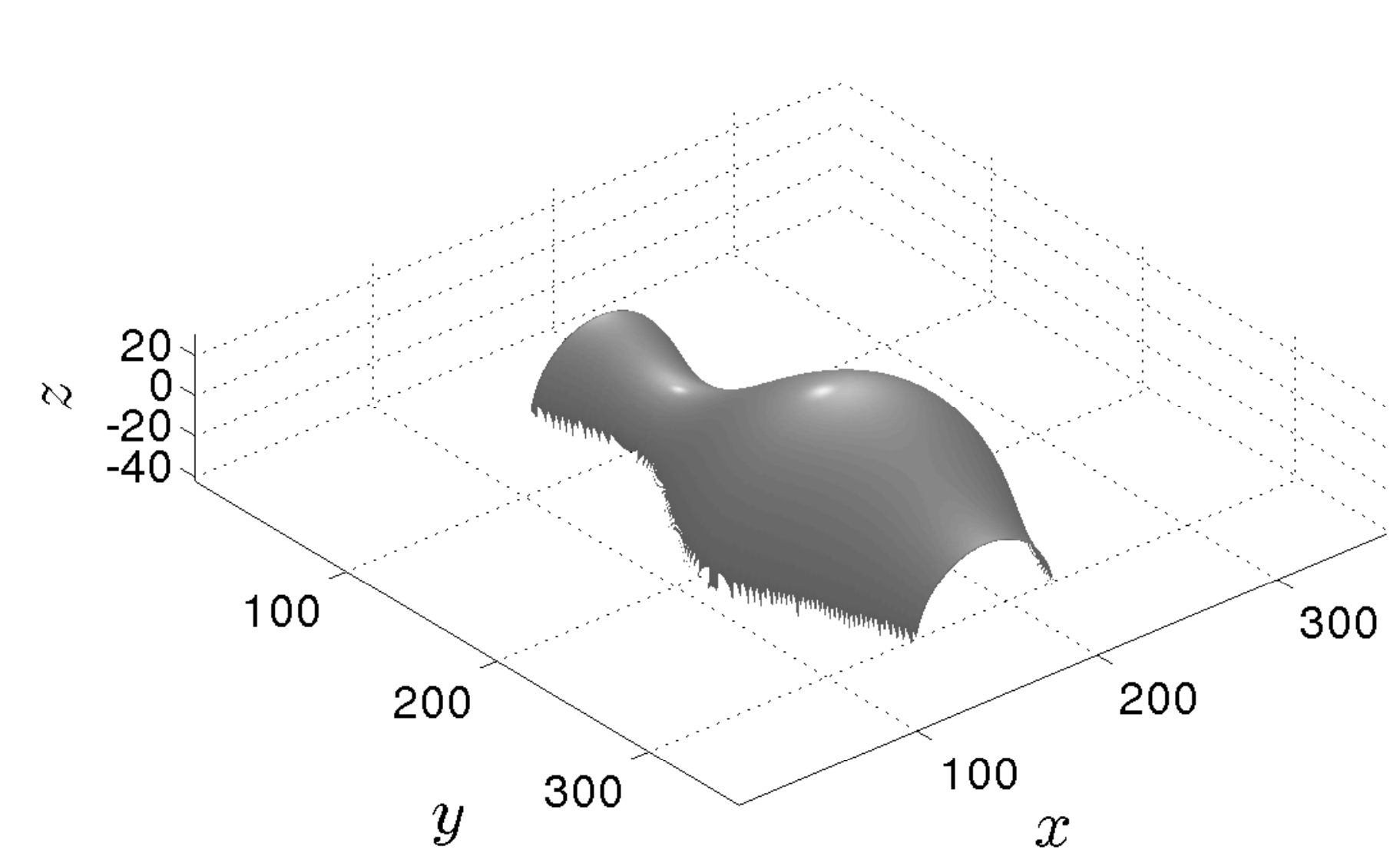} &
      \includegraphics[width = 0.3\linewidth]{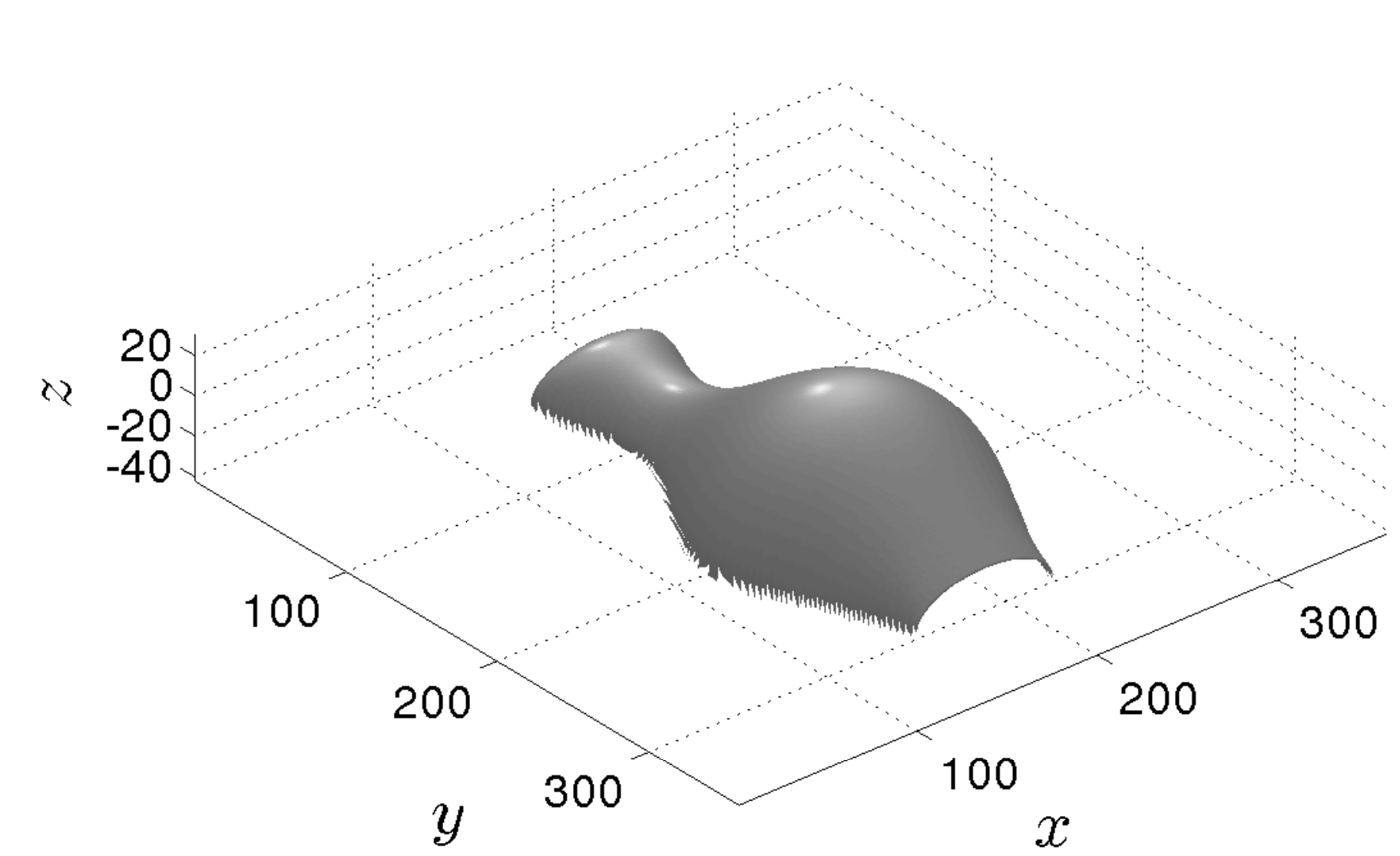} & 
       \includegraphics[width = 0.3\linewidth]{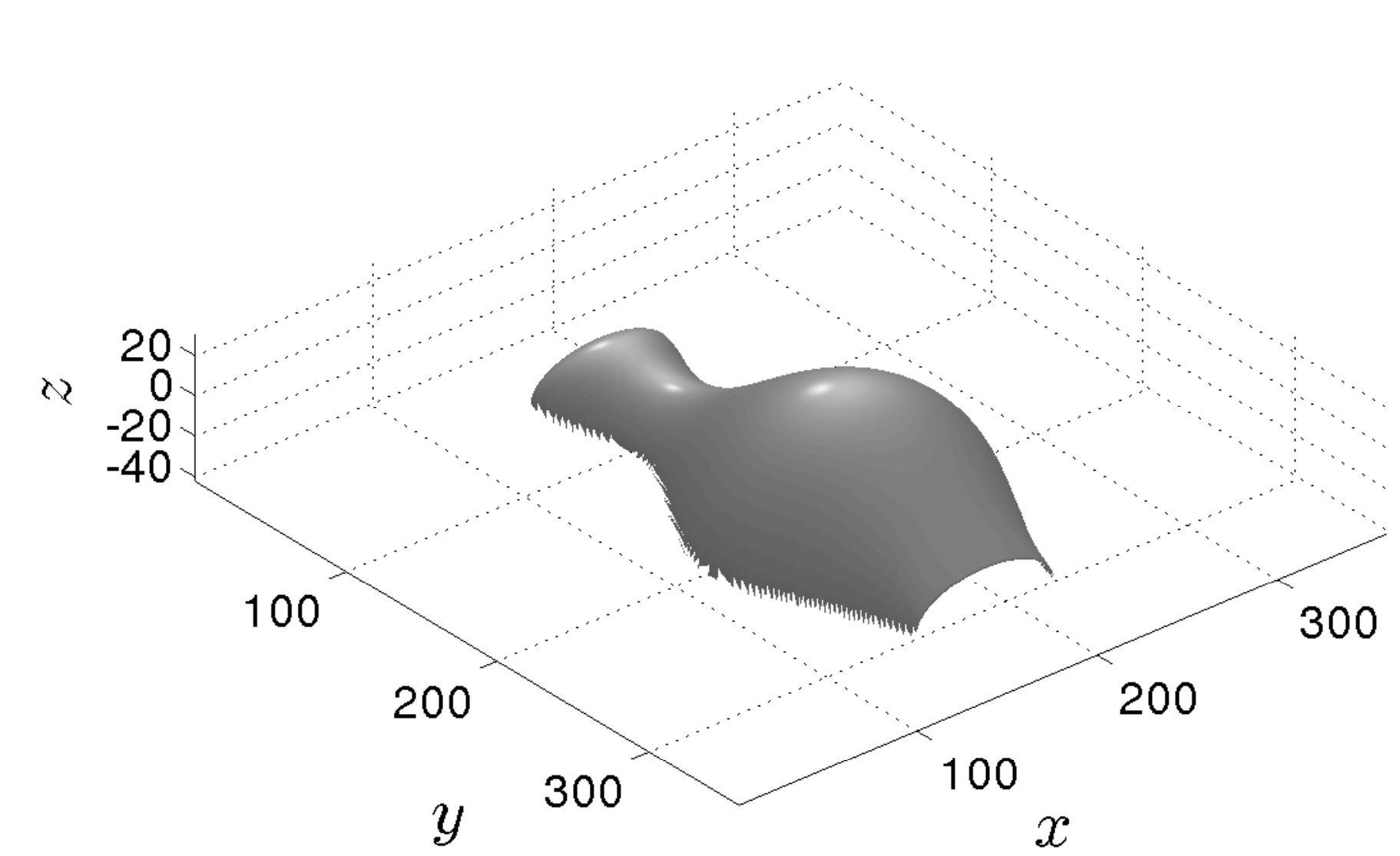}   \\
      Ground truth & FFT~\cite{Frankot1988} & DCT~\cite{Simchony1990} \\ 
     \includegraphics[width = 0.3\linewidth]{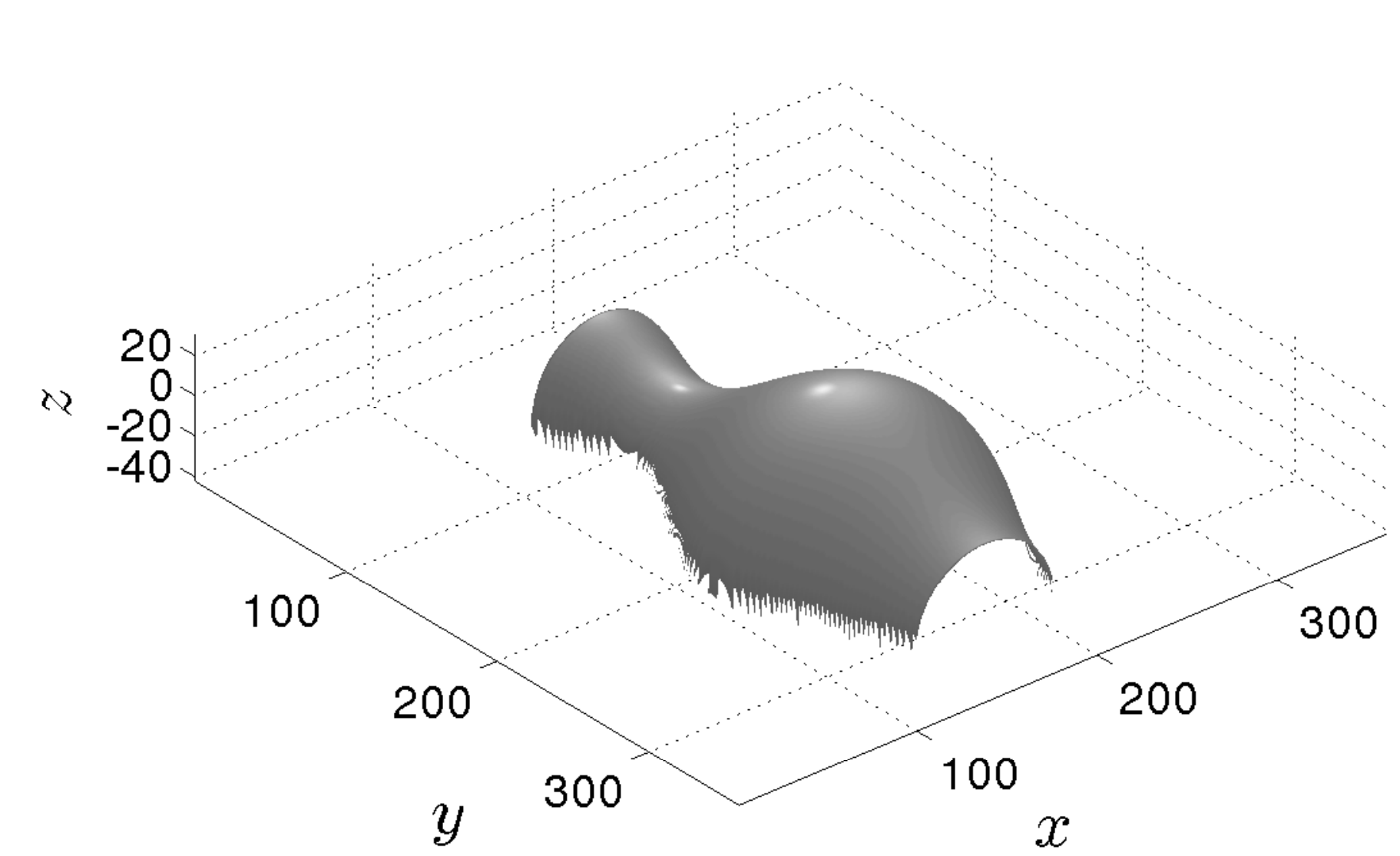} & 
     \includegraphics[width = 0.3\linewidth]{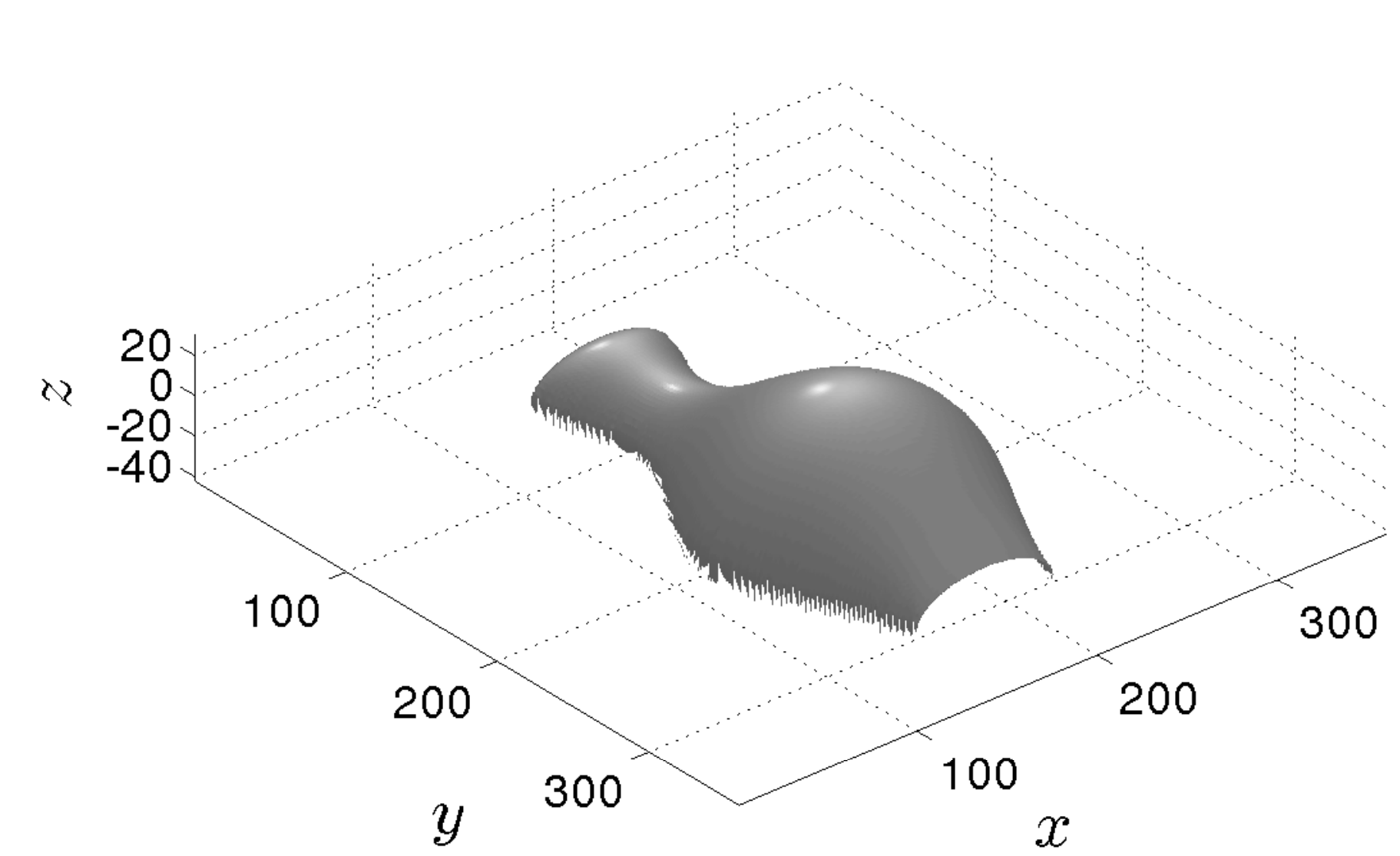} &
      \includegraphics[width = 0.3\linewidth]{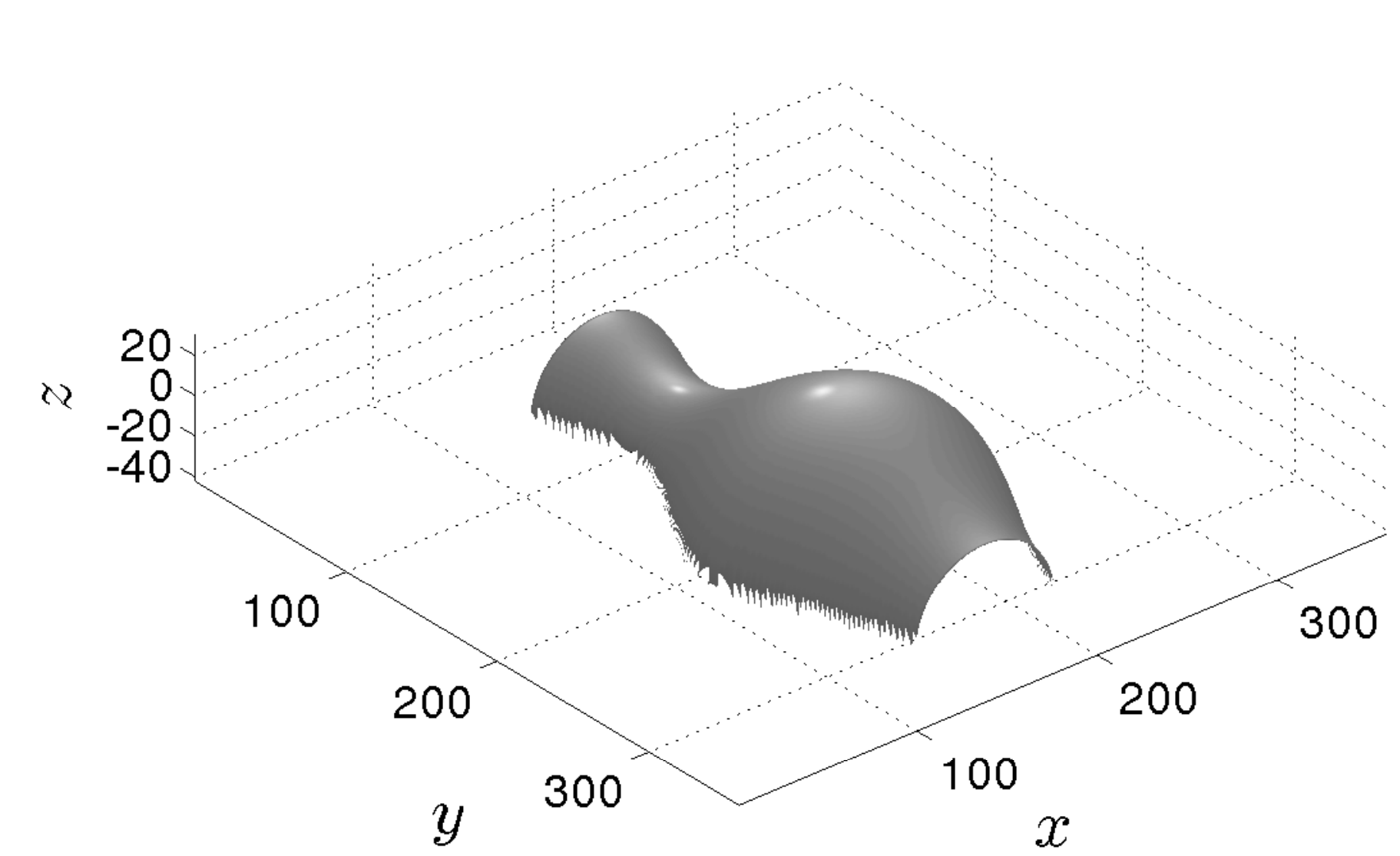} \\
      FM~\cite{Galliani2012}  & Sylvester~\cite{Harker2015} & CG-Poisson~\cite{Hestenes1952}           
    \end{tabular}
  \end{center}
\end{figure*}

The 
methods can be classified into two groups. In contrast to FM and CG-Poisson, the methods FFT, DCT and 
Sylvester, which cannot handle flexible domains, provide inaccurate reconstructions
which are not useful. The non-applicability of these methods is a considerable problem, since real-world input 
images for 3D reconstruction are typically in the centre of a photographed scene. This  
requires the flexibility to tackle non-rectangular domains, which is necessary as we have shown here
to get accurate (and efficient) reconstructions. 

As a conclusion of this experiment it can be predicted that all methods except FM and CG-Poisson are not realisable 
as an ideal, high-quality normal integrator for many future applications.

Let us note that we also have shown that FM and CG-Poisson have complementing properties and 
disadvantages -- the former is fast but inaccurate, and 
the latter is slow but accurate. At this point the combination of FM as initialisation, and a Krylov-based 
Poisson solver is clearly motivated as these should combine to a fast and accurate solver.

\subsection{Summary of the Evaluation}
In the previous experiments we tested different scenarios, which arise in real-world applications. 
It was found that boundary conditions and noisy data may have a strong effect on 
the 3D reconstructions. If rectangular domains can be considered the DCT method seems to be a realistic choice of a 
normal integrator followed by Sylvester and CG-Poisson. In fact the former one is unbeatably fast. 
However, the handling of non-rectangular domains, which is a practical issue in many 
industrial applications, cannot be underestimated. This important scenario leads to inaccuracies 
in the reconstructions of DCT and Sylvester. In this context  FM and CG-Poisson 
achieve better results. 

One can observe a certain lack in robustness w.r.t.\ noise of the FM integrator, especially
along directions not aligned with the grid structure, see for instance the results in~\cite{BQBD16}. 
This is 
because of the causality concept behind the FM scheme; errors that once appear are 
transported over the computational domain. This is not the case using Poisson reconstruction, 
which is a global approach and includes a regularising mechanism via the underlying least squares model.

Due to the possibly non-rectangular nature of the domains we aim to tackle, we cannot use fast Poisson solvers 
as e.g.\ in \cite{Simchony1990} to solve the arising discrete Poisson equations numerically. Instead, 
we explicitely constructed 
linear systems and solved them using the CG solver as often done 
by practitioners. Nevertheless, the unmodified CG-Poisson solver is still quite inefficient.

\section{Accelerating CG-Poisson}\label{PCG}

Let us now demonstrate the advantages of the proposed FM-PCG approach 
compared to other 
state-of-the-art
 methods. Thereby, we give a careful evaluation of
all the components of our novel algorithm.

\begin{table*}[!t]
    \caption{
Number of iterations and CPU time required to reach a $10^{-4}$ relative residual for the conjugate 
gradient algorithm, using the shifted Modified Incomplete Cholesky (MIC) preconditioner with different drop tolerances 
and different ``Phantom" sizes. The $10^{-3}$ drop tolerance is the one which provides the fastest results. 
Using a larger drop tolerance allows to reduce the number of required iterations, but the time used for computing the 
preconditioner dramatically increases. Note that we were unable to compute the preconditioner MIC($10^{-5})^*$ 
for the $4096^2$ dataset, because our 32GB of available memory were not sufficient.}
    \label{Tab:precond}
    \vspace {2.5mm}
    \setlength{\tabcolsep}{4.5pt}
    \centering\begin{tabular}{c|cc|cc|cc|cc}
        \toprule
        ~ & \multicolumn{2}{c}{No precond.} & \multicolumn{2}{c}{MIC($0)^*$} & \multicolumn{2}{c}{MIC($10^{-1})^*$} & \multicolumn{2}{c}{MIC($10^{-2})^*$}\\ 
        \hline
        Size & it. & CPU ($s$) & it. & CPU ($s$) & it. & CPU ($s$)& it. & CPU ($s$)\\
        \hline
        $64^2$ & 131 & 0.04 & 19 & \textbf{0.01} & 19 & \textbf{0.01} & 10 & \textbf{0.01}\\ 
        $128^2$ & 236 & 0.14 & 29 & 0.05 & 29 & 0.05 & 15 & \textbf{0.04}\\ 
        $256^2$ & 432 & 0.86 & 40 & 0.30 & 40 & 0.26 & 23 & 0.21 \\
        $512^2$ & 604 & 4.75 & 70 & 1.43 & 70 & 1.50 & 28 & \textbf{0.88} \\
        $1024^2$ & 1059 & 35.55 & 91 & 7.38 & 91 & 7.52 & 49 & 5.15 \\
        $2048^2$ & 1718 & 233.49 & 160 & 49.89 & 160 & 49.89 & 79 & 31.29 \\
        $4096^2$ & 2969 & 1577.81 & 247 & 290.93 & 247 & 290.54 & 134 & 196.1 \\
        \bottomrule
    \end{tabular}
    \vspace {4mm}
    \setlength{\tabcolsep}{4.5pt}
    \centering\begin{tabular}{c|cc|cc|cc}
        \toprule
        ~ & \multicolumn{2}{c}{MIC($10^{-3})^*$} & \multicolumn{2}{c}{MIC($10^{-4})^*$} & \multicolumn{2}{c}{MIC($10^{-5})^*$} \\
        \hline
        Size & it. & CPU ($s$)& it. & CPU ($s$)& it. & CPU ($s$)  \\
        \hline
        $64^2$ & 5 & \textbf{0.01} & \textbf{4} & \textbf{0.01} & \textbf{4} & 0.02\\
        $128^2$ & 9 & \textbf{0.04} & \textbf{6} & 0.06 & 7 & 0.12\\
        $256^2$ & 11 & \textbf{0.18} & \textbf{7} & 0.25 & 11 & 0.53 \\
        $512^2$ & 18 & 0.89 & \textbf{13} & 1.22 & 14 & 2.29 \\
        $1024^2$ & 30 & \textbf{5.02} & \textbf{19} & 5.96 & 24 & 11.83  \\
        $2048^2$ & 49 & \textbf{28.76} & \textbf{34} & 35.09 & 39 & 65.06 \\
        $4096^2$ & 80 & \textbf{171.44} & \textbf{41} & 173.5 & N/A & N/A \\        
        \bottomrule
    \end{tabular}
\end{table*}

\subsection{Preconditioned CG-Poisson}
In a first step we analyse the behaviour of the CG solver applying an additional preconditioner. 
This step is supposed to improve the condition number and convergence speed of an iterative solver 
in relation to the number $k$ of iterations and therefore the saving of computational time to reach 
the stopping criterion.

As examples of actual preconditioners, we checked\footnote{Due to the fact that pivot breakdowns for MIC($\tau$) are possible, we considered the shifted MIC version MIC($\tau,\alpha$). All methods are predefined functions in MATLAB.} IC($\tau$) and MIC($\tau,\alpha$) for the test dataset ``Phantom'' (see Figure \ref{Fig:precond}) for different input 
sizes. It was observed, that MIC($\tau,\alpha^*$) beates IC($\tau$) if we used $\alpha^*=10^{-3}$ for the global diagonal shift \footnote{This is an experimentally determined value.}.
In the following, we denoted for reasons of clarity MIC($\tau,\alpha^*$) as MIC($\tau)^*$.
The results of the MIC($0)^*$, without a fill-in strategy, are shown in Table \ref{Tab:precond} (third column) 
and illustrate its usefulness compared to the non-preconditioned CG-Poisson (second column). 
 By using MIC($0)^*$ we 
save lots of iterations and therefore we can reduce the time costs a lot: for example one can save 
around 2700 iterations and thus more than $1250~s$ for an image of size $4096\times 4096$.

\begin{figure}[!ht]
  \caption{The ``Phantom" image used in this experiment. Its gradient is unknown, hence we approximate it numerically 
by first-order forward differences. We used this dataset for comparing preconditioners, for different image sizes, 
from $64 \times 64$ to $4096 \times 4096$.}
  \label{Fig:precond}
  \vspace {2.5mm}  
\centering  \includegraphics[width = 0.33\linewidth]{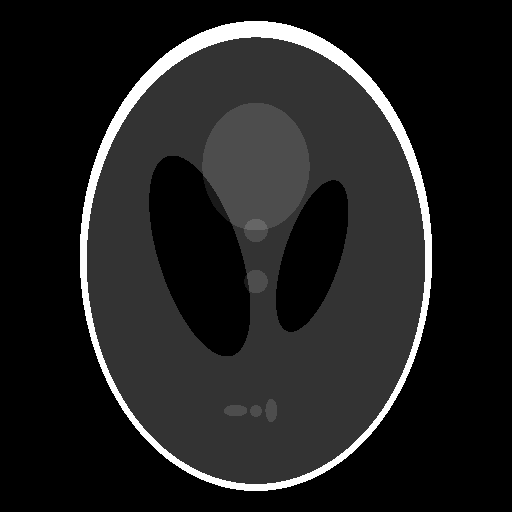}
\end{figure}

Now we want to show how useful a fill-in strategy can be. In the columns four to eight we tested different 
fill-in strategies from MIC($10^{-1})^*$ to MIC($10^{-5})^*$ 
\footnote{If the value $\tau$ in MIC($\tau)^*$ tends 
to zero then the preconditioned matrix is more dense -- having more non-zero elements. 
As a consequence, the preconditioner is better, however it costs more time to compute the preconditioner itself.}.
 A closer examination of Table \ref{Tab:precond} shows that 
MIC($10^{-3})^*$ provides the best
balance between the time to compute the 
preconditioner and the time
apply PCG. As an example, see again the image of size
 $4096\times 4096$, we can reduce the number of iterations from 247 to 80 and need around $170~s$ 
instead of $290~s$.

Thus, the application of a preconditioning, here shifted MIC, seems to be a useful solution to accelerate the CG-Poisson 
integrator, nevertheless it is not sufficient to be competitive to the common fast methods. 
However, we will see that with a proper initialisation, this standard preconditioner can already be 
considered as efficient.

\begin{table*}[!ht]
    \caption{Number of iterations and CPU time for applying the PCG algorithm, starting from the FM solution rather 
than from the trivial state. The indicated CPU time includes the time for computing the FM initialisation. Using 
FM as initial guess allows saving a lot of computations: 
the computation time to solve the $4096^2$ problem is reduced from $26$ $min$ (without FM initialisation 
nor preconditioning, see second column in Table~\ref{Tab:precond}, to $1$ $min$ (with FM initialisation and preconditioning, 
see column 6). }
    \label{Tab:init1}
    \vspace {2.5mm}
    \setlength{\tabcolsep}{4.5pt}
    \centering\begin{tabular}{c|cc|cc|cc|cc}
        \toprule
        ~ & \multicolumn{2}{c}{No precond.} & \multicolumn{2}{c}{MIC($0)^*$} & \multicolumn{2}{c}{MIC($10^{-1})^*$} & \multicolumn{2}{c}{MIC($10^{-2})^*$}\\ 
        \hline
        Size & it. & CPU ($s$) & it. & CPU ($s$) & it. & CPU ($s$)& it. & CPU ($s$)\\ 
        \hline
        $64^2$ & 119 & 0.05 & 16 & \textbf{0.02} & 16 & \textbf{0.02} & 8 & \textbf{0.02}\\ 
        $128^2$ & 210 & 0.17 & 25 & 0.10 & 25 & 0.10 & 13 & \textbf{0.09} \\  
        $256^2$ & 240 & 0.74 & 35 & 0.32 & 35 & 0.31 & 15 & 0.24 \\
        $512^2$ & 281 & 2.48 & 36 & 1.18 & 36 & 1.21 & 15 & \textbf{0.92} \\ 
        $1024^2$ & 316 & 12.80 & 40 & 5.19 & 40 & 5.20 & 18 & 4.06 \\  
        $2048^2$ & 339 & 55.44 & 45 & 23.14 & 45 & 23.28 & 19 & 13.08 \\  
        $4096^2$ & 349 & 232.95 & 46 & 98.70 & 46 & 99.30 & 19 & 76.34 \\ 
        \bottomrule
    \end{tabular}
    \vspace{4mm}
    \setlength{\tabcolsep}{4.5pt}
    \centering\begin{tabular}{c|cc|cc|cc}
        \toprule
        ~ & \multicolumn{2}{c}{MIC($10^{-3})^*$} & \multicolumn{2}{c}{MIC($10^{-4})^*$} & \multicolumn{2}{c}{MIC($10^{-5})^*$} \\
        \hline
        Size & it. & CPU ($s$)& it. & CPU ($s$)& it. & CPU ($s$)  \\
        \hline
        $64^2$ & 4 & \textbf{0.02} & \textbf{3} & \textbf{0.02} & \textbf{3} & 0.04  \\
        $128^2$ &  7 & \textbf{0.09}& \textbf{5} & 0.12 & 6 & 0.16  \\
        $256^2$ &  7 & \textbf{0.23}&   \textbf{5} & 0.30 & 6 & 0.54  \\
        $512^2$ &  9 & 0.93 & \textbf{5} & 1.07 & 6 & 2.22   \\
        $1024^2$ &  9 & \textbf{3.97}& \textbf{5} & 4.92 & 8 & 9.28   \\
        $2048^2$ &  9 & \textbf{17.12}& \textbf{5} & 21.47 & 8 & 39.91  \\
        $4096^2$ & 9 & \textbf{74.04} & \textbf{5} & 107.24 & N/A & N/A  \\        
        \bottomrule
    \end{tabular}
\end{table*}


\subsection{Appropriate Initialisation}
The suggested preconditioned CG-Poisson (PCG-Poisson) method is not widely known in computer vision, although this practicable method is surely not new and commonly used in numerical computing.
However, using an appropriate initialisation which shall decrease the number of iterations and reduce the run 
time costs, we propose a novel scheme for the 
surface normal integration (SNI) task.
Accordingly, our proposed method consists of two steps:
 in a first step the FM solution is computed 
in a fast and efficient way; after that, the Krylov-based technique with shifted modified incomplete Cholesky (MIC) is applied.

To show the effect of the new FM initialisation, the latter test for the ``Phantom" dataset is repeated 
and evaluated anew, see Table \ref{Tab:init1}. Starting from the FM solution, which needs comparatively short computation time 
see
 Table \ref{Tab:FM}) also for large images, 
gives a dramatic speed-up.

\begin{table}[!h]
    \caption{CPU time for the performing of FM on ``Phantom''  dataset  for different sizes.}
    \label{Tab:FM}
    \vspace {2.5mm}
    \setlength{\tabcolsep}{2.5pt}
    \centering\begin{tabular}{c|c|c|c|c|c|c|c}
          \toprule        
        Size & $64^2$ & $128^2$ & $256^2$ & $512^2$ & $1024^2$ & $2048^2$ & $4096^2$\\
         \hline
        CPU ($s$) & $<$ 0.01 & 0.03 & 0.08 & 0.28 & 1.20  & 5.14 & 21.79\\
\bottomrule
    \end{tabular}
\end{table}

A closer look at Tables \ref{Tab:precond} and \ref{Tab:init1} shows a significant difference, even 
without a fill-in strategy (compare both third columns). At first, let us consider the case without preconditioner: starting with 
the trivial solution leads to a constant increase of iterations (factor around 1.7) with simultaneous 
increasing of the image size. In contrast the number of iterations increases very slowly by using FM 
initialisation. The effect of this phenomenon is a notable, strong time cost reduction for large data:  
for 
 $512\times 512$ images,
 we can save more than $2~s$ 
(from $4.75~s$ to $2.48~s$),
and for 
 $4096\times 4096$ images 
the time can be reduced from $1578~s$ to $233~s$. 

The case of an additional preconditioning leads to similar results. Testing anew MIC($\tau)^*$ with MIC($10^{-1})^*$ to 
MIC($10^{-5})^*$ shows once more that MIC($10^{-3})^*$ provides the best results, 
see Table \ref{Tab:init1}. 
Using a FM initialisation reduces the required iterations to reach the stopping criterion a lot and therefore the 
combination of FM and shifted MIC leads to fast reconstructions. In case of an image of size $4096\times 4096$ the 
novel approach, including the FM performing in $21.79~s$ (see Table \ref{Tab:FM}), saves around $100~s$ 
(from $171~s$ to $74~s$) and 71 iterations compared with the trivial initialisation and MIC($10^{-3})^*$. 



Finally, using the novel approach instead of the standard CG-Poisson solver leads to an incredible speed-up, see Table \ref{Tab:init2}. 
Without considering the computation of the FM initialisation, the construction of the system and the preconditioner, 
the pure time to solve the system can be amazingly reduced from $1551.96~s$ to $18.90~s$. The findings of this experiment show impressively that choosing FM as initialisation accelerates the method a lot when it comes to 
standard preconditioners like (shifted modified) incomplete Cholesky. Thus, we believe that our novel FM-PCG method with shifted MIC preconditioning is a relevant contribution to the field of fast and accurate surface normal integrators.
\begin{table}[!h]
    \caption{ 
Repartition of CPU time between system construction, preconditioning and system resolution, for the $4096^2$ example. Knowing that the system and the 
preconditioner can often be pre-computed, this makes even more obvious the gain one can expect by choosing an appropriate 
initialisation such as by the FM result. CG refers to the resolution of the system by conjugate gradient, and +CG to the 
accelerated resolution by choosing the FM initialisation (not including the $ 21.79$ $s$ required for FM).}
   \label{Tab:init2}
   \vspace {2.5mm}
    \setlength{\tabcolsep}{4.5pt}
    \centering\begin{tabular}{ccc|cc}
        \toprule        
        ~ & Syst. constr. & Precond. & CG & +CG \\
        \hline
        No precond. & 25.85 & 0 & 1551.96 & \textbf{185.31}  \\
        MIC($10^{-3})^*$ & 25.85 & 7.50 & 138.09 & \textbf{18.90} \\                       
        \bottomrule
    \end{tabular}
\end{table}


\subsection{Evaluation of the FM-PCG Solver}
To clarify the strength of our proposed FM-PCG solver against the 
standard solvers FFT, DCT and the ``Sylvester" method of Harker and O'Leary,  
we evaluate the reconstructions of 
the datasets ``Phantom'', ``Lena'', ``Peaks'' and ``Vase'' on rectangular and non-rectangular domains with the associated values of MSE.
At first we examine the ``Phantom'', ``Lena' and ``Peaks'' datasets on a rectangular domain 
in Tables \ref{tab:results_Phantom}, \ref{tab:results_Lena} and \ref{tab:results_Peaks}. All examples contain the 
”natural“ boundary equation, moreover ``Phantom'' and ``Lena'' have sharp gradients and are more realistic.

It should be clear that FFT and DCT are the fastest methods, however the quality of FFT is not adequate 
and the results are unusable. Furthermore, it can be seen that the FM-PCG solver is the best integrator 
with respect to sharp gradients (see Table~\ref{tab:results_Lena}).

\begin{table}[!h]
    \caption{Results on the ``Phantom" dataset ($1024 \times 1024$).} 
  \label{tab:results_Phantom}    
    \vspace {2.5mm}
    \setlength{\tabcolsep}{4.5pt}
    \centering\begin{tabular}{ccc}
        \toprule
        Method & MSE ($px$) & CPU ($s$) \\ \hline
        FFT \cite{Frankot1988} & 138.6 & \textbf{0.06} \\
        DCT \cite{Simchony1990} & \textbf{127.31} & 0.13 \\
         FM \cite{Galliani2012} & 163.13 & 1.20  \\
        Sylvester \cite{Harker2015} & 169.41 & 5.78 \\
	FM-PCG & 127.89 & 4.23\\
        \bottomrule
    \end{tabular}
\end{table}

\begin{table}[!h]
    \caption{Results on the ``Lena" dataset ($512 \times 512$).} 
  \label{tab:results_Lena}    
    \vspace {2.5mm}
    \setlength{\tabcolsep}{4.5pt}
    \centering\begin{tabular}{ccc}
        \toprule
        Method & MSE ($px$) & CPU ($s$) \\ \hline
        FFT \cite{Frankot1988} & 402.37 & \textbf{0.02} \\
        DCT \cite{Simchony1990} & 132.08 & 0.03 \\
        FM \cite{Galliani2012} & 509.15 & 0.28   \\
        Sylvester \cite{Harker2015} & 113.92 & 0.71 \\
	FM-PCG & \textbf{94.07} & 1.24\\
        \bottomrule
    \end{tabular}
\end{table}

\begin{table}[!h]
    \caption{Results on the ``Peaks" dataset ($128 \times 128$).} 
  \label{tab:results_Peaks}    
    \vspace {2.5mm}
    \setlength{\tabcolsep}{4.5pt}
    \centering\begin{tabular}{ccc}
        \toprule
        Method & MSE ($px$) & CPU ($s$) \\ \hline
        FFT \cite{Frankot1988} & 7.19 & $\mathbf{<}$ \textbf{0.01} \\
        DCT \cite{Simchony1990} & 0.09 & $\mathbf{<}$ \textbf{0.01} \\
        FM \cite{Galliani2012} & 0.8  & 0.03  \\
        Sylvester \cite{Harker2015} & \textbf{0.01} & 0.07 \\
	FM-PCG & 0.02 & 0.07\\
        \bottomrule
    \end{tabular}
\end{table}


Finally, a user would employ the method with the best speed$/$quality balance, which is on rectangular domains 
probably DCT followed by Sylvester and our proposed FM-PCG solver. However, as already 
mentioned, simple rectangular domains are quite unrealistic in many applications in science and industry. Hence, we analyse in Tables~\ref{tab:results_Phantom_masked}, \ref{tab:results_Lena_masked}, \ref{tab:results_Peaks_masked} and \ref{tab:results_Vase_masked} the results
on flexible domains,  with the accompanying computational 
domains depicted in Figure \ref{fig:mask_evaluation}.

\begin{figure}[!ht]
  \caption{Masks for the evaluation of ``Phantom'', ``Lena'', ``Peaks''  and ``Vase'' datasets. It should be noted that FM and CG-Poisson work only on the part $\Omega$, represented by the white regions. By contrast, FFT, DCT and Sylvester work on the whole rectangular domain. (a) One and the same synthetic mask for ``Phantom'', ``Lena'' and ``Peaks'' datasets. (b) Expedient mask for the ``Vase'' dataset.}
  \label{fig:mask_evaluation}  
  \begin{center}
  \setlength{\tabcolsep}{0.7cm}
    \begin{tabular}{cc}
      \includegraphics[width = 0.25\linewidth]{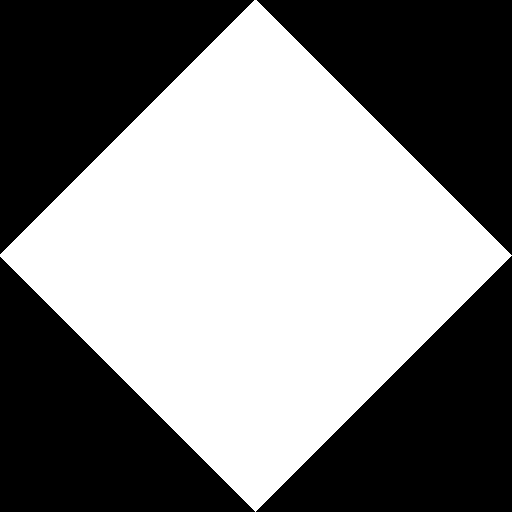} &
      \includegraphics[width = 0.25\linewidth]{mask_Vase.png} \\
      (a) & (b)      
    \end{tabular} 
  \end{center}
\end{figure}

\begin{table}[!ht]
    \caption{Results on the ``Phantom" dataset on the non-rectangular domain shown in Figure~\ref{fig:mask_evaluation}-a.} 
  \label{tab:results_Phantom_masked}    
    \vspace {2.5mm}
    \setlength{\tabcolsep}{4.5pt}
    \centering\begin{tabular}{ccc}
        \toprule
        Method & MSE ($px$)  & CPU ($s$) \\ \hline
        FFT \cite{Frankot1988} & 351.72 & \textbf{0.06} \\
        DCT \cite{Simchony1990} & 309.72 & 0.14 \\
         FM \cite{Galliani2012} & 162.43 & 0.70  \\
        Sylvester \cite{Harker2015} & 348.98 & 5.54 \\
	FM-PCG & \textbf{131.02} & 2.06\\
        \bottomrule
    \end{tabular}
\end{table}


\begin{table}[!ht]
    \caption{Results on the ``Lena" dataset on the non-rectangular domain shown in Figure~ \ref{fig:mask_evaluation}-a.} 
  \label{tab:results_Lena_masked}    
    \vspace {2.5mm}
    \setlength{\tabcolsep}{4.5pt}
    \centering\begin{tabular}{ccc}
        \toprule
        Method & MSE ($px$) & CPU ($s$) \\ \hline
        FFT \cite{Frankot1988} & 199.59 & \textbf{0.01} \\
        DCT \cite{Simchony1990} & 149.00 & 0.03 \\
        FM \cite{Galliani2012} & 444.02 & 0.19   \\
        Sylvester \cite{Harker2015} & 175.82 & 0.70 \\
	FM-PCG& \textbf{123.64} & 0.65\\
        \bottomrule
    \end{tabular}
\end{table}


\begin{table}[!ht]
    \caption{Results on the ``Peaks" dataset on the non-rectangular domain shown in Figure~ \ref{fig:mask_evaluation}-a.} 
  \label{tab:results_Peaks_masked}    
    \vspace {2.5mm}
    \setlength{\tabcolsep}{4.5pt}
    \centering\begin{tabular}{ccc}
        \toprule
        Method & MSE ($px$) & CPU ($s$) \\ \hline
        FFT \cite{Frankot1988} & 15.69 & $\mathbf{<}$ \textbf{0.01} \\
        DCT \cite{Simchony1990} & 7.23 & $\mathbf{<}$ \textbf{0.01} \\
         FM \cite{Galliani2012} & 0.86  & 0.01  \\
        Sylvester \cite{Harker2015} & 7.20 & 0.06 \\
	FM-PCG & \textbf{0.03} & 0.03 \\
        \bottomrule
    \end{tabular}
\end{table}

\begin{table}[!ht]
    \caption{Results on the ``Vase" dataset on the non-rectangular domain shown in Figure~ \ref{fig:mask_evaluation}-b.} 
  \label{tab:results_Vase_masked}    
    \vspace {2.5mm}
    \setlength{\tabcolsep}{4.5pt}
    \centering\begin{tabular}{ccc}
        \toprule
        Method & MSE ($px$) & CPU ($s$) \\ \hline
        FFT \cite{Frankot1988} & 5.71 & \textbf{0.01} \\
        DCT \cite{Simchony1990} & 5.69 & 0.02 \\
         FM \cite{Galliani2012} & 0.71  & 0.06  \\
        Sylvester \cite{Harker2015} & 5.99 & 0.38 \\
	FM-PCG& \textbf{0.03} & 0.14 \\
        \bottomrule
    \end{tabular}
\end{table}

All experiments show the expected behaviour of the employed methods. The FM-PCG solution has by far the 
best quality. It is even faster than Sylvester.
An assessment in relation to the best
performance in speed$/$quality balance is not quite easy and is depending on the exact application. 
Is the speed of secondary importance then the best choice is FM-PCG, otherwise DCT.

\subsection{Real-world Photometric Stereo Data}
The previous examples have been rather simple. For this reason we consider a more 
realistic real-world application in photometric stereo, which definitely contains 
noisy data. We used the ``Scholar" dataset\footnote{\url{http://vision.seas.harvard.edu/qsfs/Data.html}}, 
which consists of 20 images of a Lambertian surface, taken from the same angle of view but 
under 20 known, non-coplanar, lightings (see Figure~\ref{fig:results_scholar}). 

\begin{figure}[!h]
\caption{Application to photometric stereo (PS). (a-c) Three images (among 20), of size $1070 \times 1070$, 
acquired from the same point of view but under different lightings. After estimating the surface normals by PS, 
we integrated them by (d) FM, before (e) refining this initial guess by PCG iterations. The full integration 
process required a few seconds. (f-g) MSE (in pixels) 
 on the reprojected images, computed from the 
surface estimated (f) by FM  and (g) FM-PCG . (blue is $0$, and red is $>1000$). Due to the local nature of FM, radial propagation of 
errors is visible. After corrections by CG, such artefacts are eliminated. Remaining bias is due to shadows. 
These results are experimentally compared with existing methods in Table~\ref{tab:results_scholar}.}
  \label{fig:results_scholar}  
  \begin{center}
    \begin{tabular}{ccc}
      \includegraphics[width = 0.29\linewidth]{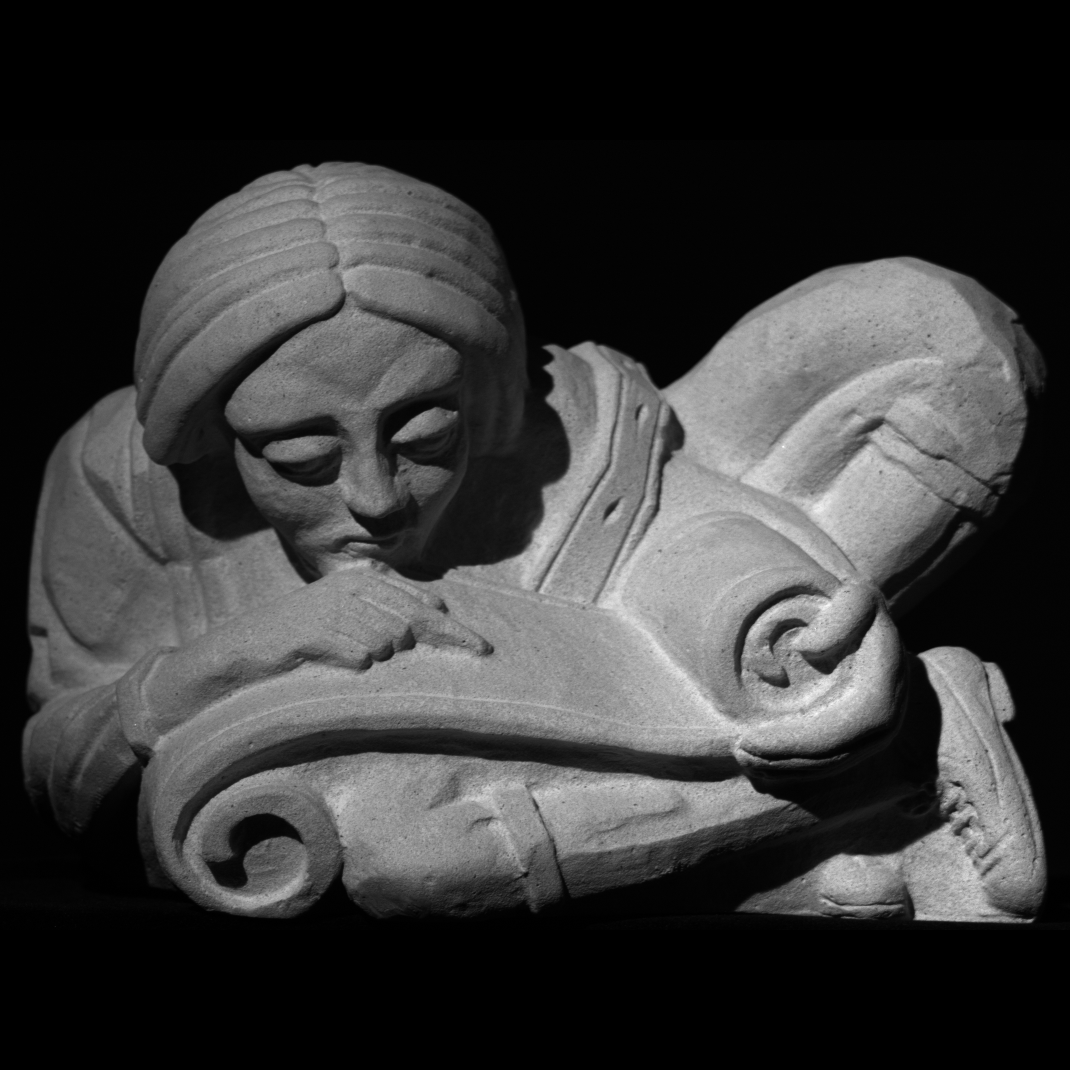} &
      \includegraphics[width = 0.29\linewidth]{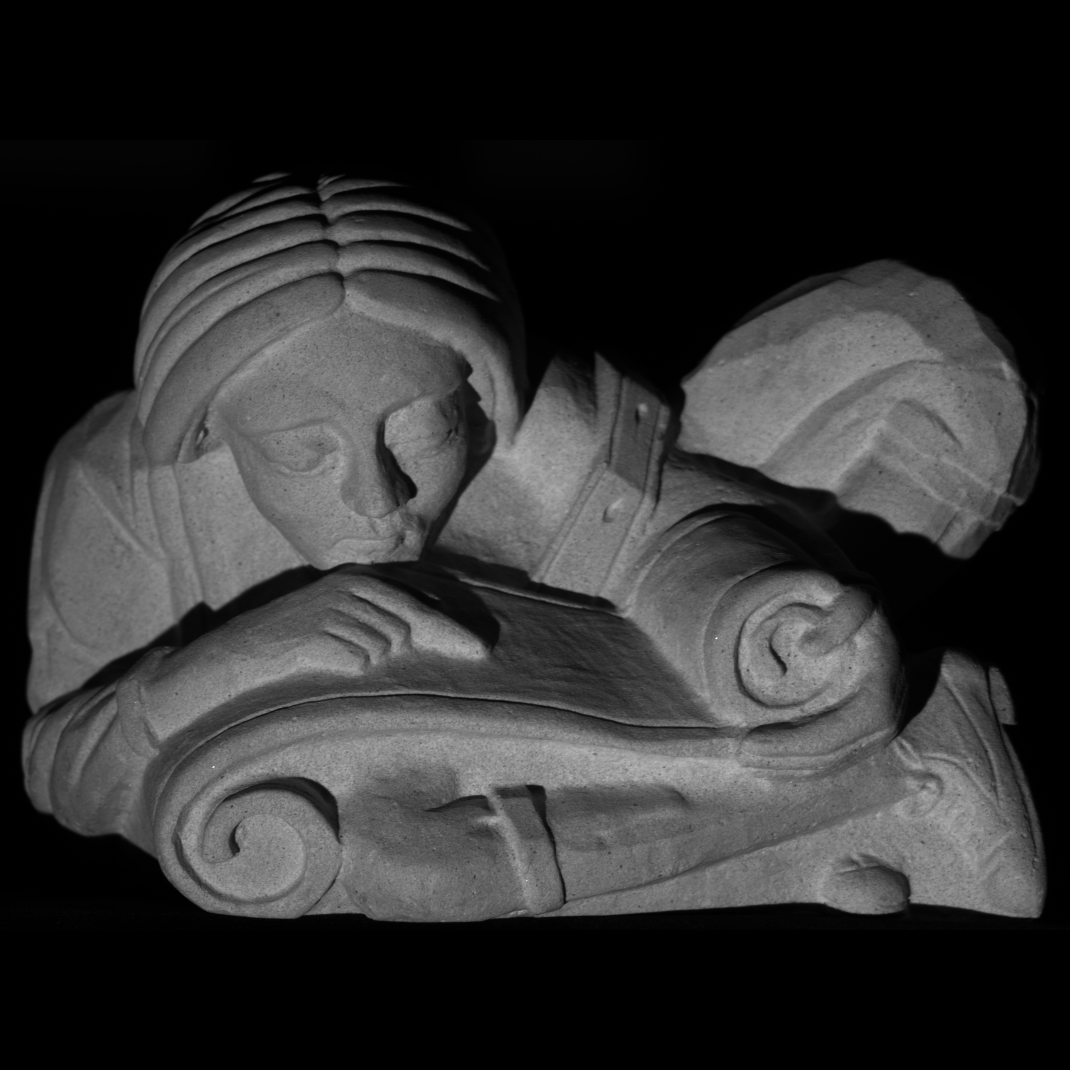} &
      \includegraphics[width = 0.29\linewidth]{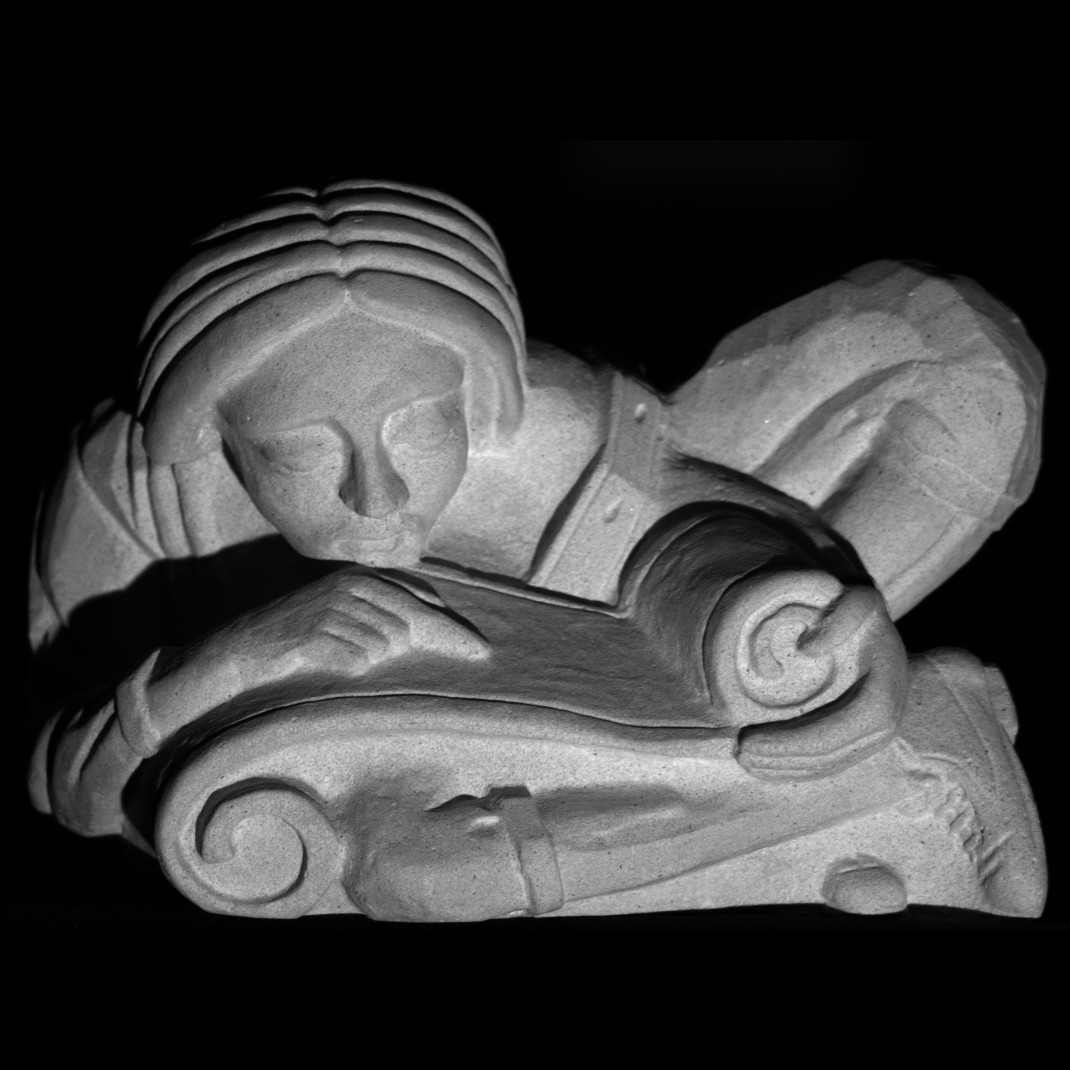} \\
    (a) & (b) & (c)   
        \vspace {2.5mm}
    \end{tabular}

    \begin{tabular}{cc}
      \includegraphics[width = 0.43\linewidth,trim=0 0 0 4cm,clip=true]{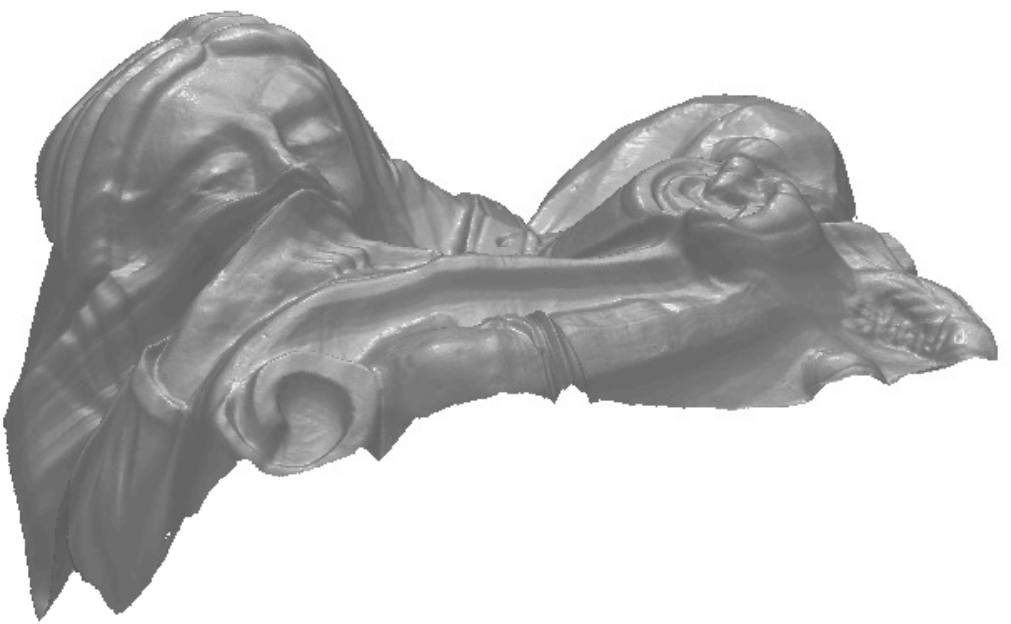} &
      \includegraphics[width = 0.43\linewidth,trim=0 0 0 4cm,clip=true]{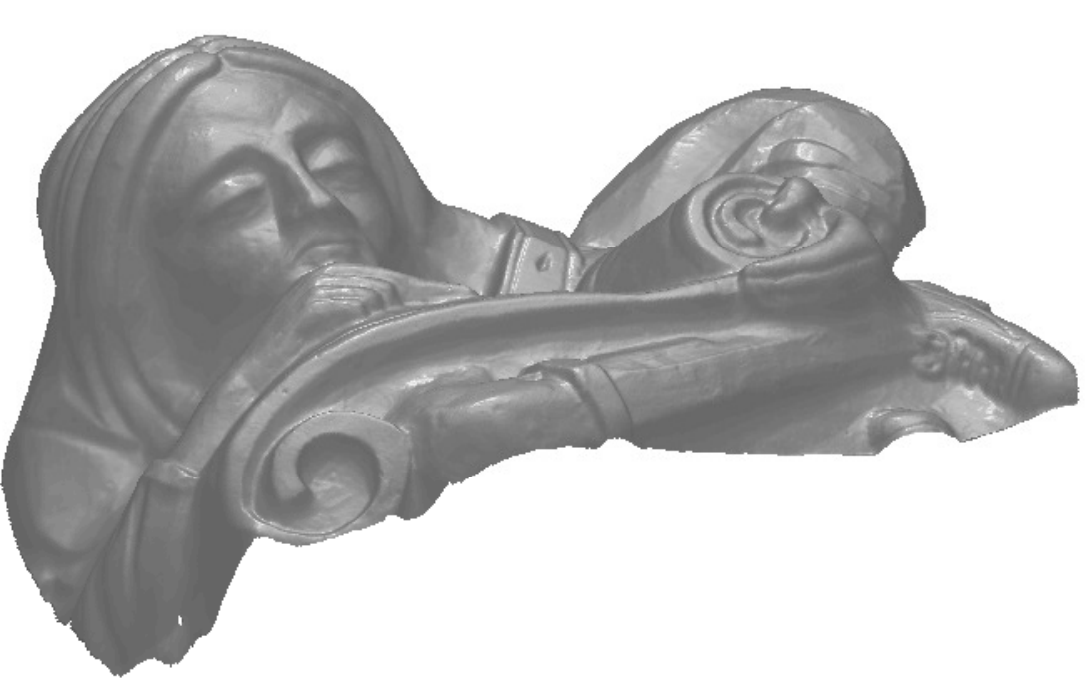} \\
      (d) & (e) \\[5pt]
      \includegraphics[width = 0.38\linewidth]{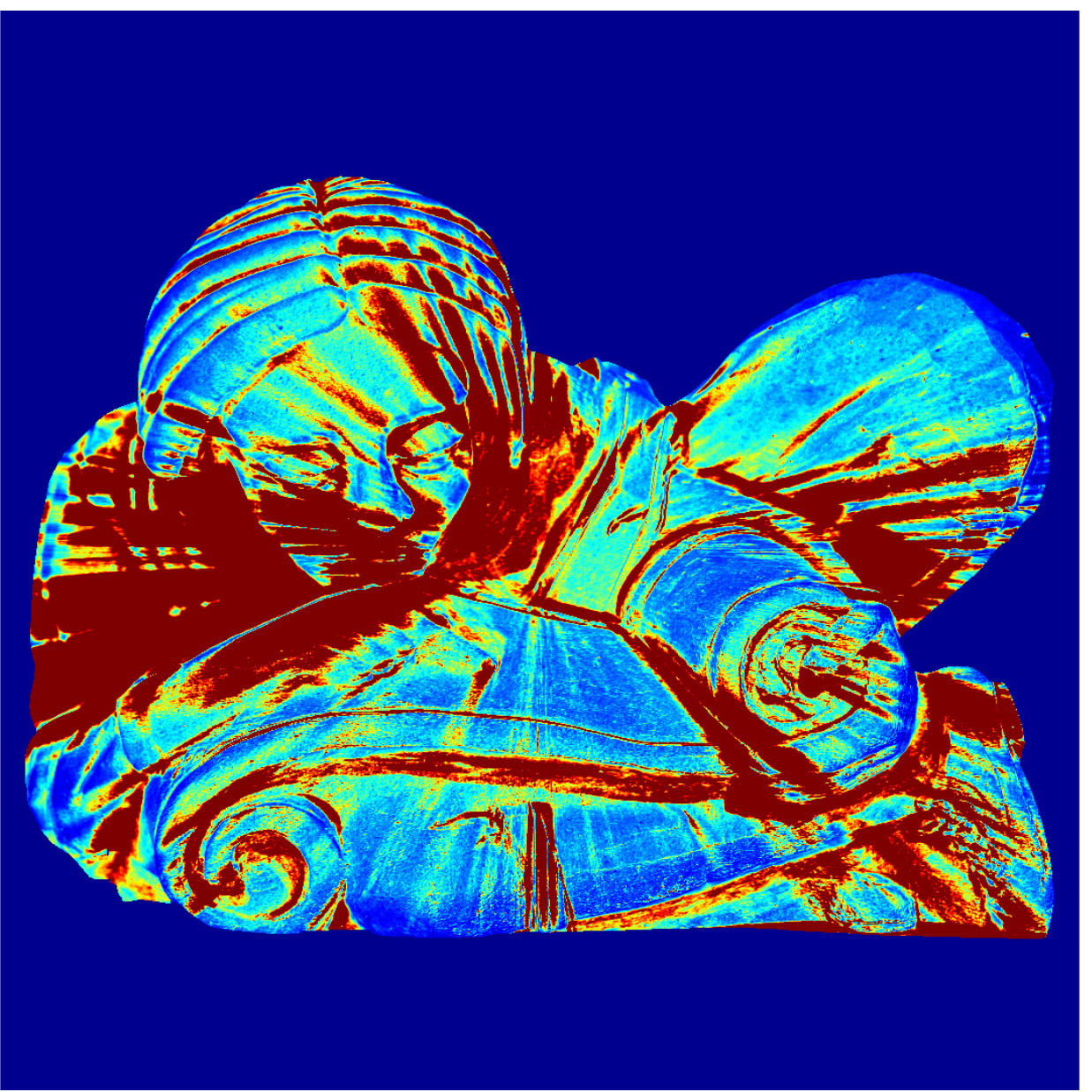} &
      \includegraphics[width = 0.38\linewidth]{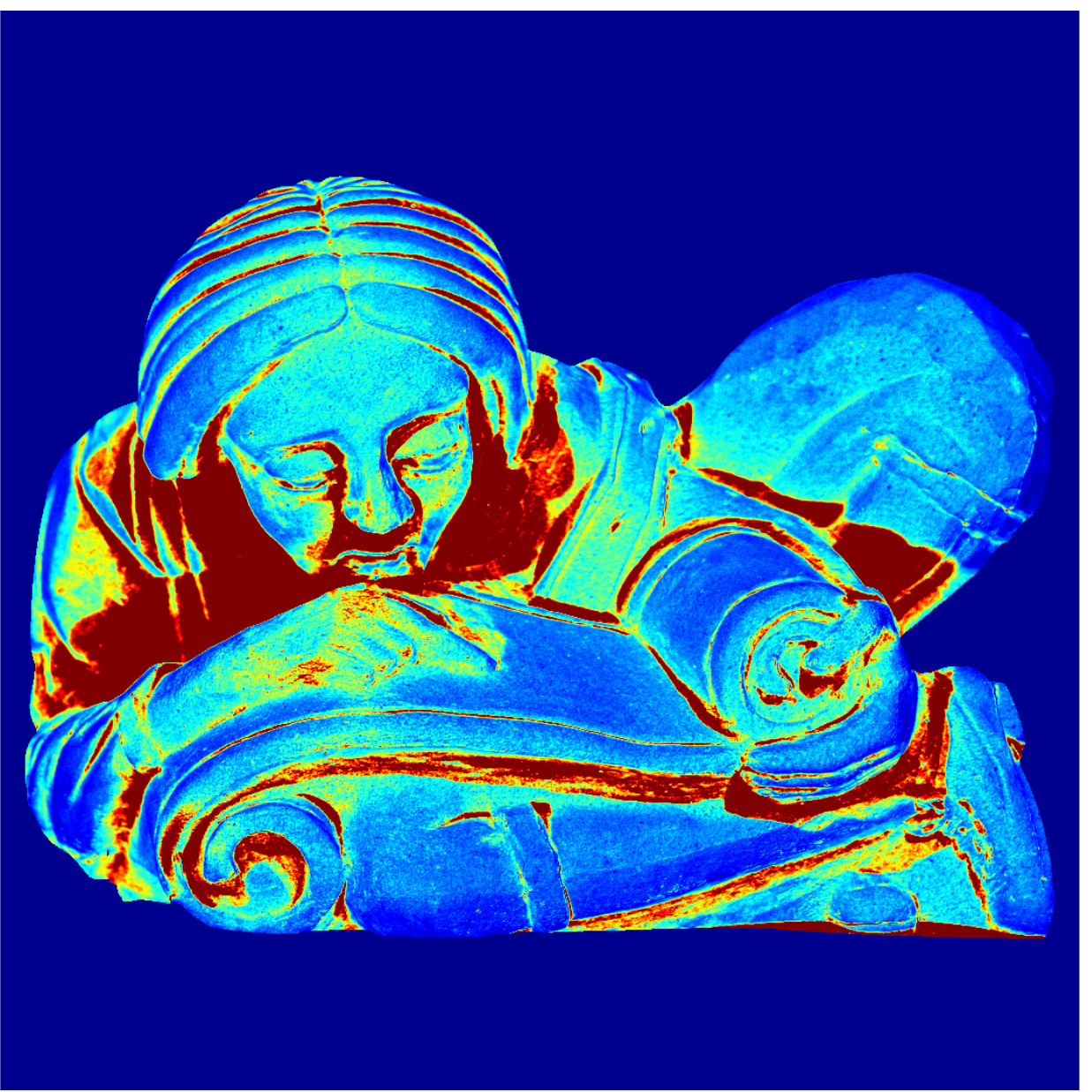} \\
      (f) & (g)      
    \end{tabular} 
  \end{center}
\end{figure}

\begin{table}[!h]
    \caption{Results on the PS dataset. Our method (initialisation by FM, then refinement by PCG from this 
initial guess) provides the most accurate results. We show the CPU time, as well as the 
mean MSE and SSIM on the 20 reprojected images.}
    \label{tab:results_scholar}
    \vspace {2.5mm}
    \setlength{\tabcolsep}{4.5pt}
    \centering\begin{tabular}{cccc}
        \toprule
        Method & MSE ($px$) & SSIM & CPU ($s$) \\ \hline
        FFT \cite{Frankot1988} & 365.43 & 0.86 & \textbf{0.09} \\
        DCT \cite{Simchony1990} & 330.55 & 0.87  & 0.15 \\        
        FM \cite{Galliani2012} & 582.65 & 0.78 & 0.45 \\
        Sylvester \cite{Harker2015} & 377.68 & 0.74 & 5.81 \\
        FM-PCG & \textbf{286.69} & \textbf{0.88} & 6.25 \\        
        \bottomrule
    \end{tabular}
\end{table}

The normals and the albedo were calculated using the classical photometric stereo approach 
by Woodham~\cite{Woodham1980}. Then, we integrated the normals using the different solvers. 
Eventually, we a posteriori recomputed the normals through finite differences from the 
recovered depth map, before ``reprojecting" the images using the estimated shape and albedo. 
By comparing the initial images with the reprojected ones, we obtain two criteria (MSE and SSIM) 
for evaluating the methods on each image. The results shown in Table \ref{tab:results_scholar} are the mean of the 20 corresponding values. 

Once again FM-PCG is the most accurate integrator and is as fast as Sylvester. Nevertheless, the 
fast computational time of DCT is unattainable.

\subsection{Handling Outliers}

Let us now
consider the case of standard photometric stereo 
applied to surfaces whose reflectance incorporates an additive off-Lambertian component (specularities). 
As can be seen from Fig.~\ref{fig:results_owl} 
and Table~\ref{tab:results_owl}
all the integration methods we consider here are
 by construction highly sensitive to 
the presence of outliers.

\begin{table}[!t]
    \caption{Results on the specular PS dataset (see Fig.~\ref{fig:results_owl}). All methods present a similar systematic  
bias due to the outliers located on the specular points.}
    \label{tab:results_owl}
    \vspace {2.5mm}
    \setlength{\tabcolsep}{4.5pt}
    \centering\begin{tabular}{cccc}
        \toprule
        Method & MSE ($px$) & SSIM & CPU ($s$) \\ \hline
        FFT \cite{Frankot1988} & 66.68 & 0.92 & $\mathbf{<}$ \textbf{0.01} \\
        DCT \cite{Simchony1990} & 46.16 & 0.95  & 0.01 \\        
        FM \cite{Galliani2012} & 94.25 & 0.90 & 0.09 \\
        Sylvester \cite{Harker2015} & 928.69 & 0.55 & 0.30 \\
        FM-PCG & \textbf{40.48} & \textbf{0.96} & 0.24 \\        
        \bottomrule
    \end{tabular}
\end{table} 

\begin{figure}[!h]
  \caption{(a-c) Three (out of 12) real-world images, of size $320 \times 320$, of a photometric 
stereo dataset. The eyes of the owl are highly specular. This induces a bias in the 
reconstructions, as shown in the reconstructions using (d) FM or (e) the proposed FM-PCG integrator.
(f-g) The corresponding MSE on the reprojected images shows that the bias is very 
localized (blue is $0$, and red is $>1000$).}
\label{fig:results_owl}  
  \vspace {2.5mm}    
  \begin{center}
    \begin{tabular}{ccc}
      \includegraphics[width = 0.29\linewidth]{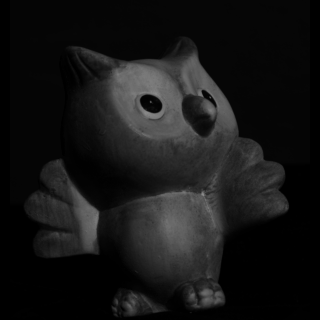} &
      \includegraphics[width = 0.29\linewidth]{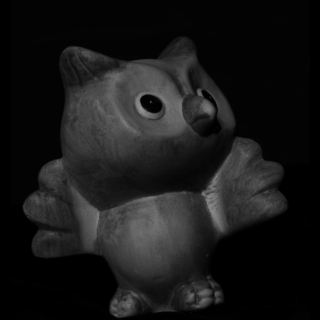} &
      \includegraphics[width = 0.29\linewidth]{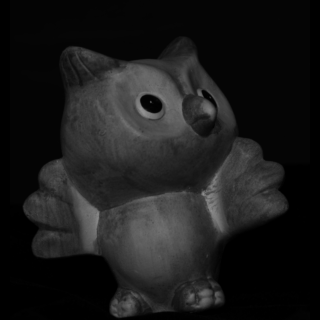} \\
    (a) & (b) & (c)    
    \end{tabular}
    \begin{tabular}{cc}
      \includegraphics[width = 0.46\linewidth,trim=0 0 0 2.5cm,clip=true]{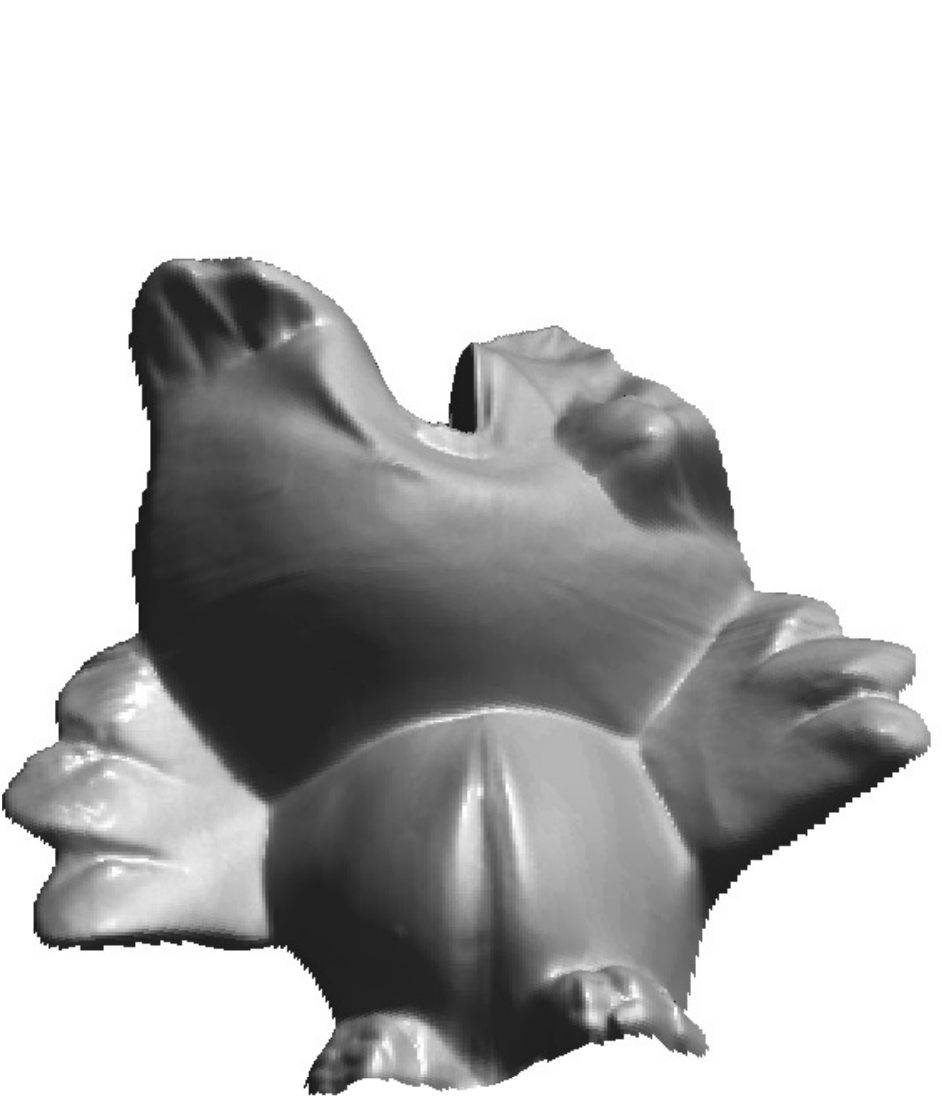} &
      \includegraphics[width = 0.46\linewidth,trim=0 0 0 2.5cm,clip=true]{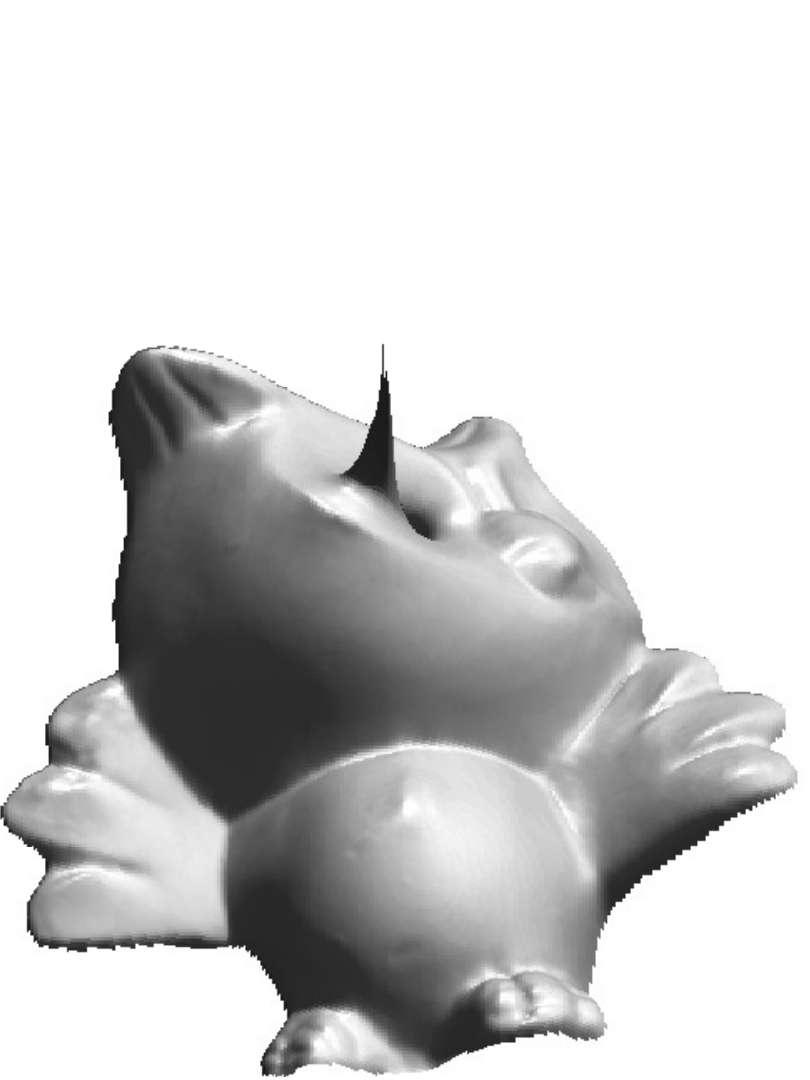} \\
      (d) & (e) \\[5pt]
      \includegraphics[width = 0.46\linewidth]{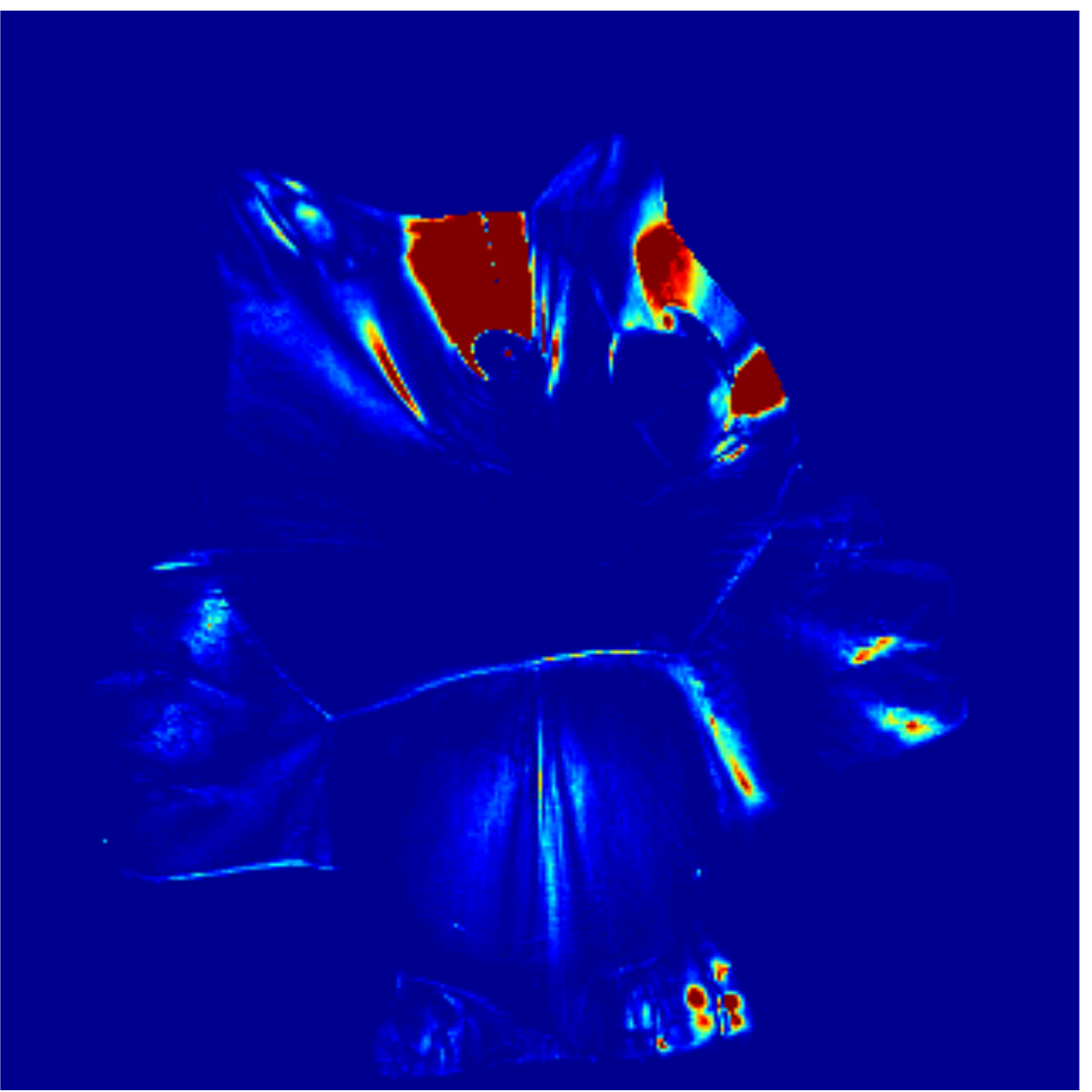} &
      \includegraphics[width = 0.46\linewidth]{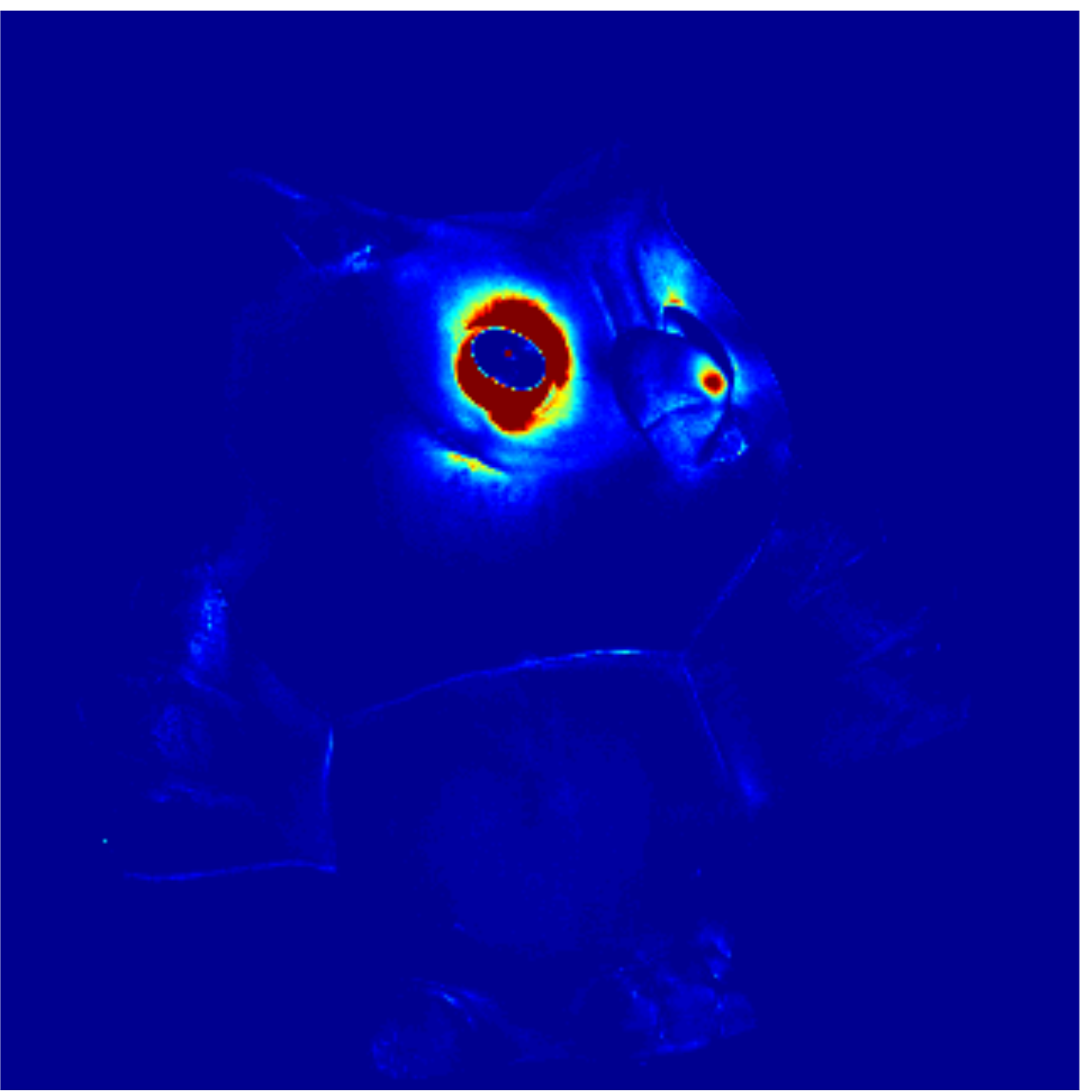} \\
      (f) & (g)      
    \end{tabular} 
  \end{center}
\end{figure}

In order to handle such outliers, we modify the Poisson integration framework according to the model~\eqref{poisson-outlier}. As we already pointed out, all the methods relying on the Poisson equation can be
adapted to this model.
Therefore, we can employ here by construction
the FFT~\cite{Frankot1988}, the DCT~\cite{Simchony1990}, and our new FM-PCG method.
We found that using these modified inputs for the other SNI methods, such 
as FM~\cite{Galliani2012} and Sylvester~\cite{Harker2015}, also yields improved results. 
Hence, our improved model can be considered as a \emph{generic} improvement for use with 
existing SNI methods, enforcing robustness to outliers. 
This is illustrated in Figure~\ref{fig:results_owlmod} and Table~\ref{tab:results_owlmod}.


\begin{figure}[!h]
    \caption{Result of the (a) improved FM and (b) 
    improved FM-PCG methods introducing a smoothness 
constraint on the outliers location. The corresponding (c) and (d) MSE maps show that errors due to the outliers are much reduced. }
\label{fig:results_owlmod}  
  \vspace {2.5mm}    
    \begin{center}
    \begin{tabular}{cc}
      \includegraphics[width = 0.45\linewidth,trim=0 0 0 2.5cm,clip=true]{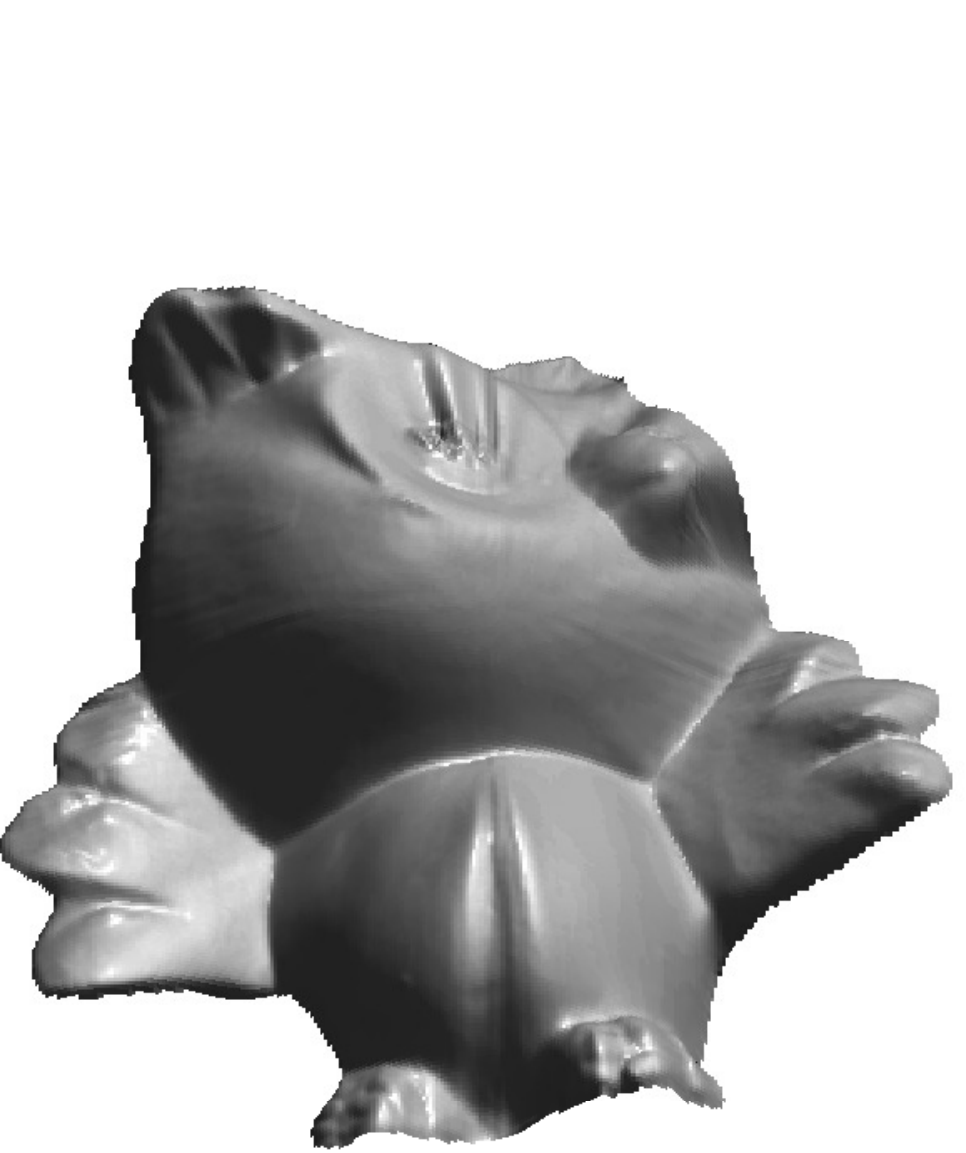} &
      \includegraphics[width = 0.45\linewidth,trim=0 0 0 2.5cm,clip=true]{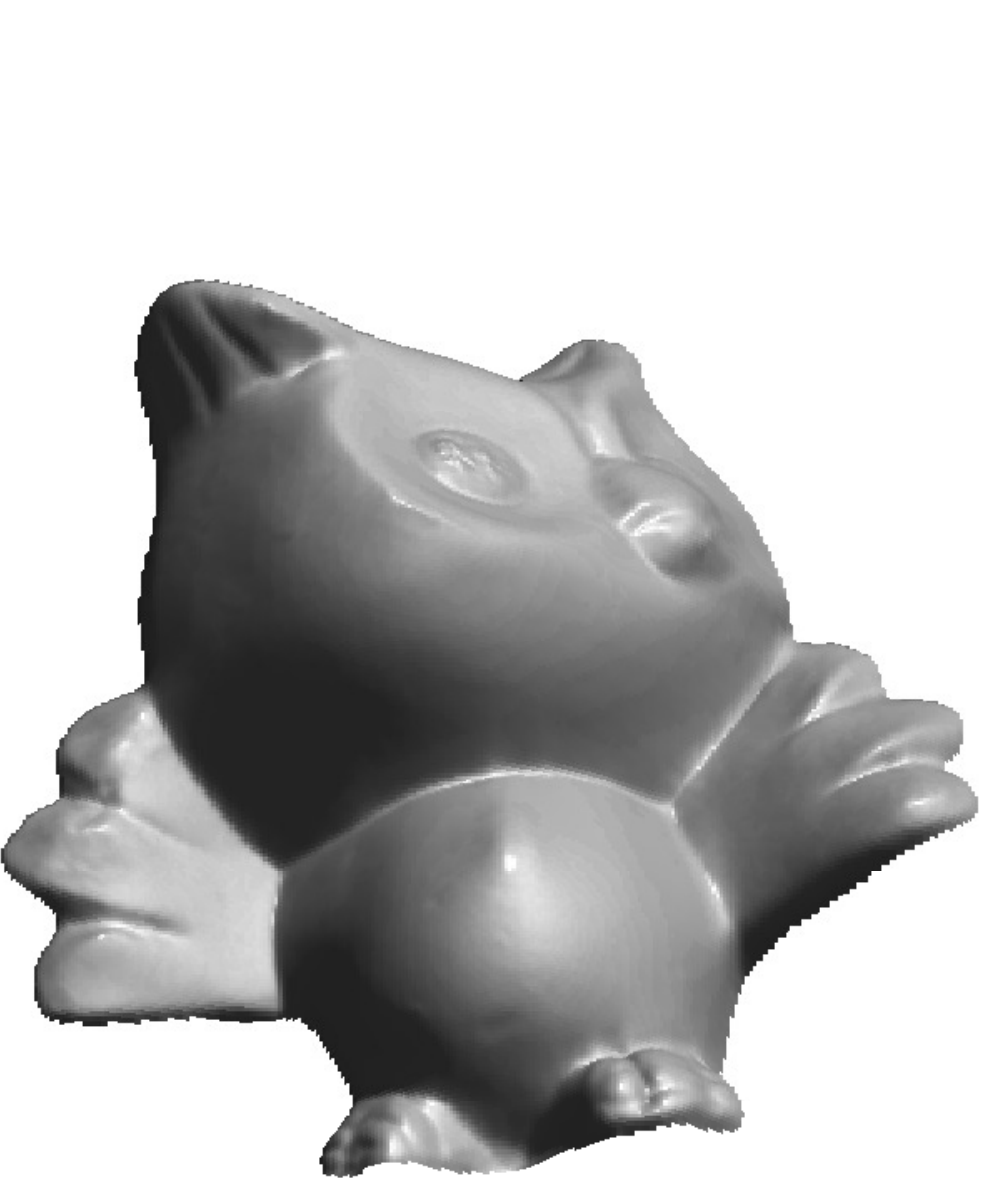}\\
    (a) & (b)  \\[5pt]
     \includegraphics[width = 0.46\linewidth]{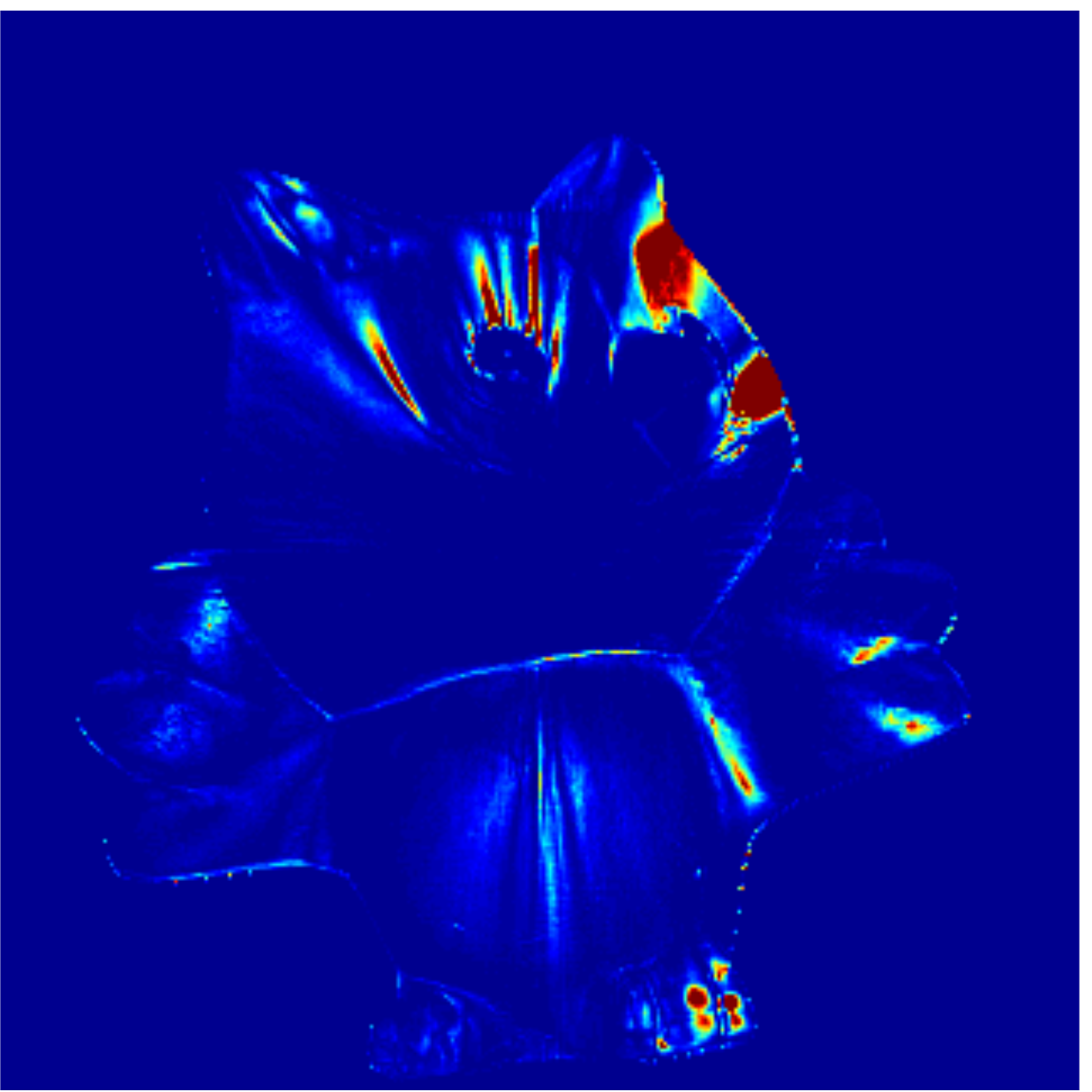} &
      \includegraphics[width = 0.46\linewidth]{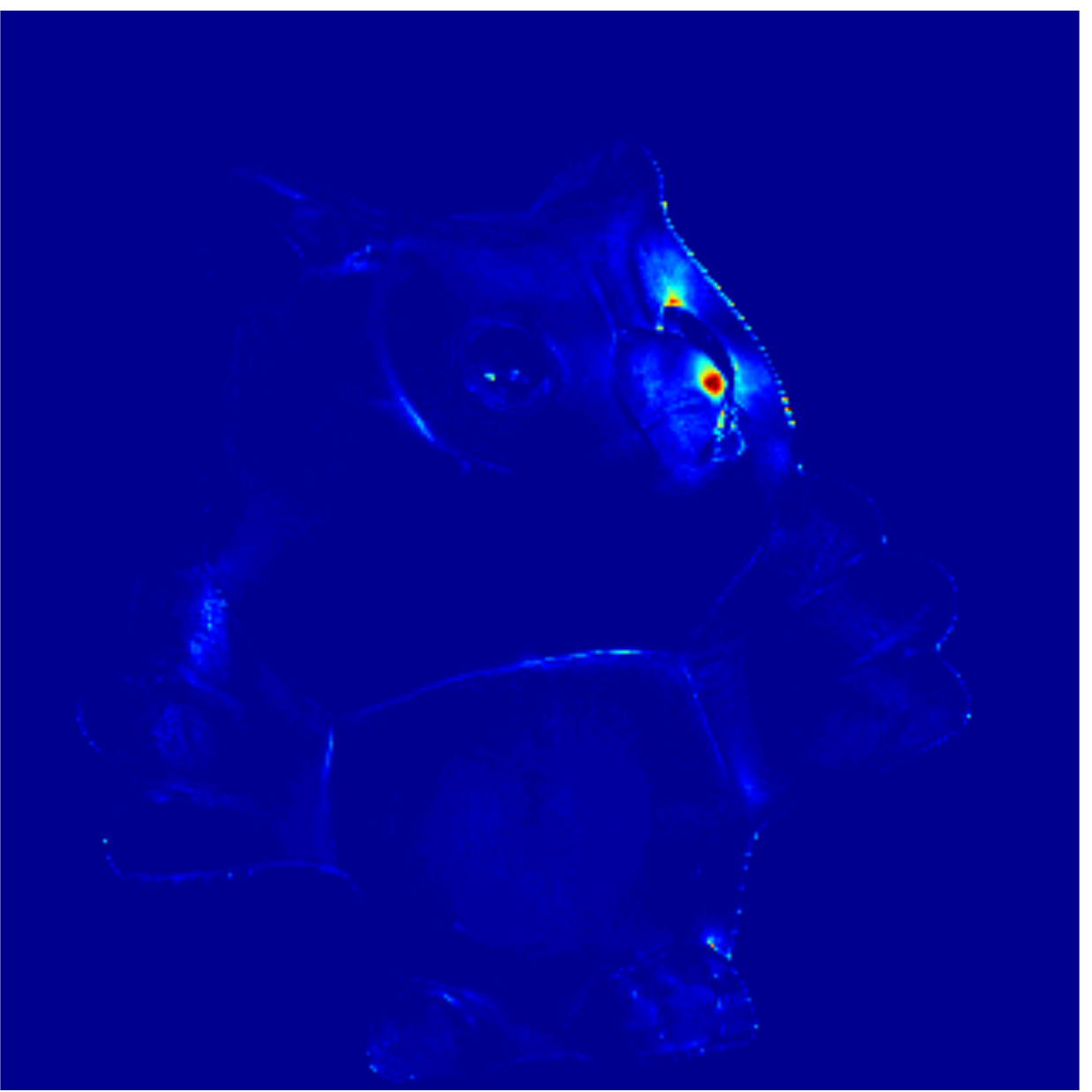} \\
      (c) & (d)      
    \end{tabular}
    \end{center}
\end{figure}

\begin{table}[!t]
    \caption{Results of the improved methods on the same data set as in Table~\ref{tab:results_owl}. 
All MSE are significantly reduced.}
    \label{tab:results_owlmod}
    \vspace {2.5mm}
    \setlength{\tabcolsep}{4.5pt}
    \centering\begin{tabular}{cccc}
        \toprule
        Method & MSE ($px$) & SSIM & CPU ($s$) \\ \hline
        FFT \cite{Frankot1988} & 17.26 & 0.95 & $\mathbf{<}$ \textbf{0.01} \\
        DCT \cite{Simchony1990} & 14.79 & \textbf{0.96}  & 0.01 \\        
        FM \cite{Galliani2012} & 34.47 & 0.91 & 0.09 \\
        Sylvester \cite{Harker2015} & 21.93 & 0.88 & 0.31 \\
        FM-PCG & \textbf{10.41} & \textbf{0.96} & 0.27 \\        
        \bottomrule
    \end{tabular}
\end{table}

\section{Conclusion and Perspectives}

We demonstrated the properties of the proposed FM-PCG surface
normal integrator. It combines all efficiency benefits of FM, Krylov-based 
and preconditioning components while retaining the robustness and accuracy
of the underlying variational approach. 

We think that all the desirable points, shown in Section \ref{PCG}, including especially the
flexibility in handling non-trivial domains are met by the proposed 
method. It can be clearly recognised that the proposed new integration scheme
generates the most accurate reconstructions independently from the underlying conditions. 
The computational costs are very low and in most cases faster 
than the recent Sylvester method of Harker and O'Leary. Only the DCT 
is much faster, however in DCT results one observes a loss of quality when the computational domain is not rectangular. 

Therefore, the FM-PCG integrator is a good choice for applications which 
require an extremely accurate and robust 3D reconstruction at relatively 
low computational costs.

Nonetheless, our integration method remains limited to \emph{smooth} surfaces. Studying the impact of appropriate preconditioning and initialisation on iterative methods which allow depth discontinuities, as for instance~\cite{Durou2009,Queau2015}, constitutes an interesting perspective. We also consider extending our study to \emph{multi-view} normal field integration~\cite{Chang2007a} as an exciting perspective, which would allow the recovery of a full 3D shape, instead of a depth map.



\bibliographystyle{plain}
{\normalsize \bibliography{biblio}}

\begin{thebibliography}{10}

\bibitem{Agrawal}
A.~Agrawal, R.~Raskar, and R.~Chellappa.
\newblock {What is the Range of Surface Reconstructions from a Gradient Field?}
\newblock In {\em Proceedings of the 9\textsuperscript{th} European Conference
  on Computer Vision}, volume 3951 of {\em Lecture Notes in Computer Science},
  pages 578--591, Graz, Austria, 2006.

\bibitem{Badri}
H.~Badri, H.~M. Yahia, and D.~Aboutajdine.
\newblock {Robust Surface Reconstruction via Triple Sparsity}.
\newblock In {\em Proceedings of the IEEE Conference on Computer Vision and
  Pattern Recognition}, pages 2291--2298, Columbus, USA, 2014.

\bibitem{BB15}
M.~B\"ahr and M.~Breu\ss{}.
\newblock {An improved eikonal method for surface normal integration}.
\newblock In {\em Proceedings of the 37\textsuperscript{th} German Conference
  on Pattern Recognition}, volume 9358 of {\em Lecture Notes in Computer
  Science}, pages 274--284, Aachen, Germany, 2015.

\bibitem{Benzi2002}
M.~Benzi.
\newblock {Preconditioning Techniques for Large Linear Systems: A Survey}.
\newblock {\em Journal of Computational Physics}, 182(2):418--477, 2002.

\bibitem{BQBD16}
M.~Breu\ss{}, Y.~Qu\'eau, M.~B\"ahr, and J.-D. Durou.
\newblock {Highly efficient surface normal integration}.
\newblock In {\em Proceedings of the Algoritmy Conference on Scientific
  Computing}, pages 204--213, Podbanske, Slovakia, 2016.

\bibitem{CCF14}
S.~Cacace, E.~Cristiani, and M.~Falcone.
\newblock {Can local single-pass methods solve any stationary
  Hamilton-Jacobi-Bellman equation?}
\newblock {\em SIAM Journal on Scientific Computing}, 36(2):A570--A587, 2014.

\bibitem{Chang2007a}
J.~Y. Chang, K.~M. Lee, and S.~U. Lee.
\newblock {Multiview Normal Field Integration Using Level Set Methods}.
\newblock In {\em Proceedings of the IEEE Conference on Computer Vision and
  Pattern Recognition, Workshop on Beyond Multiview Geometry: Robust Estimation
  and Organization of Shapes from Multiple Cues}, Minneapolis, USA, 2007.

\bibitem{Coleman1982}
E.~N. Coleman and R.~Jain.
\newblock {Obtaining 3-dimensional shape of textured and specular surfaces
  using four-source photometry}.
\newblock {\em Computer Graphics and Image Processing}, 18(4):309--328, 1982.

\bibitem{Du}
Z.~Du, A.~Robles-Kelly, and F.~Lu.
\newblock {Robust Surface Reconstruction from Gradient Field Using the L1
  Norm}.
\newblock In {\em Proceedings of the 9\textsuperscript{th} Biennial Conference
  of the Australian Pattern Recognition Society on Digital Image Computing
  Techniques and Applications}, pages 203--209, Glenelg, Australia, 2007.

\bibitem{Durou2009}
J.-D. Durou, J.-F. Aujol, and F.~Courteille.
\newblock {Integrating the normal field of a surface in the presence of
  discontinuities}.
\newblock In {\em Proceedings of the 7\textsuperscript{th} International
  Workshop on Energy Minimization Methods in Computer Vision and Pattern
  Recognition}, volume 5681 of {\em Lecture Notes in Computer Science}, pages
  261--273, Bonn, Germany, 2009.

\bibitem{Durou2}
J.-D. Durou and F.~Courteille.
\newblock {Integration of a normal field without boundary condition}.
\newblock In {\em Proceedings of the 1\textsuperscript{st} International
  Workshop on Photometric Analysis For Computer Vision}, Rio de Janeiro,
  Brazil, 2007.

\bibitem{Esteban2008}
C.~H. Esteban, G.~Vogiatzis, and R.~Cipolla.
\newblock {Multi-view photometric stereo}.
\newblock {\em IEEE Transactions on Pattern Analysis and Machine Intelligence},
  30(3):548--554, 2008.

\bibitem{Frankot1988}
R.~T. Frankot and R.~Chellappa.
\newblock {A method for enforcing integrability in shape from shading
  algorithms}.
\newblock {\em IEEE Transactions on Pattern Analysis and Machine Intelligence},
  10(4):439--451, 1988.

\bibitem{Galliani2012}
S.~Galliani, M.~Breu\ss{}, and Y.-C. Ju.
\newblock {Fast and Robust Surface Normal Integration by a Discrete Eikonal
  Equation}.
\newblock In {\em Proceedings of the 23\textsuperscript{rd} British Machine
  Vision Conference}, Guildford, UK, 2012.

\bibitem{Golub96}
G.~H. Golub and C.~F. Van~Loan.
\newblock {\em {Matrix Computation}}.
\newblock Johns Hopkins, 3\textsuperscript{rd} Edition, 1996.

\bibitem{Haque2014}
S.~M. Haque, A.~Chatterjee, and V.~M. Govindu.
\newblock {High Quality Photometric Reconstruction Using a Depth Camera}.
\newblock In {\em Proceedings of the IEEE Conference on Computer Vision and
  Pattern Recognition}, Columbus, USA, 2014.

\bibitem{Harker}
M.~Harker and P.~O'Leary.
\newblock {Least squares surface reconstruction from measured gradient fields}.
\newblock In {\em Proceedings of the IEEE Conference on Computer Vision and
  Pattern Recognition}, Anchorage, USA, 2008.

\bibitem{Harker2015}
M.~Harker and P.~O'Leary.
\newblock {Regularized Reconstruction of a Surface from its Measured Gradient
  Field}.
\newblock {\em Journal of Mathematical Imaging and Vision}, 51(1):46--70, 2015.

\bibitem{Helmsen}
J.~J. Helmsen, E.~G. Puckett, P.~Colella, and M.~Dorr.
\newblock {Two new methods for simulating photolithography development in 3D}.
\newblock In {\em Optical Microlithography IX}, volume 2726 of {\em Proceedings
  of SPIE}, pages 253--261, Santa Clara, USA, 1996.

\bibitem{Hestenes1952}
M.~R. Hestenes and E.~Stiefel.
\newblock {Methods of Conjugate Gradients for Solving Linear Systems}.
\newblock {\em Journal of Research of the National Bureau of Standards},
  6(49):46--70, 1952.

\bibitem{Ho}
J.~Ho, J.~Lim, M.~H. Yang, and D.~Kriegmann.
\newblock {Integrating Surface Normal Vectors using Fast Marching Method}.
\newblock In {\em Proceedings of the 9\textsuperscript{th} European Conference
  on Computer Vision}, volume 3953 of {\em Lecture Notes in Computer Science},
  pages 239--250, Graz, Austria, 2006.

\bibitem{Horn}
B.~K.~P. Horn and M.~J. Brooks.
\newblock {The variational approach to shape from shading}.
\newblock {\em Computer Vision, Graphics and Image Processing}, 33(2):174--208,
  1986.

\bibitem{Kaasschieter1988}
E.~F. Kaasschieter.
\newblock {Preconditioned conjugate gradients for solving singular systems}.
\newblock {\em Journal of Computational and Applied Mathematics},
  24(12):265--275, 1988.

\bibitem{Karacali}
B.~Kara\c{c}ali and W.~Snyder.
\newblock {Noise Reduction in Surface Reconstruction from a Given Gradient
  Field}.
\newblock {\em International Journal on Computer Vision}, 60(1):25--44, 2004.

\bibitem{Kershaw1978}
D.~S. Kershaw.
\newblock {The incomplete Cholesky-conjugate gradient method for the iterative
  solution of systems of linear equations}.
\newblock {\em Journal of Computational Physics}, 26(1):43--65, 1978.

\bibitem{kimmelsethian01}
R.~Kimmel and J.~A. Sethian.
\newblock {Optimal algorithm for shape from shading and path planning}.
\newblock {\em Journal of Mathematical Imaging and Vision}, 14:237--244, 2001.

\bibitem{Klette1996}
R.~Klette and K.~Schl\"uns.
\newblock {Height data from gradient fields}.
\newblock In {\em Proceedings of the Machine Vision Applications,
  Architectures, and Systems Integration}, volume 2908 of {\em Proceedings of
  the International Society for Optical Engineering}, pages 204--215, Boston,
  USA, 1996.

\bibitem{Kovesi2005}
P.~Kovesi.
\newblock {Shapelets correlated with surface normals produce surfaces}.
\newblock In {\em Proceedings of the 10\textsuperscript{th} IEEE International
  Conference on Computer Vision}, volume~2, pages 994--1001, Beijing, China,
  2005.

\bibitem{Manteuffel}
T.~A. Manteuffel.
\newblock {An Incomplete Factorization Technique for Positive Definite Linear
  Systems}.
\newblock {\em Mathematics of Computation}, 34(150):473--497, 1980.

\bibitem{Meijerink1977}
J.~A. Meijerink and H.~A. van~der Vorst.
\newblock {An Iterative Solution Method for Linear Systems of which the
  Coefficient Matrix is a Symmetric $M$-Matrix}.
\newblock {\em Mathematics of Computation}, 31(137):148--162, 1977.

\bibitem{Meister2}
A.~Meister.
\newblock {Comparison of Different Krylov Subspace Methods Embedded in an
  Implicit Finite Volume Scheme for the Computation of Viscous and Inviscid
  Flow Fields on Unstructured Grids}.
\newblock {\em Journal of Computational Physics}, 140:311--345, 1998.

\bibitem{Meister14}
A.~Meister.
\newblock {\em {Numerik linearer Gleichungssysteme -- Eine Einf\"uhrung in
  moderne Verfahren}}.
\newblock Springer Spektrum, 2014.

\bibitem{Meurant1999}
G.~Meurant.
\newblock {\em {Computer Solution of Large Linear Systems}}.
\newblock Elsevier Science, First Edition, 1999.

\bibitem{Meurant06}
G.~Meurant.
\newblock {\em {The Lanczos and Conjugate Gradient algorithms, from theory to
  finite precision computations}}.
\newblock SIAM, 2006.

\bibitem{Noakes2003a}
L.~Noakes and R.~Kozera.
\newblock {Nonlinearities and Noise Reduction in 3-Source Photometric Stereo}.
\newblock {\em Journal of Mathematical Imaging and Vision}, 18(2):119--127,
  2003.

\bibitem{Perez2003}
P.~P{\'e}rez, M.~Gangnet, and A.~Blake.
\newblock {Poisson image editing}.
\newblock {\em ACM Transactions on Graphics}, 22(3):313--318, 2003.

\bibitem{Queau2015}
Y.~Qu\'eau and J.-D. Durou.
\newblock {Edge-Preserving Integration of a Normal Field: Weighted Least
  Squares, TV and L1 Approaches}.
\newblock In {\em Proceedings of the 5\textsuperscript{th} International
  Conference on Scale Space and Variational Methods in Computer Vision}, volume
  9087 of {\em Lecture Notes in Computer Science}, pages 576--588, L\`ege
  Cap-Ferret, France, 2015.

\bibitem{Reddy}
D.~Reddy, A.~K. Agrawal, and R.~Chellappa.
\newblock {Enforcing integrability by error correction using L1-minimization}.
\newblock In {\em Proceedings of the IEEE Conference on Computer Vision and
  Pattern Recognition}, pages 2350--2357, Miami, USA, 2009.

\bibitem{Robles-Kelly}
A.~Robles-Kelly and E.~R. Hancock.
\newblock {A Graph-spectral Method for Surface Height Recovery}.
\newblock {\em Pattern Recognition}, 38(8):1167--1186, 2005.

\bibitem{Saad}
Y.~Saad.
\newblock {\em {Iterative Methods For Sparse Linear Systems}}.
\newblock Society for Industrial and Applied Mathematics, Second Edition, 2003.

\bibitem{Sethian2}
J.~A. Sethian.
\newblock {A fast marching level set method for monotonically advancing
  fronts}.
\newblock {\em Proceedings of the National Academy of Sciences of the United
  States of America}, 93(4):1591--1595, 1996.

\bibitem{Sethian}
J.~A. Sethian.
\newblock {\em {Level Set Methods and Fast Marching Methods}}.
\newblock Cambridge University Press, Second Edition, 1999.

\bibitem{Simchony1990}
T.~Simchony, R.~Chellappa, and M.~Shao.
\newblock {Direct analytical methods for solving Poisson equations in computer
  vision problems}.
\newblock {\em IEEE Transactions on Pattern Analysis and Machine Intelligence},
  12(5):435--446, 1990.

\bibitem{Smith2000}
M.~L. Smith and R.~J. Stamp.
\newblock {Automated inspection of textured ceramic tiles}.
\newblock {\em Computers in Industry}, 43(1):73--82, 2000.

\bibitem{Tang2007}
J.~M. Tang and C.~Vuik.
\newblock {Acceleration of Preconditioned Krylov Solvers for Bubbly Flow
  Problems}.
\newblock In {\em Proceedings of the 7\textsuperscript{th} International
  Conference on Computational Science}, volume~1, pages 874--881, Beijing,
  China, 2007.

\bibitem{Tsitsiklis}
J.~N. Tsitsiklis.
\newblock {Efficient algorithms for globally optimal trajectories}.
\newblock {\em IEEE Transactions on Automatic Control}, 40(9):1528--1538, 1995.

\bibitem{Wang2004}
Z.~Wang, A.~C. Bovik, H.~R. Sheikh, and E.~P. Simoncelli.
\newblock {Image Quality Assessment: From Error Visibility to Structural
  Similarity}.
\newblock {\em IEEE Transactions on Image Processing}, 13(4):600--612, 2004.

\bibitem{WeiKlette2001}
T.~Wei and R.~Klette.
\newblock {A Wavelet-Based Algorithm for Height from Gradients}.
\newblock In {\em Proceedings of the International Workshop on Robot Vision},
  volume 1998 of {\em Lecture Notes in Computer Science}, pages 84--90,
  Auckland, New Zealand, 2001.

\bibitem{WeiKlette2003}
T.~Wei and R.~Klette.
\newblock {Depth Recovery from Noisy Gradient Vector Fields Using
  Regularization}.
\newblock In {\em Proceedings of the 10\textsuperscript{th} International
  Conference on Computer Analysis of Images and Patterns}, volume 2756 of {\em
  Lecture Notes in Computer Science}, pages 116--123, Groningen, The
  Netherlands, 2003.

\bibitem{Woodham1980}
R.~J. Woodham.
\newblock {Photometric Method For Determining Surface Orientation From Multiple
  Images}.
\newblock {\em Optical Engineering}, 19(1):134--144, 1980.

\bibitem{Wu}
Z.~Wu and L.~Li.
\newblock {A line-integration based method for depth recovery from surface
  normals}.
\newblock {\em Computer Vision, Graphics and Image Processing}, 43(1):53--66,
  1988.

\bibitem{Yatziv2006}
{Yatziv, L. and Bartesaghi, A. and Sapiro, G.}
\newblock {O(N) implementation of the fast marching algorithm}.
\newblock {\em Journal of Computational Physics}, 212(2):393--399, 2006.

\bibitem{zafeiriou2013}
S.~Zafeiriou, G.~A. Atkinson, M.~F. Hansen, W.~A.~P. Smith, V.~Argyriou,
  M.~Petrou, M.~L. Smith, and L.~N. Smith.
\newblock {Face recognition and verification using photometric stereo: The
  photoface database and a comprehensive evaluation}.
\newblock {\em IEEE Transactions on Information Forensics and Security},
  8(1):121--135, 2013.

\bibitem{ZBVBWRS08}
H.~Zimmer, A.~Bruhn, L.~Valgaerts, M.~Breu\ss{}, J.~Weickert, B.~Rosenhahn, and
  H.-P. Seidel.
\newblock {PDE-based Anisotropic Disparity-driven Stereo Vision}.
\newblock In {\em Proceddings of the 13\textsuperscript{th} International Fall
  Workshop Vision, Modeling, and Visualization}, pages 263--272, Konstanz,
  Germany, 2008.

\end{thebibliography}

%
%

\end{document}